\documentclass[12pt,fleqn,emtex]{article}

\usepackage{graphicx}
\usepackage{framed}
\usepackage{times}
\usepackage{mathptmx}

\usepackage[german, american]{babel}  %Deutsch

\usepackage{amsfonts}       %
\usepackage{amsmath} %[leqno]=Gleichungen werden links numeriert;
\usepackage{amssymb}        %
\usepackage{amstext}        %
\usepackage{amsthm}         %theorem-umgebungen
\usepackage[matrix,curve]{xypic}

\usepackage{epsfig}
\psfigdriver{dvips}
\usepackage{latexsym}
\usepackage[matrix,curve]{xypic}
\usepackage{rotating}
\rotdriver{dvips}

\newcommand{\Z}{\mathbb{Z}}
\newcommand{\Q}{\mathbb{Q}}
\newcommand{\R}{\mathbb{R}}
\newcommand{\C}{\mathbb{C}}

\newcommand{\Gm}{{{\mathbb{G}}_m}}

\newtheorem{Lemma}{Lemma}
\newtheorem{Theorem}{Theorem}
\newtheorem{Corollary}{Corollary}
\newtheorem{Proposition}{Proposition}

\DeclareMathSizes{12}{11}{8}{5} 
\usepackage{xypic,dsfont}

\begin{document}
%\section{}
%\section{}

\centerline{\Large\bf Why certain Tannaka groups attached to abelian}

\smallskip

\centerline{\Large\bf  varieties are almost connected} 

\bigskip\noindent
\centerline{ R. Weissauer }

\bigskip\noindent
 
\medskip 

{\small {\bf Abstract}: Attached to a smooth projective algebraic variety $Y$ over $\mathbb C$ or a finite field we define reductive algebraic groups $G(Y)$ over certain algebraically closed coefficient fields $\Lambda$.  We show that the groups $G(Y)$ are Zariski connected and explain how their structure is related to Brill-Noether theory.}

\section{Introduction} \label{S1}

\medskip
Let $k$ be an algebraically closed field, where we assume that $k$ is either the field $\C$  of complex numbers or the algebraic closure of a finite field $\kappa$. 
For a smooth projective variety $Y$ over $k$ and a $k$-valued point 
$y_0$ of $Y$ let $$f:Y \longrightarrow Alb(Y)\ $$
be the Albanese mapping,  
normalized such that $f(y_0)=0$. 
The Albanese mapping $f$ and also the induced morphisms $f_i: S^i(Y) \to Alb(Y)$ from the symmetric powers $S^i(Y)$ of $Y$ are proper morphisms, so their images $W_i(Y)$ are closed subvarieties of $Alb(Y)$. In  the simplest case of curves $Y$ their study leads to Brill-Noether theory.
 For curves $Y$ these images ($i$ copies) $$W_i(Y)\ =\ f(Y)+ \cdots +f(Y) \ $$ allow to recover many interesting
properties of the curve $Y$ and its line bundles. 
For higher dimensional varieties $Y$ the extended Albanese mappings $f_i$ are not understood and there may be no analog of Brill-Noether theory, in particular if $i \cdot \dim(Y)$ rapidly exeeds $\dim(Alb(Y))$. Already for surfaces
these questions are mostly unexplored.

\medskip
We propose a more sophisticated approach 
to replace or extend the \lq{classical}\rq\ point of view of Brill-Noether theory with its focus on the images $W_i(Y)$: 
Let us first consider the case of curves $Y$.
For a smooth projective curve $Y$ the constant sheaf with support in $f(Y)$, shifted as a sheaf complex by degree 1, is a perverse sheaf $K=\delta_Y$ on $X= Alb(Y)$. Instead of considering the iterated sums $W_i(Y)$ of the support subvariety $f(Y)$ in $X$, we consider the
iterated convolution products  of $K$ and its Tannaka dual $K^\vee=(-id_{X})^*(DK)$, where the convolution product of two sheaf complexes $M$, $N$ is defined by $M*N= Ra_*(M \boxtimes N)$ for the group law $a: X\times X\to X$ of $X=Alb(Y)$. Here $DK$ denotes the Verdier dual of $K$. 
The convolution product of two semisimple complexes
is a semisimple complex. The simple perverse summands of the iterated convolution products  considered above clearly define denumerably many perverse sheaves on $Alb(Y)$.  For example, all the intersection cohomology sheaves $IC_{W_i(Y)}$ of the Brill-Noether varieties $W_i(Y)$ of the curve $Y$ appear as simple perverse constituents in this way. Indeed, they are among the simple constituents of the $i$-th convolution product ($i$ copies) $$ \delta_Y \ * \ \cdots \ *\ \delta_Y $$ of $\delta_Y$ for $i=1,..,g=\dim(Alb(Y))$. Moreover, the Brill-Noether stratification of the varieties $W_i(Y)$ by their well known  subvarieties $W_i^r(Y)$ can be recovered from the perverse sheaves  $IC_{W_i(Y)}$ as the supports of the cohomology sheaves ${\cal H}^{2r-i}(IC_{W_i(Y)})$. So the above approach via convolution powers obviously contains the input for the classical Brill-Noether theory, at the cost of introducing infinitely many perverse sheaves that in addition to the $IC_{W_i(Y)}$ appear in the iterated convolution powers, and that do not seem to be needed and also do not seem to be easy to understand.  
This should feel less disturbing, once it is understood that all irreducible perverse sheaves
constructed in this way can be naturally interpreted as the irreducible representations of a single reductive group defined over $\Lambda$, so that the $IC_{W_i(Y)}$ correspond to fundamental representations.

\medskip
To be more precise: In the case of a non-hyperelliptic curve $Y$ of genus $g$, the above group turns out to have the commutator group $G(Y)=Sl(2g-2)$. This latter group has $g-1$ fundamental representations defined by the alternating powers $\Lambda^i(V)$ of the $2g-2$ dimensional standard representation $V=V(Y)$, and via the mentioned correspondence  these fundamental representations
correspond  to the $g-1$ irreducible perverse sheaves $IC_{W_i(Y)}$. Hence
for $i\leq g-1$ we get 
$$  \Lambda^i(V)   \ \ \ \longleftrightarrow  \ \ \ IC_{W_i(Y)} \ \ \ \longleftrightarrow  \ \ \ W_i(Y) \ .$$
The remaining fundamental representations $\Lambda^{2g-2-i}(V)$ for $i=0,...,g-1$ similarly correspond
to $\kappa - W_i(Y)$ for the Riemann constant $\kappa$.
In this way one recovers the essential geometric objects $W_i(Y)$ of Brill-Noether theory from $\delta_Y$, using the representation theory of $G(Y)$ via the tensor structure defined by the convolution product. See [BN] and [KrW3].

\medskip
If $Y$ is a projective smooth variety, as in the curve case 
%choose a perverse sheaf on $Y$, the most natural and certainly most interesting candidate being 
consider the perverse intersection cohomology complex $IC_Y = \Lambda_Y[\dim(Y)]$. By our assumptions on $Y$,
this is just the constant sheaf on $Y$ with coefficients in an algebraically closed coefficient field
$\Lambda$ of characteristic zero, up to a shift as a sheaf complex 
by the dimension of $Y$ that makes $IC_Y$ into a perverse sheaf.  The direct image complex $K=Rf_*(IC_Y)$ of the perverse sheaf $IC_Y$  is a sheaf complex on the Albanese variety $X$. By the decomposition theorem it is a direct sum of semisimple perverse sheaves up to complex shifts. In what follows such a sheaf complex  will be called a semisimple complex.  
Notice, one could allow $Y$ to be a normal projective but not necessarily smooth variety, in which case we still dispose over the intersection cohomology sheaf $IC_Y$ except that the correct notion of an Albanese morphism may not be obvious [VS]. Furthermore, we can extend the construction in the following way: 
For $L=IC_Y$, or more generally some simple perverse sheaf $L$ on $Y$, consider the direct image $K=Rf_*(L)$ on $X$ and study the tensor category defined 
from perverse constituents of convolution powers of $K \oplus K^\vee$ and the unit 
in the same way. Again from these data one can define a reductive super group $G(L,y_0)$ over the coefficient field $\Lambda$
together with a natural $\Lambda$-rational super representation $V(L)$ of $G(L,y_0)$, as in the curve case discussed above. These general  
results rely on the construction of an abelian tensor category  attached to the abelian variety $X$ and   
the formalism of Tannaka duality [KrW].

\medskip
As the  discussion above also shows, in the curve case this convolution power approach may not lead far beyond the classical Brill-Noether approach and hence perhaps looks dispensible.
However, for higher dimensional varieties $Y$ this new point of view could turn out to be essential since the convolution products $K^{*i}$ of the perverse sheaf $K=Rf_*(L)$,
with $i$ tensor factors, still contain important information beyond the critical bound $i> \dim(Alb(Y))/\dim(Y)$, where in contrast the Brill-Noether images $f_i(Y)$ in $Alb(Y)$ usually become  uninteresting.  

\medskip
The groups $G(Y,y_0)$ we introduce and study are reductive quotient groups of some large Tannakian super group 
${\bf G}(X)$. In [KrW] this large pro-algebraic reductive super group ${\bf G}(X)$ has been attached to an abelian variety $X$ over $k$. 
It is defined over certain algebraically closed coefficient 
fields $\Lambda$. For this we assume that the coefficient field is $\Lambda=\overline\Q_l$ or $\C$ if $k=\C$, or $\Lambda=\overline\Q_l$ if $k$ is the algebraic closure of a finite field.
The groups $G(L,y_0)$ are algebraic quotient super groups of this large super group scheme ${\bf G}(X)$.    We briefly review the definitions. Let $D(X)$ denote the bounded derived category $D_c^b(X,\Lambda)$ for $k=\C$, respectively for $char(k)\neq 0$ we define $D(X)$ to be the full subcategory of sheaf complexes in $D_c^b(X,\Lambda)$ which are defined over some finite subfield of $k$. 
Then $D(X)$  is a symmetric monoidal rigid $\Lambda$-linear tensor category [BN]. Its tensor structure is defined by the convolution product $K*L$ of complexes $K$, $L$ on $X$. Let $P(X)$ denote the full subcategory of $D_c^b(X,\Lambda)$ of semisimple perverse sheaves in $D(X)$, and let $D^{ss}(X)$ be the full tensor subcategory  of semisimple complexes in $D(X)$; so every object in $D^{ss}(X)$  is a finite direct sum of objects $P[\nu]$ for simple perverse sheaves
$P\in P(X)$ and complex shifts by $\nu\in \Z$. 
Objects in $D^{ss}(X)$, whose simple constituents $P[\nu]$ have Euler characteristic zero, will be called \lq{negligible}\rq  (see the appendix). By [KrW], the negligible objects
define a tensor ideal $N_{Euler}$ of the tensor category $D^{ss}(X)$. As shown in [KrW], the quotient tensor category $\overline D^{ss}(X)$ of $(D^{ss}(X),*)$ obtained from $D^{ss}(X)$ by
dividing out this tensor ideal $N_{Euler}$ is a semisimple neutral super-Tannakian category  and in particular  is an abelian category. Let $\overline P(X)$ denote the image of $P(X)$ in $\overline D^{ss}(X)$. As a tensor category, the abelian quotient category $\overline D^{ss}(X)$ is equivalent to the semisimple tensor category $sRep_\Lambda({\bf G}(X),\mu)$ of finite dimensional super representations over $\Lambda$
$$ \overline D^{ss}(X) \cong sRep_\Lambda({\bf G}(X),\mu) \ ,$$
of some pro-algebraic affine reductive super group ${\bf G}(X)$, where the representations are to satisfy a certain central constraint condition defined by $\mu$. 

\medskip
Any
object $K$ in $\overline D^{ss}(X)$  generates a tensor subcategory, and by the Tannakian formalism therefore defines an algebraic reductive super group ${\bf G}(K)$ over $\Lambda$. 
Via Tannaka duality the sheaf complex $K$  in $D^{ss}(X)$
also defines a finite dimensional faithful super representation $V(K)$ of the group ${\bf G}(K)$, and 
$$ ({\bf G}(K),V(K)) $$ is an irreducible super representation if $K$ is a simple perverse sheaf on $X$.
Attached to the semisimple complex $K=Rf_*(L)$ for some complex $L$ on $Y$, this defines the affine super group scheme $G(L,y_0)$. In the special case $L=IC_Y$ it will be denoted $G(Y,y_0)$, 
and obviously $G(Y,y_0)$ defines some new invariant of $(Y,y_0)$.

\medskip 
From a certain point of view to be explained now, in fact all the supergroups $G(Y,y_0)$ and also the $G(L,y_0)$   are {\it affine reductive groups} over $\Lambda$.  Indeed, there exists a  canonical surjective  group homomorphism from the super group scheme ${\bf G}(X)$ to a multiplicative group $\Gm$ that we
call the {\it cohomological quotient}:  This name comes from the fact that the irreducible representations of the cohomological quotient $\Gm$ keep track of the different perverse cohomology
groups of a semisimple complex; i.e. $t\in \Gm$ acts on the $i$-th perverse cohomology
${}^p H^i(M)$ of a semisimple complex $M$ via multiplication by $t^i$ and this defines and determines the quotient map. As shown in [KrW] there exists a {\it canonical product decomposition} of super group schemes 
 $${\bf G}(X)\ =\ {\bf G}_{geom}(X) \times  {\Gm} $$ so that the projection onto the second factor
 defines the cohomological quotient. This decomposition corresponds to an isomorphism $$sRep_\Lambda({\bf G}(X),\mu)\ \cong\   Rep_\Lambda({\bf G}_{geom}(X)) \ \otimes_\Lambda \ sRep_\Lambda(\Gm,\mu) $$ of the underlying tensor categories.  The splitting of ${\bf G}(X)$ above comes from the fact that as a tensor category $Rep_\Lambda({\bf G}_{geom}(X))$ is equivalent  to the Tannakian  subcategory $\overline P(X)$ in $\overline D^{ss}(X)$ of  perverse  multipliers as defined in [KrW]. The category $sRep_\Lambda(\Gm,\mu) $ is a semisimple tensor category whose simple objects correspond to the complex shifts $\delta_0[i]$ for $ i \in \Z$  of the skyscraper sheaf $\delta_0$ at zero.  The morphism $\mu: \{ \pm 1\} \to \Gm$ in the sense of [KrW] is the canonical inclusion. So the categorial dimension of $\delta_0[i]$ in this tensor category is $(-1)^i$. In other words,
$sRep_\Lambda(\Gm,\mu) $ is equivalent as a tensor category to the category of finite dimensional
$\mathbb Z$-graded $\Lambda$-vector spaces, with a twisted tensor product using the Koszul rule. For $X=0$ notice $sRep_\Lambda(\Gm,\mu) = {\bf G}(X)$. On the other hand,  ${\bf G}_{geom}(X)$ is a projective limit of affine reductive groups over $\Lambda$ (not a super group) so that
$Rep_\Lambda({\bf G}_{geom}(X))$ is the semisimple tensor category  of finite dimensional algebraic representations of ${\bf G}_{geom}$ over $\Lambda$.  This follows from the main result of [W2] together with the characterization of Tannakian categories in terms of dimensions given by Deligne in [D2]. So we may rephrase all this in down to earth terms as follows:  {\it Objects in 
$\overline D^{ss}(X) \cong sRep_\Lambda({\bf G}(X),\mu)$ are $\mathbb Z$-graded finite dimensional algebraic representations of the reductive pro-algebraic group ${\bf G}_{geom}(X)$}. 

\medskip
This being said, we mention that the above decomposition of ${\bf G}(X)$ is closely related to the following theorem of [KrW], [W].

\begin{Theorem} \label{relVT}
Let $f: X \to Y$ be a homomorphism between abelian varieties over $k$ and let $K$ be a perverse
sheaf in $P(X)$. Then
for \lq{most}\rq\ characters $\chi$ of $\pi_1(X,0)$ the direct image $Rf_*(K_\chi)$ of a character twist $K_\chi$ of $K$ is a perverse sheaf in $P(Y)$. 
\end{Theorem} 

In the formulation of this theorem,
$\chi$ is a character $\chi: \pi_1(X,0) \to \C^*$ of the topological resp. etale fundamental group $\pi_1(X,0)$ of $X$, and
$K_\chi$ denotes the perverse sheaf $K \otimes_\Lambda L_\chi$ where $L_\chi$ is the rank one etale local system on $X$ attached to the character $\chi$ of the fundamental group $\pi_1(X,0)$. For the definition of the notion \lq{most}\rq\ we refer to [W, section 1]; if the abelian variety $Y$ is zero, \lq{most}\rq\ means:  $\chi$ is in the complement of a suitable finite union of translates $\chi_i \cdot K(X_i)$ for nontrivial abelian subvarieties $X_i\subseteq X$, where $K(X_i)$ denotes the group of characters $\chi: \pi_1(X,0)\to \Lambda^*$ whose restriction to the subgroup $\pi_1(X_i,0)$ of $\pi_1(X,0)$ vanishes. For further details see  
[KrW].

\medskip  
Before we describe the main result of this paper, we sketch some of the basic properties of the groups and super groups introduced so far: First,
it is very easy to see that the naive commutator subgroup $$G(L):= [G(L,y_0),G(L,y_0)]$$ of the super group $G(L,y_0)$ is independent of the chosen base point $y_0\in Y$. By the splitting of the geometric quotient, $G(L)$ turns out to be an ordinary reductive algebraic group over $\Lambda$. Thus the reductive group $G(Y):=G(Rf_*(IC_Y))$ is an intrinsic invariant of the algebraic variety $Y$ not depending on the base point $y_0$. The group $G(Y)$ is trivial if 
$Alb(Y)=0$. Notice, in this case one can study
the similar groups for finite ramified coverings of $Y$ by extending the equivariant approach explained in  [BN, section 7.8]. %If $Alb(Y)\neq 0$, then $G(Y,y_0)=0$ implies that the topological Euler-Poincare characteristic $\chi(Y)$ of $Y$ vanishes resp. more generally that $\chi(IC_Y)=0$ holds
%if $Y$ is not smooth. 

%\medskip
%Concerning the functoriality of $G(Y,y_0)$, suppose that
%$$  g: Y_1 \longrightarrow Y_2 $$
%is a morphism between smooth projective varieties $Y_1$ and $Y_2$ that respects the base points
%$y_1,y_2$, i.e. $g(y_1)=y_2$.
%Then we have a commutative diagram
%$$ \xymatrix@+0,5cm{ Y_1 \ar[d]^{f_1}\ar[r]^g &  Y_2 \ar[d]^{f_2}\cr
%Alb(Y_1) \ar[r]^{Alb(g)} &  Alb(Y_2) \cr } $$
%For any simple perverse sheaf $L_1$ or more generally semisimple sheaf complex on $Y_1$,
%the direct image complex $L_2=Rg_*(K_1)$ is a semisimple sheaf complex on $Y_2$.
%So we have the reductive super group $G(L_1,y_1)$ defined by $L_1$ and also the reductive super group
%$G(L_2,y_2)$ defined by $L_2$. These groups are related in the following way.  If $Alb(g)$ is surjective, then $G(L_2,y_2)$ is a closed subgroup of $G(L_1,y_1)$; otherwise it is only a subquotient. In general  there exists a contravariant 2-functorial algebraic super group homomorphism (coming from a 2-functor it is uniquely defined only up to conjugation)
%$$  G(g): G(L_2,y_2) \longrightarrow G(L_1,y_1) \ $$ 
%that naturally relates $G(L_1,y_1)$ and $G(L_2,y_2)$. 

\medskip
As already explained, the groups we are interested in are reductive quotient groups of the large Tannakian super group 
${\bf G}(X)$.  Let $\pi_1^{et}(\hat X,0)$ denote the etale profinite fundamental group of the dual abelian variety $\hat X$ of $X$. In [KrW]
we constructed an epimorphism of super groups $$ \pi_X: {\bf G}(X)\ \twoheadrightarrow \ \pi_1^{et}(\hat X,0)(-1)\ $$ whose
kernel contains the connected component ${\bf G}(X)^0$ of ${\bf G}(X)$, defined as the projective limit of the connected components of all algebraic quotient groups of ${\bf G}(X)$. 
Obviously ${\bf G}(X)^0={\bf G}_{geom}(X)^0$. The main result of this paper is

\begin{Theorem} \label{complexcase} {\it The kernel  
of the group homomorphism $ \pi_X: {\bf G}(X) \twoheadrightarrow \pi_1^{et}(\hat X,0)(-1)$ is the \lq{Zariski connected}\rq\ component  ${\bf G}(X)^0$ of the pro-algebraic supergroup ${\bf G}(X)$}.
\end{Theorem}

Theorem 2 immediately implies that the
groups $G(Y)$ are connected algebraic reductive groups.
Hence their finite dimensional irreducible representations over $\Lambda$ can be described in terms of highest weights. In all known examples the representation attached to
$L=IC_Y$ is  fundamental or another distinguished \lq{small}\rq\ representations
of $G(Y,y_0)$.  This irreducible representation, as well as $G(Y,y_0)$ itself, are interesting invariants of $(Y,y_0)$.  
 For the cases where the Albanese morphism $f$ is a generically finite morphism, it is tempting  to believe that the direct image complex $K=Rf_*(L)$ defines some small and distinguished representation of $G(Y,y_0)$ similar as in the curve case.

\medskip
Since the super groups $G(Y,y_0)$ resp. $G(K,y_0)$ are quotients of ${\bf G}(X)$, theorem 2 implies that their group of connected components is a finite quotient group of the profinite group $\pi_1^{et}(\hat X,0)(-1)$. In section \ref{S2} before lemma \ref{newlemma} we explain how this can be made more explicit.

\medskip This paper is organized as follows: The sections \ref{S1}, \ref{S2} and \ref{Groth} contain general results
on the groups considered. The sections \ref{S4},\ref{Drinfeld},7,9  contain the proof of theorem \ref{complexcase}. In section \ref{S4} its proof  will be reduced to  the statement of theorem \ref{abelian}: Perverse sheaves $K\cong K^\vee$ with the property $K*K \cong d \cdot K$ and without negligible constituents are skyscraper sheaves on $X$
with support  contained in the torsion points of $X(k)$. This statement is shown by induction on the dimension of $X$ (end of section \ref{S4}). The induction step uses reduction mod $p$ in order to deduce the assertion from the finite field case (section \ref{Drinfeld}). Finally we apply a theorem of Drinfeld [Dr2] on the existence of compatible $l$-adic systems for perverse sheaves (section \ref{congruences}) in combination with the Fourier transform (section \ref{Fourier}).
That the underlying group structure of $X$ is abelian then allows  to deduce the induction step from some elementary assertion on the existence of certain torsion points (lemma \ref{Jacobi}). This uses the Fourier transform defined in section \ref{Fourier}. Further results, needed for [W] and [W2] (in particular lemma \ref{14} and corollary \ref{key}), are shown \lq{en passant}\rq\   in the remaining  sections \ref{CEB} and \ref{Plancherel}. This involves the Fourier transform and  arguments based on the Cebotarev density theorem. 
Finally, in section \ref{invertible} we classify the invertible objects of the category $\overline P(X)$ and show that these are the irreducible skyscraper sheaves. The appendix deals with general facts on translation invariant sheaf complexes.

\goodbreak

\section{Homomorphisms between abelian varieties}
\label{S2}

We start with some remarks on descent for perverse sheaves.

\smallskip
\begin{Lemma} \label{descender} For a finite group $H\subset X(k)$ of torsion points consider the etale isogeny $f: X \to X'=X/H$. 
Then for a perverse sheaf $K$ on $X$ the following assertions are equivalent:
\begin{enumerate}
\item $\sigma_x: T^*_x(K)\cong K$ holds for the translations $T_x$ by $x\in H$ with respect to isomorphisms $\sigma_x$ satisfying  cocycle descent conditions in [KW, def.III.15.1]. 
\item $K = f^*(K')$ holds for a perverse sheaf $K'$ on $X'$.
\end{enumerate}
\end{Lemma}

{\it Proof}. The direction 2. $\Rightarrow$ 1. is obvious from $f\circ T_x =f$. 
For the converse first suppose that 
$K$ is an irreducible perverse sheaf on $X$. Then $K = i_*j_{!*} L[dim(U)]$ for a smooth etale sheaf $L$ on $U$, where $j: U\hookrightarrow \overline U$ is an open embedding and 
$i: \overline U \to X$ is an irreducible closed subvariety $\overline U$ of $X$. The image $U'$ of $U$ in $X'$ is an irreducible subvariety $j'\! : U'\hookrightarrow X'$, and by shrinking $U$ we may assume that $f: U \to U'$ is an etale
covering. The sheaf $L$ defines a local system, in other words corresponds to an irreducible finite dimensional representation
$\rho: \pi_1(U,x_0) \to Gl(r,\Lambda)$. 
By the descent conditions
the representation $\rho$ can be extended to a representation $\rho'$ in the sense above. 
Indeed this follows by descend theory for etale sheaves. Notice that $U\times_{U'} U \cong
\bigcup_{g\in H} U$ since $U\to U'$ is a an etale Galois covering with Galois group $H$.
So the conditions 1. are equivalent to the usual descent conditions.
Therefore $L \cong f^*(L')$ holds
for some smooth etale sheaf $L'$ on $U'$. Hence $K \cong f^*(K')$ for
$K' = i'_* j'_{!*} L'$ for the inclusion $i'\!: \overline U' \hookrightarrow X'$. 
\qed

\begin{Corollary} \label{c1} For an isogeny $f: X\to X'$ of abelian varieties over $k$ and
for a simple perverse sheaf $K'$ on $X'$ the pullback $K=f^*(K')$ is semisimple and the translations $T^*_x$, $x\in Kern(f)$ act transitively on the isotopic components of the simple constituents of $K$. 
\end{Corollary}

{\it Proof}. For inseparable isogenies this follows from well known properties of the etale topology.
So, assume that $f$ is separable. Then we can apply lemma \ref{descender} to obtain isomorphisms
$\sigma_x: T^*_x(K)\cong K$ for $x\in Kern(f)$ with descend conditions.
The maximal semisimple subobject $S$ of $K$ is nontrivial and stable under these isomorphisms
$\sigma_x$. By lemma \ref{descender} it is therefore of the form $S=f^*(S')$ for some perverse sheaf $S'$ on $X'$. 
By the perverse adjunction formula $0 \neq Hom(f^*(K'),S) = Hom(K',f_*(S)) = Hom(K',\bigoplus_\chi S'_\chi)$ it follows that the simple perverse sheaf $K'$ injects into some $S'_\chi$, where  $\chi$ is a character of $\pi_1(X',0)$ whose restriction to $\pi_1(X,0)$ is trivial. Hence $K=f^*(K')\hookrightarrow f^*(S'_\chi) \cong f^*(S')=S$. On the other hand $S$ injects into $K$.
Comparing the length, we conclude $S=K$ and therefore $K$ is semisimple. Let $\tilde S$
be the span of all translates $T_x^*(L')$, $x\in Kern(f)$ of an isotopic constituent $L'$ corresponding to some simple constituent $L$ of $K$. Then the same argument as above
implies $\tilde S=K$.  This proves the last assertion.
\qed

\medskip
Recall from [KrW] that any surjective homomorphism $f:X\to Y$  between abelian 
varieties induces an embedding ${\bf G}(f): {\bf G}(Y) \hookrightarrow {\bf G}(X)$
of the supergroups ${\bf G}(X)$ and ${\bf G}(Y)$ attached to $X$ and $Y$, as explained in the introduction. In the first part of this section we consider the case where $f$ is an isogeny and we show in 
lemma \ref{isog} that ${\bf G}(Y)$ is a subgroup of ${\bf G}(X)$ of finite index and describe this explicitely.

\medskip 
An  isogeny $f:X \to Y$ of abelian varieties over $k$ with kernel $F$
$$ \xymatrix{ 0 \ar[r] & F \ar[r] &  X \ar[r]^f &  Y \ar[r] &  0 } $$
factorizes into a composite of an inseparable isogeny and a separable finite etale
isogeny. The  degree of the latter will be denoted  $deg_s(f)$. Attached to
$f$ there is an embedding of groups ${\bf G}(f) : {\bf G}(Y)\hookrightarrow {\bf G}(X)$,
as defined in [KrW, lemma 13.7].  

\begin{Lemma}\label{permanent}
For an isogeny $f$ the direct image $f_*=f_!$ and the pullback $f^*=f^!$ are adjoint functors that
preserve the category of perverse sheaves
$$  f_*: Perv(X,\Lambda) \to Perv(Y,\Lambda) \quad , \quad
 f^*: Perv(Y,\Lambda) \to Perv(X,\Lambda) \ .$$
\end{Lemma} 
 
This is clear if $f$ is completely inseparable. Hence it suffices
to consider the case where $f$ is finite and etale, where this follows from [BBD] cor. 2.2.6.
Furthermore, $f_*$ maps semisimple perverse sheaves to semisimple perverse sheaves
by the decomposition theorem. $f^*$ preserves semisimplicity  by corollary \ref{c1}. 
Hence $f_*: P(X)\to P(Y)$ and $f^*: P(Y)\to P(X)$. \qed

\medskip
If we now use lemma \ref{7} from the appendix, we see that the functors $f_*$ and $f^*$ also preserve  the subcategories $N_{Euler}$ of negligible perverse sheaves.
Hence we obtain induced functors $\overline f_*: \overline D^{ss}(X) \to \overline D^{ss}(Y)$
and $\overline f^*: \overline D^{ss}(Y) \to \overline D^{ss}(X)$. 
Recall, by Tannaka duality $\overline D^{ss}(X) \cong sRep({\bf G}(X),\mu)$ and similarly for $Y$. Hence we obtain functors
$$R:=\overline f_* \ \ \mbox{  and } \ \  I:=\overline f^*$$ 
$$ R: sRep_\Lambda({\bf G}(X),\mu) \to sRep_\Lambda({\bf G}(Y),\mu) \quad , \quad I: sRep_\Lambda({\bf G}(Y),\mu) \to sRep_\Lambda({\bf G}(X),\mu)
\ .$$
As shown above, $f$ induces functors $\overline f_*: \overline P(X)\to \overline P(Y)$
resp. $\overline f_*: \overline P(Y)\to \overline P(X)$. Hence, the functors $R$ and $I$ restrict to functors between
the subcategories 
$$ R: Rep_\Lambda({\bf G}_{geom}(X),\mu) \to Rep_\Lambda({\bf G}_{geom}(Y),\mu) $$ $$ I: Rep_\Lambda({\bf G}_{geom}(Y),\mu) \to Rep_\Lambda({\bf G}_{geom}(X),\mu)
\ .$$
Obviously $f_*(K[1])=f_*(K)[1]$ and $f_*(\delta_0)=\delta_0$, and in particular $f_*(\delta_0[i])=\delta_0[i]$ holds.
Hence $R$ induces the identity functor on the tensor subcategory $sRep_\Lambda(\mathbb G_m,\mu)$ 
corresponding to the cohomological quotient groups $\mathbb G_m$. This implies the next lemma.

\begin{Lemma} \label{canonic} For isogenies $f$ the inclusion
${\bf G}(f): {\bf G}(Y) \hookrightarrow {\bf G}(X)$ respects the canonical product decompositions
${\bf G}(X) = {\bf G}_{geom}(X)\times {\mathbb G}_m$ and ${\bf G}(Y) = {\bf G}_{geom}(Y)\times \mathbb G_m$, i.e. is of the form  $${\bf G}(f)={\bf G}_{geom}(f) \times id_{{\mathbb G}_m}\ .$$
\end{Lemma}

\begin{Lemma}
The functors $R$ and $I$ are the restriction an the induction functor
with respect an inclusion of supergroups ${\bf G}(f): {\bf G}(Y) \hookrightarrow {\bf G}(X)$ 
of finite index $deg_s(f)$.
\end{Lemma}

\medskip
{\it Proof}. The adjunction formula $Hom_{D(X)}(f^*(K),L)= Hom_{D(Y)}(K,f_*(L))$
for the adjoint functors $f^*$ and $f_*$ implies $Hom_{\overline D^{ss}}(X)(f^*(K),L)= Hom_{\overline D^{ss}}(Y)(K,f_*(L))$, hence by Tannaka duality we obtain
$$ Hom_{{\bf G}(X)}(I(K),L) = Hom_{{\bf G}(Y)}(K,R(L)) \ .$$
Since $R$ is given by the restriction of representations from ${\bf G}(X)$ to the subgroup ${\bf G}(Y)$ via ${\bf G}(f)$,
the functor $I$ corresponds to induction of representations. 

\medskip
Let $\pi_1(X,0)$  denote the topological fundamental group
for $k=\C$ respectively the etale fundamental group for fields $k$ of characteristic $>0$.
The group $F(k)$ can be identified with the cokernel of the homomorphism
$\pi_1(f): \pi_1(X,0) \to \pi_1(Y,0)$ and has order $deg_s(f)$. The Pontryagin dual $F^*$ of $F$
can be identified with the group of characters $\chi: \pi_1(Y,0)\to \Lambda^*$
that vanish on the image of $\pi_1(X,0)$. For a \lq{continuous}\rq\ character $\chi: \pi_1(X,0)\to \Lambda^*$ and the corresponding local system $L_\chi$ on $X$, for  $K\in D(X)$ let $K_\chi=K\otimes_\Lambda L_\chi$  denote the twisted complex. Since $R(I(K)) =  \overline f_*(\overline f^*(K))$ is represented by
$f_*(f^*(K)) = f_*(\Lambda_X) \otimes_\Lambda K$, therefore $Rf_*(\Lambda_X)= \bigoplus_{\chi\in F^*} K_\chi$ and this implies
$$   R(I(K)) = \bigoplus_{\chi\in F^*}  K_\chi \ .$$
Notice, for the categorial dimensions $\dim(K_\chi)\!=\!\dim(K)$ holds (see [KrW, cor.10]). Since $R=\overline f_*$  is a tensor functor, the functor $R$ preserves categorial dimensions. Hence $  \dim(I(K))\! =\! \# F(k) \cdot \dim(K)$ or $$  \dim(I(K)) = deg_s(f) \cdot \dim(K) \ .$$ 
This implies that ${\bf G}(Y)$ is a super subgroup of ${\bf G}(X)$ of index $deg_s(f)$.
\qed

\begin{Lemma} \label{isog} The inclusion ${\bf G}(f)\!: {\bf G}(Y) \to {\bf G}(X)$ identifies ${\bf G}(Y)$ with a normal super subgroup of ${\bf G}(X)$
of index $[{\bf G}(X):{\bf G}(Y)]=deg_s(f)$ and there exists a commutative diagram
$$    \xymatrix{0\ar[r] & Kern(\pi_Y) \ar[r]\ar[d]_\cong &  {\bf G}(Y) \ar[r]^-{\pi_Y} \ar@{^{(}->}[d]^{{\bf G}(f)} &  \pi_1^{et}(\hat Y,0)(-1)  \ar@{^{(}->}[d]^{\pi_1^{et}(\hat f)} \ar[r] & 0 \ \cr
0\ar[r]  & Kern(\pi_X) \ar[r] &  {\bf G}(X) \ar[r]^-{\pi_X} &  \pi_1^{et}(\hat X,0)(-1) \ar[r] & 0  .\cr } $$
\end{Lemma}

\medskip
{\it Proof}. For the proof the pro-algebraic supergroups ${\bf G}(X)$ and ${\bf G}(Y)$ can be replaced by 
the pro-algebraic groups ${\bf G}_{geom}(X)$ and ${\bf G}_{geom}(Y)$, and these can be replaced
by algebraic quotient groups, denoted $G(X)$ and $G(Y)$ in this section.

\medskip
Let $L$ denote an irreducible object in $Rep_\Lambda(G(X))$. 
Then $Hom_{G(X)}(I({\bf 1}),L)\!=\! 0$ implies $Hom_{G(Y)}({\bf 1},R(L))\!=\! 0$. Since
$I({\bf 1})$ contains only  finitely many irreducible subrepresentations $L$,
there are only finitely many isomorphism classes of irreducible representations $L$ of $G(X)$ whose restriction $R(L)$
to the subgroup $G(Y)$ contains the trivial representation ${\bf 1}$.  
Using induction and restriction between $G(Y)$ and  $G^0(Y)$, the analogous assertion  also follows for the Zariski connected components $G^0(Y)$ and $G^0(X)$. Indeed, if
for an irreducible representations $L^0$ of $G^0(X)$ its restriction to  the subgroup $G^0(Y)$ contains the trivial representation, then the  restriction of $L=Ind_{G^0(X)}^{G(X)}(L^0)$
to $G(Y)$ contains the trivial representation. 
This being said, let $\Sigma$ denote the finite set of highest
weights of the finitely many irreducible representations $L^0$ of $G^0(X)$ whose restriction
to $G^0(Y)$ contains the trivial representation. 
Since $char(\Lambda)=0$, now classical invariant theory implies $G^0(X) = G^0(Y)$. Indeed, the coordinate ring $\Lambda[G^0(X)/G^0(Y)] = \Lambda[G^0(X)]^{G^0(Y)}$ of
the quotient $G^0(X)/G^0(Y)$ is affine, since $G^0(Y)$ is a reductive algebraic group [M]. For irreducible $G^0(X)$-subcomodules $L$ in the coordinate ring $\Lambda[G^0(X)]$ then
$L^{G^0(Y)} =0$ holds for $\alpha\notin \Sigma$ if $L$ defines an irreducible representation of $G^0(X)$ of highest weight $\alpha$. This implies $dim_{\Lambda}(\Lambda[G^0(X)/G^0(Y)]) < \infty$ for the the coordinate ring $\Lambda[G^0(X)]$, since  any irreducible comodule $L$ occurs with finite multiplicity in $\Lambda[G^0(X)]$ up to isomorphism. Since $G^0(X)/G^0(Y)$ is Zariski connected, therefore as an algebra
$\Lambda[G^0(X)/G^0(Y)]$  is obtained by adjoining nilpotent elements to $\Lambda$. So $\Lambda[G^0(X)/G^0(Y)]=\Lambda$ by [GIT, thm. 1.1]; see also page 5 of loc. cit. and this proves 
$$  G^0(Y)= G^0(X) \ .$$
Hence the reductive group $G(Y)$  has finite index in the reductive group $G(X)$. 
Therefore $\dim(I(K)) = [\pi_0(G(X)):\pi_0(G(Y))] \cdot \dim(K)$, and we get
$$   [\pi_0(G(X)):\pi_0(G(Y))] \ = \ deg_s(f) \ .$$
To show that $G(Y)$ is a {\it normal} subgroup of $G(X)$, it now
suffices to show for the abelian quotient groups  $\pi_0(G(X))$ and $\pi_0(G(Y))$ 
$$ [\pi_0(G(X))^{ab}:\pi_0(G(Y))^{ab}] \geq deg_s(f)  \ .$$
For this we are allowed to make the assumption that $G(X)$, as a quotient of ${\bf G}(X)$, is large enough
in the following sense: $G(X)$ is the Tannaka group of a tensor subcategory ${\cal T}(K)$ generated by some object $K$ of $(\overline D^{ss}(X),*)$. We can suppose that the images of $\delta_x, x\in F(k)$
are contained in ${\cal T}(K)$ by enlarging $K$, without restriction of generality.
Observe that  any $x\in F(k)$ defines a one dimensional character
of $G(X)$, whose restriction to $G(Y)$ defines
the trivial representation of $G(Y)$ by $Rf_*(\delta_x)=\delta_0$. Thus $ [\pi_0(G(X))^{ab}:\pi_0(G(Y))^{ab}] \geq
\# F(k)$. Since $\# F(k)=deg_s(f)$ and hence $ [\pi_0(G(X))^{ab}:\pi_0(G(Y))^{ab}] \leq 
\# F(k)$,  this immediately proves our claim by passing to the projective limits ${\bf G}(X)$ resp. ${\bf G}(Y)$. \qed

\medskip
In the last steps of the proof we constructed an isomorphism between $F(k)$ and $Hom({\bf G}(X)/{\bf G}(Y), \Lambda^*)$. Notice that by lemma \ref{canonic} we have $${\bf G}(X)/{\bf G}(Y) = {\bf G}_{geom}(X)/{\bf G}_{geom}(Y)\ .$$

\medskip {\it Conjugation action}. For a given isogeny $f:X\to Y$ the group ${\bf G}(X)$ acts on its normal subgroup ${\bf G}(Y)$ by conjugation. For a representation $P$ of ${\bf G}(Y)$ and $g\in {\bf G}(X)$, let $P^g$ denote the representation of ${\bf G}(Y)$ twisted by conjugation with $g$. Up to isomorphism $P^g$ depends only on the coset $g\in {\bf G}(X)/{\bf G}(Y)$.  Since by Mackey's lemma the conjugates $P^g$ of an irreducible representation $P$ of ${\bf G}(Y)$  are the constituents of the module $R(I(P))$, they are isomorphic to simple modules $P_\chi$ obtained from $P$ by a character twist as follows from $$\bigoplus_{g\in {\bf G}(X)/{\bf G}(Y)} P^g \ \cong \ R(I(K)) \ \cong \ \bigoplus_{\chi\in F^*} \ P_\chi\ .$$
The elements $\chi$ in the sum on the right side correspond to the characters in the Pontryagin dual $F^*$ of $F$.
The elements in ${\bf G}(X)/{\bf G}(Y)$ are described by $ \pi_1^{et}(\hat X,0)(-1)/\pi_1^{et}(\hat Y,0)(-1) \cong \hat F(-1)$ for the Cartier dual $\hat F$ of $F$, and $\hat F(-1)$ can be identified
with the Pontryagin dual $F^*$. This is obvious for $char(k)=0$, and left as an exercise for $char(k)>0$.
Thus we have shown for the isogeny $f: X\to Y$ with kernel $F$ the following

\begin{Lemma}
The conjugates of the super representation corresponding to the perverse sheaf $P$ on $Y$
under elements $g \in {\bf G}(X)/{\bf G}(Y)$ are the super represenations of ${\bf G}(Y)$ corresponding to the twisted perverse $P_\chi$ sheaves for
$\chi \in F(k)^*$.
\end{Lemma}

\medskip
{\it Translations}. Let $T_x: X\to X$ denote translation with respect to $x\in X(k)$.
For a simple perverse perverse sheaf $K$ on $X$ we denote its stabilizer by $$ H=\{ x\in X \ \vert \ T_x^*(K)\cong K\} \ .$$

\begin{Lemma} \label{finit} If the simple perverse sheaf $K$ is not negligible, $H$ is a finite group.
\end{Lemma}
 
{\it Proof}. The Tannaka dual $K^\vee = (-id_X)^*D(K)$ (for the Verdier dual $D$) 
also is a simple perverse sheaf on $X$, and $K^\vee * K$ is a semisimple  sheaf complex on $X$.
Suppose $K$ is not negligible. Then $K^\vee * K$ contains the skyscraper perverse sheaf $\delta_0$ (unit object) as a direct summand
with multiplicity one (Schur's lemma). Furthermore, for a simple perverse sheaf $L$ the convolution
product $K^\vee*L$ contains $\delta_0$ as a direct summand if and only if $K$ and $L$ are isomorphic perverse sheaves. Since $T_x^*(K^\vee*L) \cong K^\vee * T^*_x(L)$ for $x\in X(k)$, under the assumption that
$K$ is simple and not negligible this implies for all $x\in H$  that $T_x^*(\delta_0)=\delta_{-x}$
is a perverse constituent of $K^\vee*K$. But $K^\vee * K$ has finite length. Therefore
the stabilizer $H$ of a simple not negligible perverse sheaf $K$ is a finite subgroup of $X$. \qed

\medskip 
In the situation of lemma \ref{finit}, for the stabilizer $H$ of $K$  consider now the separable isogeny $$\pi: X\to Y=X/H$$ and  the semisimple perverse sheaf $L=\pi_*(K)$.  

\begin{Lemma} For a fixed simple constituent $P$ of $L$ all
other constituents are of the form $P_\chi$ for some $\chi\in H^*$. 
The pullback $\pi^*(P) \cong m \cdot K$ is isotypic for some multiplicity $m$. Furthermore $m^2$
is the cardinality of the group
$$ \Delta = \{ \chi \in H^* \ \vert \ P_\chi \cong P \} \ .$$
\end{Lemma}

\medskip
{\it Proof}. 
$L^\vee * L \cong \pi_*(K^\vee * K)$ contains $\pi_*(\bigoplus_{x\in H} \delta_{-x}) \cong \# H \cdot \delta_0$ 
as a direct summand.  Furthermore $\delta_0$ occurs in $L^\vee *L$ with
precise multiplicity $\# H$. Indeed, any simple summand in $L^\vee *L$ is a summand of $\pi_*(S)$ for some simple perverse sheaf $S$ on $X$, and 
$Hom_{D(Y)}(\delta_0, \pi_*(S))= Hom_{D(X)}(\pi^*(\delta_0),S)$ implies that
$\delta_0$ is contained in $\pi_*(S)$ if and only if $S \cong \delta_{-x}$ holds for some $x\in H$. 
Then $S\cong \delta_{-x}$ and $\delta_{-x}$ occurs in $K^\vee*K$.
But $\delta_{-x}$ occurs in $K^\vee*K$ (and then with multiplicity one) if and only if $\delta_0$ occurs in $K^\vee * T_{-x}^*(K)$, or in other words if and only if $T_{-x}^*(K) \cong K$ holds. Here 
we used Schur's lemma for $K$.
Since $L$ is perverse and semisimple and 
$0\neq Hom_{P(Y)}(P,L) = Hom_{P(Y)}(P,\pi_*(K)) =Hom_{P(X)}(\pi^*(P), K)$,
there exists an a nontrivial morphism from $\pi^*(P)$ to $K$. Since $K$ is a simple perverse sheaf, this is an epimorphism.
Since $\pi^*(P)$ is semisimple, hence $K$ is one of its constituents.
So $$\pi^*(P) \cong  m \cdot K$$ is isotopic by corollary \ref{c1} because all translates $T_x(K)$, $x\in H$ are isomorphic to $K$ by the definition of $H$.  
Then $m \cdot \pi_*(K)\cong \pi_*\pi^*(P) \cong
\bigoplus_{\chi \in H^*} P_\chi$, and we obtain 
$$ m \cdot \pi_*(K) \cong \bigoplus_{\chi\in H^*} P_\chi  .$$
Hence all constituents of $L=\pi_*(K)$ are isomorphic to twists $P_\chi$ of $P$.
On the other hand $End_{P(X)}(m\cdot K) = End_{P(X)}(\pi^*(P)) = Hom_{P(Y)}(P,\pi_*\pi^*(P))
= Hom_{P(Y)}(P, m \cdot \pi_*(K))$ and 
$\dim_\Lambda(End_{P(X)}(m \cdot K))=m^2$
imply $\dim_\Lambda(Hom_{P(Y)}(P,\pi_*K))=m$.  For
$\chi\in H^*$, 
the same holds for $P$ replaced by $P_\chi$. Hence every constituent of $L$ appears with exact multiplicity $m$
and all $P_\chi$, $\chi\in H^*$ occur as constituents.
Therefore $L^\vee * L$ contains $\delta_0$ with multiplicity equal to
$m^2 \cdot \# H^* /\# \Delta$. As shown above, this multiplicity is 
$\# H = \# H^*$.
It follows that
 $\# \Delta = m^2 $. \qed

\bigskip
{\it Group theory}. We will now apply the following well known facts from group theory.
Let $N$ be a normal subgroup of finite index of an algebraic group $G$ over $\Lambda$
and let $\rho$ be an irreducible algebraic representation of $G$ on a finite dimensional $\Lambda$-vector space $V$. Then by
[S, theorem 16], either the restriction $V\vert_N$ of $V$ to $N$ is isotypic,
or there exists a proper subgroup $H$ of $G$ containing $N$ and an irreducible representation
$W$ of $H$ that induces $V$. 
If $V=Ind_H^G(W)$,  then for all characters $\psi: G/H \to \Lambda^*$ by Frobenius reciprocity $Hom_G(V,V\otimes \psi) \cong Hom(W,R(V\otimes\psi)) =  Hom(W,R(V)) \neq 0$ and hence $V$ is isomorphic to $V\otimes \psi$.  So, if $G/N$ is a finite abelian group and
$V \not\cong V \otimes \psi$ holds for all nontrivial characters $\psi: G/N \to \Lambda^*$, then the restriction of $V$ to $N$ must be  isotypic, say
$V\vert_N \cong m\cdot W$ for some irreducible representation $W$ of $N$. Again by Frobenius reciprocity, all the nonisomorphic irreducible representations $V\otimes \psi$, $\psi: G/N \to \Lambda$
are constituents of the induced representation $\tilde V = Ind_N^G(W)$. Since $\tilde V\vert_N = 
[G:N] \cdot W$ and $\bigoplus_\psi  (V\otimes \psi)\vert_N = \#\{ \psi\} \cdot m\cdot W$,
the inclusion $\bigoplus_\psi  (V\otimes \psi) \subseteq \tilde V$ and $\#\{ \psi\} =[G:N]$ imply 
$m=1$ and $\tilde V= \bigoplus_\psi  (V\otimes \psi)$. In particular, $V\vert_N =W$ remains irreducible
as a representation of $N$. In other words, this shows that the restriction of an irreducible representation $V$ to the connected component $N=G^0$ of $G$ remains irreducible if $G/G^0$ is abelian and 
$V \not\cong V\otimes \chi$ holds for all characters $\chi: G/N\to \Lambda^*$. We want to apply these remarks in the following context.

\medskip
For a simple perverse sheaf $P$ on $X$ that is not negligible we explained in the introduction how one can define a reductive group ${\bf G}(P)$  and a faithful irreducible representation $V=V(P)$ of ${\bf G}(P)$ over $\Lambda$. In this context we may apply the above mentioned group theoretic considerations. Indeed, as a consequence of the later theorem
\ref{abelian} we can apply the group theoretic facts above
for the normal subgroup $N={\bf G}^0(P)$ defined by the connected component of ${\bf G}(P)$ since the quotient group is abelian. If  $V \not\cong V\otimes \chi$ holds for all characters $\chi: {\bf G}(P)/{\bf G}^0(P)\to \Lambda^*$, hence the restriction of the representation $V$ to the connected component ${\bf G}^0(P)$ remains irreducible.
If we apply this in the special case $P=IC_Y$ 
for some irreducible subvariety $Y$ of $X$, we obtain

\begin{Lemma} \label{remainirr} For
an irreducible subvariety $Y$ of the abelian variety $X$ let $H(Y)=\{ x\in X(k)\ \vert \ Y+x =Y\}$ be
the stabilizer. If $H(Y)=0$, then the canonical representation $V=V(IC_Y)$ of the reductive algebraic group $G={\bf G}(IC_Y)$  attached to $Y$ restricts to an irreducible
and faithful representation of the Zariski connected component $G^0$ of $G$. Hence the center of $G^0$ is either a finite cyclic group or a torus isomorphic to ${\mathbb G}_m$.    
\end{Lemma}

%\medskip
%{\it Proof of lemma \ref{remainirr}}. The natural map
%$\pi_1(Y,y_0) \to \pi_1(X,0)$ has finite kokernel $Q$ and there exists an isogeny $f: \tilde X \to X$
%with covering group $Q$ dual to $F=kern(f)$ such that the inverse image $\tilde Y= f^{-1}(Y)$ of $Y$ decomposes into a disjoint union $\bigcup_{x\in F} Y'$ of irreducible varieties all isomorphic to $Y$. It easily follows that $G(\tilde Y,0) \cong
%F \times G(Y',0)$. Furthermore $G(Y',0) \cong G(Y,0)$, so we can identify these two super groups
%and their canonical representations $V'$ and $V$.
%By lemma 4 in [W2], for $P=IC_Y$ the condition $P_\chi \cong  P$ implies that $\chi: \pi_1(X) \to \Lambda^*$ is a character that is trivial on the image of $\pi_1(Y,y_0)$. Hence $P'_\psi \cong P'$
%for $P'=IC_{Y'}$ implies $\psi=1$. So the canonical representation
%$V'$ restricts to an irreducible and faithful representation of the connected reductive group ${\bf G}^0(P')$. Hence the same holds for $Y$: {\it The canonical representation
%$V$ of $G(Y)$ restricts to an irreducible and faithful representation of the connected reductive group ${\bf G}^0(Y)$.}
%Therefore the center of ${\bf G}^0(Y)$ is either finite, or isomorphic to a torus $\mathbb G_m$ by Schur's lemma. In the latter case the maximal abelian quotient group ${\bf G}^0(Y)^{ab}$ is a torus
%isomorphic to $\mathbb G_m$ and embeds into ${\bf G}(Y)^{ab}$. \qed

\medskip
We remark that in the situation of the last lemma the center can be ${\mathbb G}_m$ only for $k=\mathbb C$, but not in the case when $k$ is the algebraic closure of a finite field.   This follows from the later proposition \ref{units} resp. [KrW]  prop. 10.1(b). This references furthermore imply that for a given irreducible perverse sheaf $P$  the group ${\bf G}(K)$ is semisimple for a suitable translate $K=T^*_x(P)$ of $P$.

\begin{Lemma} \label{l3}
For an irreducible subvariety $Y$ of $X$
there exists a translate $K$ of $P=IC_Y$ so that the group
${\bf G}(K)$ is semisimple.
\end{Lemma}

\medskip
{\it Perfect quotients}. For an abelian variety $X$ let now ${\cal T}$ denote the Tannakian category
generated by a fixed semisimple complex $K$. 
The subcategory ${\cal T}_0$ of $\cal T$ generated by 
skyscraper sheaves $\delta_x$ in $\cal T$  for torsion points $x\in X$ then is stable
under convolution. Hence it defines a Tannakian subcategory ${\cal T}_0$  of $\cal T$. 
Its Tannaka group can be identified with a finite abelian quotient group $F_0$ of the reductive Tannaka supergroup
${\bf G}(K)$ attached to ${\cal T}$. The invertible  objects $\delta_x\in {\cal T}_0$ define one dimensional characters of ${\bf G}(K)$ that are trivial on the image ${\bf G}^1(K)$ of $Kern(\pi_X: {\bf G}(X)  \to \pi_1^{et}(\hat X,0)(-1))$ in ${\bf G}(K)$. By Pontryagin duality and $\delta_x *\delta_y = \delta_{x+y}$, this  allows to identify
$F_0$ with the corresponding finite group $F$ of torsion points in $X$. 

\medskip
For the isogeny $f: X\to Y=X/F$, the Tannakian category ${\cal T}'$ generated by $K'=Rf_*(K)$
is the image of the tensor category ${\cal T}$ under the tensor functor $f_*=Rf_*$. 
We claim, ${\cal T}'$ does not contain skyscraper sheaves $\delta_y$ other than $\delta_0$.
Indeed any  $\delta_y\in {\cal T}'$ is a retract of $f_*(P)$ in $P(Y)$, for some $P\in {\cal T}$. Of course, we may assume that $P$ is a simple perverse sheaf.  Since $ \delta_y \hookrightarrow f_*(P)$, $P$ is a perverse constituent of $f^*(\delta_y)$ by adjunction
$Hom_{P(X)}(f^*(\delta_y), P)= Hom_{P(Y)}(\delta_y, f_*(P))$ and corollary \ref{c1}.
Since $f$ is an isogeny, therefore $\dim(supp(P)) =0$. Therefore
$P\cong \delta_x$ for some closed point $x\in X$. This and $P\in {\cal T}$ implies $x\in F$,  hence $y = f(x) = 0$.  

\medskip
By proposition \ref{units} (finite field case) and [KrW, prop.10.1 b)] (complex field case) the invertible objects of the  category $\overline P(X)$ are precisely the skyscraper sheaves.   
We conclude

\begin{Lemma} \label{newlemma} For an abelian variety over $k$ let $K$ be a semisimple complex  on $X$.
If $\pi_0({\bf G}(K))$ is nonabelian, then  $\pi_0({\bf G}(K'))$ is a nontrivial perfect finite group. 
In particular ${\bf G}_{geom}(Y)$ admits a nontrivial
perfect finite quotient group $H$.
\end{Lemma}

We later use this last lemma in section \ref{Drinfeld}
for the proof of theorem \ref{complexcase}.
Of course, theorem \ref{complexcase} finally implies that
$\pi_0({\bf G}(K))$ does not admit perfect quotients, and this
shows in the above argument that $F$ is isomorphic to $\pi_0({\bf G}(K))$.

\goodbreak

\section{Grothendieck rings} \label{Groth}

\medskip
For an abelian variety $X$ of dimension $g$ over $k$,
let  $K^0(P(X))$ be the Grothendieck group  of the category $P(X)$ of semisimple perverse sheaves on $X$. Let $K_*^0(X)$ denote the tensor product  of $K^0(P(X))$ with the polynomial ring $\Z[x,x^{-1}]$. We define a $\Z[x,x^{-1}]$-linear homomorphism 
$$ h: K_*^0(X) \to \Z[x,x^{-1}] \ ,$$
given on generators (i.e. perverse sheaves $K$) by $h(K) = \sum_{\nu \in \Z} \dim_\Lambda(H^\nu(X,K)) \cdot x^\nu $.

\medskip  
The Grothendieck ring $\overline K_*(X)$ of the tensor category 
$\overline P(X)$ tensored with the polynomial ring $\Z[x,x^{-1}]$ 
is the quotient of $K_*^0(X)$ divided  by the classes of negligible perverse
sheaves. If $X$ is simple, by [KrW], [W2] simple negligible perverse sheaves on $X$ are character twists $L_\chi[g]$ of the constant
perverse sheaf $\delta_X = \Lambda_X[g]$ for $g=\dim(X)$. Since $h(L_\chi[g])=0$ holds for nontrivial $\chi$, hence the composite homomorphism $h_*$ of $h$ and the  quotient map by the principal ideal generated by $h(\delta_X)=(x+x^{-1}+2)^{g} $ factorizes as follows
 $$ \xymatrix{   K_*^0(X) \ar@{->>}[d]\ar[r]^-{h} \ar[dr]^-{h_*} &  \Z[x,x^{-1}] \ar@{->>}[d] \cr
\overline K_*(X) \ar@{.>}[r]^-{\overline h_*}  & \Z[x,x^{-1}]/(x+x^{-1}+2)^{g}    } $$ 

Notice $$ \Z[x,x^{-1}] /(x+x^{-1}+2)^{g}  \ \cong \ \Z[y]/y^{2g} $$ 
for $y=1+x$.
Notice,  $-x=1-y$ is inverse to $1+y+ \cdots + y^{2g-1}$ in  $\Z[y]/y^{2g}$.

\medskip
{\it Products of simple abelian varieties}. For a product $X = \prod_{j=1}^r A_j$ of abelian varieties $A_j$ of dimensions $g_j$ 
consider the projections $f_i:X_i =\prod_{j=1}^{i}A_j  \to  X_{i-1} =\prod_{j=1}^{i-1}A_j $. In analogy with the simple case we define
$$ K_*^0(X_i) = K_0(Perv(X_i)) \otimes \Z[x_1,x_1^{-1},...,x_i,x_i^{-1}] $$
and the  ring $R \ =\  \bigotimes_{i=1}^r  \ \Z[x_i,x_i^{-1}]/(x_i  + x_i^{-1} + 2)^{g_i}\ \cong\ \Z [y_1,..,y_r]/(y_1^{2g_1},..,y_r^{2g_r})$.
Again $R_\Q=R\otimes\Q$ is a local Artin ring. 
Define $$h=h_r: K_0^*(X) \to  \Z[x_1,x_1^{-1},...,x_r,x_r^{-1}]$$ and more generally
$$ h_i: K_0^*(X_i) \to  \Z[x_1,x_1^{-1},...,x_i,x_i^{-1}] \quad , \quad i=1,..,r \ $$
by $h_i =  h_{i-1} \circ  p_i$, where
now $  p_i: K_0^*(X_i) \to K_0^*(X_{i-1}) \otimes \Z[x_i,x_i^{-1}]$ 
will be defined as the $\Z[x_1,x_1^{-1},...,x_i,x_i^{-1}]$-linear  extension of the mapping that for perverse sheaves $K\in Perv(X_i)$ is given in terms of the perverse cohomology  $$p_i(K)= \sum_{\nu\in\Z} {}^p H^\nu(Rf_{i*}(K))\cdot x_i^\nu\ .$$
The $p_i$ and therefore the $h_i$ are {\it ring homomorphisms}.
For simple objects $K$ and $L$ this follows from ${}^p H^k(Rf_{i*}(K*L)) \cong \bigoplus_{\nu+\mu=k} {}^p H^\nu(Rf_{i*}(K)) * {}^p H^{\mu}(Rf_{i*}(L))$ in $\overline P(X)$, using the decomposition theorem.  Now, as in the case of simple abelian varieties, we define
$ h_*: K_*^0(X) \longrightarrow R $
as the composite of $h=h_r$ and the projection to $R$. 

\medskip

\begin{Lemma}\label{hom}
{\it For a product $X=\prod_{j=1}^r A_j$  of simple abelian varieties $A_j$ over $k$ the ring morphism $ h_*: K_*^0(X) \longrightarrow R $ 
 $$ \xymatrix{   K_*^0(X) \ar@{->>}[d]\ar[r]^-{h} \ar[dr]^-{h_*} &  \Z[x_1,..,x_r,x_1^{-1},..,x_r^{-1}] \ar@{->>}[d] \cr
\overline K_*(X) \ar@{.>}[r]^-{\overline h_*}  &  R } $$ 
factorizes over a ring homomorphism
$  \overline h_*: \overline K_*(X) \longrightarrow R $ for the quotient ring $\overline K_*(X)$ of $K_*^0(X)$  where $\overline K_*(X)= K^0(\overline P(X)) \otimes \Z[x_1,..,x_r,x_1^{-1},..,x_r^{-1}]$ and $K^0(\overline P(X))$ denotes the Grothendieck ring  of the tensor category
$\overline P(X)$.}
\end{Lemma}

\medskip
{\it Proof}. We have to show $h_*(K)=0$ for all simple negligible perverse sheaves $K$ on $X$. By [W] any negligible perverse sheaf in $P(X)$ is a direct sum  of simple perverse sheaves $K$, each of which is
$A$-invariant for a certain simple abelian variety $A\subset X$ of dimension $\dim(A)>0$, in the sense that  $T_x^*(K) \cong K$ holds for all closed points $x$ in $A$. 
Then all simple perverse components
of ${}^pH^i(Rf_{r*}(K))$ are $f_r(A)$-invariant, by lemma  \ref{dep}. If $f_r(A)$ is not zero, the complex
$Rf_{r*}(K)$ is neligible and our claim follows by induction on $r$ in this case. If $f_r(A)$ is zero, then $A=A_r$ and $K \cong K_{r-1} \boxtimes L_{\chi_r}$ for some perverse sheaf $K_{r-1}$ on $X_{r-1}$ and some twisting perverse character sheaf
$L_{\chi_r}$ on $A_r$ attached to a one dimensional character of the fundamental group
$\pi(A_r,0)$. Then $Rf_{r*}(K) = H^\bullet(A_r,L_{\chi_r}) \otimes_\Lambda K_{r-1}$. Hence the image
$h_*(K)$ of $h(K)$ is zero in $R$, since $p_r(K) = (x_r + x_r^{-1} + 2)^{g_r}\cdot K_{r-1}$
or $pr_r(K)$ is zero. \qed

\medskip
For a homomorphism $f: X\to Y$ between abelian varieties and a pure perverse sheaf
$K$ on $X$ the hard Lefschetz theorem [BBD] implies ${}^p H^\nu(Rf_*(K)) \cong {}^p H^{-\nu}(Rf_*(K))$.  (We here ignore the Weil sheaf structure and hence ignore Tate twists in positive characteristic).
This shows
$$  h_*(K)(x_1,...,x_r) = h_*(x_1,...,x_{i-1},x_i^{-1},x_{i+1},...,x_r) $$
for any $i=1,...,r$. Poincare duality implies ${}^p H^\nu(Rf_*(K^\vee)) \cong {}^p H^{-\nu}(Rf_*(K))^\vee$. Hence $  h_*(K^\vee)(x_1,...,x_r) = h_*(K)(x_1^{-1},...,x_r^{-1})$ and therefore
$$  h_*(K^\vee) = h_*(K) \ .$$

\section{Finite component groups $H$}\label{S4}

%In this section $k$ is either $\C$, or the algebraic closure of a finite field.

\medskip
From  the decomposition ${\bf G}(X) \cong {\bf G}_{geom}(X) \times \mathbb G_m $ the group $\pi_0({\bf G}(X))$ of connected components of ${\bf G}(X)$ is isomorphic to the corresponding group $\pi_0({\bf G}_{geom}(X))$ for the algebraic reductive group ${\bf G}_{geom}(X)$. By definition $\pi_0({\bf G}_{geom}(X))$ is the projective limit of the finite groups $\pi_0(G)$ of Zariski connected components, the limit being taken over all algebraic quotient groups $G$ of ${\bf G}_{geom}(X)$. Notice, $\pi_0({\bf G}_{geom}(X))$ is an abelian group if all these finite component groups $\pi_0(G)$ are abelian as shown in theorem \ref{abelian}, and the converse of course is also true.

\medskip
The representation categories
of the finite quotient groups $H$ of $G$  or ${\bf G}_{geom}(X)$ are full tensor subcategories of the representation category of ${\bf G}_{geom}(X)$. 
Since  $Rep_\Lambda({\bf G}_{geom}(X))$ is a neutral Tannakian category by [KrW] (complex case) resp. [W2], theorem \ref{complexcase} (finite field case), this immediately follows from

\begin{Theorem}\label{abelian}{\it For an abelian variety $X$ over $k$ the irreducible representations of the group $\pi_0(G)$ of connected components of an algebraic quotient group $G$ of ${\bf G}_{geom}(X)$ correspond to skyscraper sheaves $\delta_x$ for certain torsion points $x$ in $X(k)$. In particular
$\pi_0(G)$ is abelian.} 
\end{Theorem}

\medskip
To prove Theorem \ref{abelian}, by Tannakian arguments  it suffices that any full tensor subcategory  in $Rep_\Lambda({\bf G}_{geom}(X))$ isomorphic to $Rep_\Lambda(H)$ for  a finite group $H$  is generated as an additive category by (perverse) skyscraper sheaves $\delta_x$, where the points $x$ are contained in  a finite torsion subgroup  $ H\subseteq X(k)$.

\medskip
{\it  The regular representation of $H$}. Let $H$ be a finite group.
Over the algebraically  closed field $\Lambda$ of characteristic zero the category $Rep_\Lambda(H)$ of representations of $H$ contains the regular representation $\overline K$ of $H$ that is defined by left multiplication of $H$ on the group ring
$ \overline K \ = \ \Lambda[H]$.  This regular representation is self dual in the sense $\overline K \cong \overline K^\vee$ and its dimension is  
the cardinality $d$ of $H$. By Wedderburn's theorem $\Lambda[H]$
 contains all irreducible representations of $H$ up to isomorphism: 
$ \overline K \ \cong \ \bigoplus_{i=1}^h \ d_i \cdot \overline K_i$ holds with $\dim_\Lambda(\overline K_i)=d_i $ for a set of representatives $\overline K_i$ of the isomorphism classes of the  irreducible representations
of $H$ over $\Lambda$. Furthermore
$$ \overline K \otimes_\Lambda \overline K \ \cong \ d \cdot \overline K \ .$$

\medskip
Suppose  $Rep_\Lambda(H)$ is equivalent to a tensor subcategory ${\cal T}_H$ of the Tannaka category $\overline P(X) = Rep_\Lambda({\bf G}_{geom}(X))$ as a tensor category. 
Then representations of $H$ can be considered as objects of $\overline P(X) = Rep_\Lambda({\bf G}_{geom}(X))$. So the regular representation  of $H$ is represented by a semisimple perverse sheaf denoted $K$, and if chosen as a clean perverse sheaf (see the appendix for the notion clean) this perverse sheaf $K\in Perv(X,\Lambda)$ is uniquely determined by $\overline K$. Then $K$ is Tannaka  self dual and $$ K*K \ \oplus\ T'\ =\ d\cdot K\ \oplus\ T \quad , \quad \chi(K)=\dim(K)=d $$ holds for some finite direct sums $T$, $T'$
of complex shifts of negligible perverse sheaves on $X$.  %It is easy to see $T'=0$. 
The perverse constituents $K_i \in P(X)$ 
representing the irreducible constituents $\overline K_i \in \overline P(X)$   are uniquely defined by the $\overline K_i$ (up to isomorphism), and
$$K \ = \ \bigoplus_{i=1}^h \ d_i \cdot K_i \quad , \quad \chi(K_i)=d_i > 0$$  is a semisimple perverse sheaf on $X$.   

\medskip
{\bf Definition}. {\it A semisimple perverse sheaf, obtained in the way described above, will be called a {\it $H$-regular} perverse sheaf on $X$.}

\medskip
Since twisting with a character sheaf $L_\chi$
defines a tensor functor (see [KrW]),
character twists $K_\chi$ of $H$-regular perverse sheaves $K$ are $H$-regular perverse sheaves.

%\section{Strongly $d$-regular perverse sheaves}\label{num6}

\medskip {\bf Definition}.  For an integer $d\geq 1$ a semisimple sheaf complex $K\in D_c^b(X,\Lambda)$ will be  called {\it weakly $d$-regular} if $K^\vee \cong K$ and $$K*K \ \oplus \ T' \ \cong \ d \cdot K \ \oplus\ T$$ hold for direct sums $T$ and $T'$ of complex shifts of negligible perverse sheaves on $X$.  We say, $K$ is {\it $d$-regular} if $T'=0$, and {\it strongly $d$-regular} if $K$ is a perverse sheaf and $T'=T=0$. A $H$-regular and strongly $d$-regular perverse sheaf will be called
{\it strongly $H$-regular}. 

\medskip
For any weakly $d$-regular sheaf complex $\chi(K)=d\geq 1$ holds.
If $K$ is  a $H$-regular perverse sheaf $K$, then $m\cdot K$ is weakly $d$-regular for $d=m\cdot \# H$.

\medskip % In the following assume $k=\C$.
Suppose $K$ is weakly $d$-regular on $X$. Then the subcategory of $\overline P(X)$ defined by finite direct sums of the images $\overline K_i$ of the simple constituents $K_i$ of $K$ defines a tensor subcategory of $\overline P(X)$ with finitely many isomorphism classes of irreducible objects, thus corresponds to some finite quotient group $H$ of ${\bf G}(X)$. 

\begin{Lemma} \label{per} Let $K\in D_c^b(X,\Lambda)$ be semisimple and weakly $d$-regular. If $K=M \oplus R$ is the decomposition into a clean complex $M\in D_c^b(X,\Lambda)$ and a negligible complex $R$, then $M$ is a perverse sheaf. 
\end{Lemma}

{\it Proof}. $K=\bigoplus_\nu M_\nu[-\nu]$ with $M_\nu={}^p H^\nu(K)$. Suppose $\nu$ is maximal (or minimal) such that $M_\nu \notin N_{Euler} $. Then $M_\nu[-\nu]$ defines a representation of  ${\bf G}(X)= {\bf G}_{geom}(X) \times \Gm $. Its restriction to $\Gm$ is an isotopic multiple of the character $t\mapsto t^{-\nu}$. The representation associated to the tensor product $K*K$ therefore has a nontrivial eigenspace for the character $t\mapsto t^{-2\nu}$ of $\Gm$. Since $0\neq {}^p H^{2\nu}(K*K) \cong d\cdot {}^p H^{2\nu}(K)$ holds in $\overline D^{ss}(X)$, this contradicts the maximality of $\nu$
unless $\nu=0$.  Hence $K$ is the sum of a clean perverse sheaf and a negligible complex. \qed

\begin{Lemma} \label{reduce}
Any weakly $d$-regular semisimple sheaf complex $K$ is of the form $$K\ = \ M \ \oplus\ R$$ for a 
a negligible sheaf complex $R$ and  clean $d$-regular perverse sheaf $M$.
\end{Lemma}

{\it Proof}.  The decomposition $K= M\oplus R$ of lemma \ref{per} implies $M^\vee \cong M$  and $$M^{*2} \oplus (2\cdot M*R\oplus R^{*2}\oplus T')\ \cong\ d\cdot M \oplus (d\cdot R\oplus T)\ .$$  Negligible object define a tensor ideal, hence $(2\cdot M*R\oplus R^{*2}\oplus T')$
is negligible and since $d\cdot M$ is clean, it must occur as a direct summand of $d\cdot R\oplus T$. This implies $M^{*2} \cong d\cdot M \oplus T''$  
for some negligible complex $T''$. \qed

\begin{Lemma} \label{dirim}  For an isogeny $g:X\to Y$ % of  complex 
% abelian varieties 
the direct image $g_*(K)$ of a $d$-regular semisimple perverse sheaf $K$ on $X$ is a $d$-regular semisimple perverse sheaf on $Y$. Twists $K_\chi$  of $d$-regular perverse sheaves $K$ are $d$-regular perverse sheaves.
\end{Lemma}

{\it Proof}. By the decomposition theorem  $g_*(K)$ is semisimple. Since $g_*$ is a tensor functor, the claim is obvious because  $g_*(T)$ is negligible by lemma \ref{7}.5.  \qed  

\medskip
This being said, we come to a key step in the proof of theorem \ref{complexcase}. By lemma \ref{reduce} the proof of theorem \ref{complexcase}
amounts to an inductive proof (induction on $n$) of the following two assertions
\label{induction}

\medskip\noindent
{\bf Assertion $Tor(n)$}. {\it 
For all % complex 
abelian varieties $X$ of dimension $\dim(X)\leq n$ any clean $d$-regular semisimple sheaf complex $K$ on $X$ is a perverse skyscraper sheaf concentrated in torsion points of $X(k)$.}

\medskip\noindent
{\bf Assertion $Reg(n)$}. {\it 
For % complex 
abelian varieties $X$ of dimension $\dim(X)\leq n$ clean $d$-regular semisimple sheaf complexes $K$  on $X$ are strongly $d$-regular, i.e. satisfy $$ \fbox{$ K^{*2} \ \cong \ d \cdot K \quad , \quad K \cong K^\vee $}\ .$$}

\medskip Obviously $Tor(n)$ implies $Reg(n)$. More importantly, we have a partial converse given by the next proposition \ref{regular}. Notice, Proposition \ref{regular} allows us to assume  for the
proof of the induction step $Tor(n) \Longrightarrow Tor(n+1)$ (that will be later given in the sections \ref{Drinfeld} and \ref{congruences}) that the assertion $Reg(n+1)$ already holds.  

\begin{Proposition} \label{regular} $Tor(m)$ for all $m < n$ implies
$Reg(n)$. 
\end{Proposition}

{\it Proof}. Assume $\dim(X)=n$.  By lemma \ref{reduce} a clean semisimple $d$-regular complex $K$ on $X$ is a perverse sheaf so that $K*K\cong d\cdot K \oplus T$ and $K\cong K^\vee$ holds. We have to show $T=0$. We freely use the notion of the appendix on translation invariant perverse sheaves.

\medskip
{\it First step} (where $X$ is a simple abelian variety). Then we 
claim %that the cohomology groups of $K$ and of all its character twist $K_\chi$
%vanish in the degrees $\neq 0$
$$ H^\nu(X,K_\chi) = 0  \quad , \quad   (\forall \chi \ , \forall \nu\neq 0) \ .$$
Since the assumptions on $K$ are inherited by the $K_\chi$, we may assume $\chi=1$.
In terms of the ring homomorphism
$  h: K_*^0(X) \to   \Z[x,x^{-1}] $ 
defined in section \ref{Groth}, we have to show that $h(K)$ is a constant. To this end, we may extend coefficients from $\Z$ to $\Q$ and (by abuse of notation)
consider the extended map $h: K_*^0(X)\otimes \Q \to   \Q[x,x^{-1}]$.  
We first show that the image of $K$ under  $$\Q[x,x^{-1}] \mapsto  R_\Q=\Q[x,x^{-1}]/(x+x^{-1} +2)^{n} $$ becomes constant. By lemma \ref{hom} this allows to consider the ring homomorphism    
$$  h_*: K_*^0(X)\otimes \Q \to   R_\Q \ .$$ 
Since $h_*$ is trivial on negligible
complexes, $P^2 =d\cdot P$ holds for $P=h_*(K)$ in $R_\Q$. Since $R_\Q$ is a local Artin ring, hence $P$ is either in the nil radical or $P$ is a unit. In the first case $P$ is $0$ in $R_\Q$, since $P^2=d\cdot P$ implies $d^{2n-1}P=P^{2n} = 0$ in $R_\Q$. In the second case $P$ is invertible in $R_\Q$, hence $P=d$ in $R_\Q$. So in both cases, $h_*(K)$ is constant.

\medskip
To finish the proof for our claim it remains to show that this implies that $h(K)$ is either zero or equal to $d$ in $\Q[x,x^{-1}]$. Since
$K$ is semisimple perverse sheaf without translation invariant constituents, the perverse cohomological dimension bounds [BBD, 4.2.4] %or [KW], theorem III.11.3
give
$$ h(K)= \sum_{\nu=-n+1}^{n-1} \dim( H^\nu(X,K)) \cdot x^\nu \ \ \in \ \Q[x,x^{-1}]\ .$$
Since $h_*(K)=a$ is constant, there exists a constant $a\in \Q$ (either $a=0$ or $a=d$ in $\Q$) and a Laurent polynomial $f$ in $x$ so that
$ h(K) - a = f(x,x^{-1}) \cdot (x+1)^{2n}$.
If  $h(K)\neq a$ were not constant, then $h(K)-a = x^{-j} \cdot h(x)$ holds for some polynomial $h(x)$
of degree $\leq 2n-2$, where $j$ is determined by the normalization condition $h(0) \neq 0$. If we write $f(x,x^{-1})= x^{-i} \cdot g(x)$
for a polynomial $g(x) \in \Q[x]$ such that $g(0) \neq 0$,   
then $x^{-j}h(x) = x^{-i} (x+1)^{2n} g(x)$. Therefore 
$x^{i-j} h(x) =  (x+1)^{2n} g(x)$ implies $i=j$, by putting $x=0$. Hence $h(x)=(x+1)^{2n} g(x)$.  Since $deg(h) \leq 2n-2$, comparing degrees gives a contradiction unless $h(x)=g(x)=0$.
Thus $h(K)=a$ and our claim $  H^\nu(X,K) = 0$, $ \nu\neq 0 $  follows,
and the same for $K$ replaced by $K_\chi$.

\medskip
By the K\"unneth formula, the property $H^\nu(X,K_\chi)=0$ for all $\chi$ and all $\nu\neq 0$ and $K^{*2}=d\cdot K \oplus T$ implies 
$H^\nu(X,T_\chi)=0$ for all $\nu\neq 0$. 
Since $\chi(T_\chi)=0$, $H^0(X,T_\chi)=0$ follows for all $\chi$. Hence all $T_\chi$ are acyclic and all their negligible simple perverse
constituents. Since an irreducible translation invariant semisimple  complex is acyclic for all its character twists if and only if it is zero,
this forces $T=0$. This completes the proof in the case of simple abelian varieties.

\medskip
{\it Second step}.  Let $X$ be isogenous to a product $\prod_{i=1}^r A_i$
of abelian varieties. Let $g: X \to \prod_{i=1}^r A_i$ be such an isogeny.  
To show $T=0$ for a $d$-regular semisimple clean perverse sheaf on $X$, we may replace
$X$ by $\prod_{i=1}^r A_i$, $K$ by $g_*(K)$ and $g_*(T)$ without restriction of generality
using lemma \ref{dirim} and lemma \ref{7}. So, let us assume $X=\prod_{i=1}^r A_i$ and $r>1$
for a simple abelian variety $A_1$.

\medskip
In this second step we show $T=0$  under {\it the assumption $Stab(T)=X$}. For the projection  $q: X \to B= A_1$ with kernel
$A=\prod_{i=2}^r A_i$ let $L=Rq_*(K) = M \oplus R$  be the direct image complex 
with its decomposition into a clean perverse sheaf $M$ and a negligible complex $R$. The assumption $K*K = d\cdot K \oplus T$ and lemma \ref{reduce} imply that $M$ is a clean semisimple $d$-regular
complex 
$$   M^{*2} \ \oplus\ \bigl( 2 \cdot M*R \oplus  R^{*2} \bigr) \ =\  d\cdot M\ \oplus\ \bigl( d\cdot R \oplus Rq_*(T) \bigr) \ $$
where the sheaf complexes  within the brackets are negligible complexes. By our induction 
assumption (for the induction start use step 1) the assertion
$Reg(m)$ for $m< n$ implies $M^{*2} \cong d\cdot M$. This leads to the equation (*) 
$$ (2 \cdot M*R) \ \oplus\ R^{*2}  \ =\   d\cdot R \ \oplus\ Rq_*(T) \ .$$
$Stab(T)=X$ implies $Stab(Rq_*(T))=B$ by lemma \ref{7}, part 4).  
Since $B$ is simple and $R$ is negligible, also 
$Stab(R)=B$. For certian graded $\Lambda$-vectorspace $V_\chi^\bullet$ we get a decomposition  $R= \bigoplus_\chi R_\chi$ into $\chi$-blocks
with  $$R_\chi := V_\chi^\bullet \otimes_\Lambda L_\chi $$ where $L_\chi$ are the rank one negligible 
perverse sheaves 
$L_\chi = \Lambda_B[\dim(B)]_\chi$ 
for the characters $\chi$ of $\pi_1(B,0)$.
Similarly $Rq_*(T)= \bigoplus_\chi S_\chi$
for 
$$S_\chi := W_\chi^\bullet \otimes_\Lambda H_A^\bullet \otimes_\Lambda L_\chi \ .$$
Again, $W_\chi^\bullet$ are certain graded $\Lambda$-vectorspaces and $H_A^\bullet$ abbreviates $H^\bullet(A,\Lambda)[\dim(A)]$. 
Since $L_{\chi} *L_{\chi'} =0$ for $\chi\neq \chi'$ and $L_\chi^{*2} = H_B^\bullet \otimes_\Lambda L_\chi$ for $H_B^\bullet = H^\bullet(B,\Lambda)[\dim(B)]$, equation (*) completely decouples into the 
corresponding equations (*)${}_\chi$
$$ (2\cdot M*R_\chi) \ \oplus\ R_\chi^{*2}  \ =\   d\cdot R_\chi\ \oplus\ S_\chi \ $$
for each $\chi$-block.  
Notice $$M*R_\chi := (H^\bullet(B,M_\chi) \otimes_\Lambda V_\chi^\bullet) \otimes_\Lambda L_\chi \quad \mbox{ and }\quad R_\chi^{*2} = (V_\chi^{\otimes 2} \otimes_\Lambda H_B^\bullet)\otimes_\Lambda L_\chi \ .$$  
This decoupling allows us to suppress the index $\chi$ for rest of the argument in this step:
For the Poincare polynomials $P_U(x) \in \Z[x,x^{-1}]$  of graded $\Lambda$-vector spaces $U=U^\bullet$, the equations (*)${}_\chi$ for each $\chi$ in the Laurent ring $\Q[x,x^{-1}]$  become 
$$ P_V(x) \cdot \Bigl( 2 P_M(x) - d + P_V(x)\cdot (2+x+x^{-1})^{\dim(B)} \Bigr) \ =\ P_W(x)
\cdot \Bigl(2+x+x^{-1}\Bigr)^{\dim(A)} \ .$$
 %We now proceed as in step 1).  
Since $K$ is perverse and clean and $R$ is the negligible component of $L=Rq_*(K)$, there are bounds  $$ 1 -\dim(A) \leq  deg_x(P_V) \leq \dim(A) -1\ .$$ Clearing the powers of $x$ in the denominators,
gives an equation in the polynomial ring $\Q[x]$. %We compare the factorizations with respect to the prime polynomial $(x+1)$ on both sides of the identity above. 
The right side is divisible by $(x+1)^{2\dim(A)}$. The polynomial obtained from $2P_M(x) - d$ is not divisible by $(x+1)$. Indeed $P_M(-1) = \chi(M) =
\chi(Rq_*(K))=\chi(K) = d$ implies $2P_M(-1)-d = d \neq 0$. On the other hand $dim(B)>0$, so the bracket on the left side is not divisible by $(x+1)$. Thus $P_V(x)$ must be divisible by $(x+1)^{2\dim(A)}$. After clearing denominators, $P_V(x)$ becomes a polynomial in $x$ of degree at most $2\dim(A)-2$ in $\Q[x]$. Being divisible $(x+1)^{2\dim(A)}$ it must be zero. This implies $P_W(x)=0$. Since this holds for all $\chi$-blocks, we obtain
$S=Rq_*(T)=0$.  Finally, using twists by all characters of $\pi_1(A,0)$ and the assumption $Stab(T)=X$ similarly $Rq_*(T)=0$ then implies $T=0$. This completes the second step.

\medskip
{\it Last step} (general case). We now show that the assumption $Stab(T)=X$ of the second step holds unconditionally. If this were not true, there exists  a constituent  $N\neq 0$ of $T$ with $Stab(N)^0\neq X$ and $depth(T)=dim(Stab(N)^0)>0$. 
Choose an abelian variety $A\subset X $ and an isogeny $A \times Stab(N)^0 \to X$ so that
the quotient map $q: X \to B=X/A$ defines an isogeny $q: Stab(N)^0 \to B$.
As in the preliminary remarks made in step 2 of the proof, we can easily reduce the proof to the case where $X$ is a direct product $$ X=A \times B \quad , \quad  B=Stab(N)^0$$ so that $q$ is the projection onto the second factor. Then by lemma \ref{nice} there exists a character $\chi_0$ of $\pi_1(X,0)$ such that for {\it most characters} $\chi$ of $\pi_1(A,0)$ the complex $$ S(\chi)=Rq_*(T_{\chi_0\chi})$$ is is {\it not acyclic}, but {\it negligible} of depth $\dim(B)$ on $B$. In particular, the semisimple complex $S(\chi)$ is a direct sum of complex shifts of perverse sheaves that are all invariant under $B$. Since $B$ is not acyclic, it contains an irreducible constituent of $S:=S(\chi)$ of the form $\delta_B[-\nu]$ for the constant perverse sheaf $\delta_B = \Lambda_B[\dim(B)]$ on $B$
$$ \delta_B[-\nu] \hookrightarrow S(\chi) \ .$$ 
Notice, the complex shift $\nu=\nu(\chi)$ may depend on $\chi$. 
Since $K$ and all its twists are $d$-regular, $M:=Rq_*(K_{\chi_0\chi})$ again is a $d$-regular complex: 
$$    M*M \ \cong  d\cdot M \oplus S  \quad \mbox{ for } \quad        S:=   Rq_*(T_{\chi_0\chi}) \ .$$
By lemma \ref{reduce} the complex $M$ decomposes into a clean perverse sheaf $M_{clean}$ and a negligible complex $R$ 
$$M= M_{clean} \oplus R $$  
depend on $\chi$. So let us fix $\chi$.
Since $\dim(B) < \dim(X)$, we can apply the induction assumption $Reg(m)$ for $m<n$ to show
that $M_{clean}^{*2} = d \cdot M_{clean}$.  If we insert this into $M*M \ \cong  d\cdot M \oplus S$, we obtain 
$$    (2 \cdot M_{clean}*R) \ \oplus \ R^2 \ =\ d\cdot R \ \oplus S \ .$$ 
By the induction assumption $Tor(m)$ the perverse sheaf $M_{clean}$ is a skycraper sheaf and $\chi(M_{clean})=d$. Hence, $2\cdot  M_{clean} *R$ is a direct sum of $2d$ translates of $R$. If $R\neq 0$, therefore the left side can not be completely contained in the summand $d\cdot R$ on the right side. So $R\neq 0$ implies that some translate of each simple constituent of $R$ also occurs as constituent in $S$. Since $S$ is translation invariant under $B$, for $R\neq 0$ this implies
$$ Stab(R)=B \ .$$ Hence, up to complex shifts, all simple constituents of $R$ are of the form $L_\psi$ for characters $\psi:\pi_1(B,0)\to \Lambda^*$. However,
we claim that this negligible complex $R$ is {\it not acyclic}. 
If it were, $(2 \cdot M_{clean}*R) \ \oplus \ R^2$ and therefore also $ d\cdot R \ \oplus S$ would be direct sums of the $L_\psi[\mu]$ for $\psi\neq 1$, because $L_\psi * P  = H^\bullet(B,P_\psi) \cdot L_\psi$ holds for all $P$. But this can not be, since then
$S=S(\chi)$ would be acyclic as opposed to our  construction. Indeed, the  semisimple negligible but  not acyclic complex $S$ contains a summand of the form $\delta_B[-\nu_0]$. Since $R$ is not acyclic, we get
$$  \  M \ = Rq_*(K_{\chi_0\chi}) \ = \ M_{clean} \ \oplus\  R \ = \ M_{clean} \ \oplus\ \delta_B[-\nu_0] \ \oplus \ ...  \ $$
and by  the hard Lefschetz theorem we may assume $\nu_0=\nu_0(\chi)\geq 0$. Hence
the Leray spectral sequence implies 
$$ \bigoplus_{i\geq \dim(B)} H^{i}(A\times B,K_{\chi_0\chi}) \ \cong\  \bigoplus_{i\geq \dim(B)} H^{i}(B,M) \neq 0  \ .$$
If we compute $H^{i}(A\times B,K_{\chi_0\chi})$ by the perverse Leray 
spectral sequence for the morphism $p$ (which degenerates by the decomposition theorem) 
$$  \xymatrix{  X= A\times B \ \ar[d]_-{p}\ar[r]^-{q} &  \ \ B \ \ar[d]^-{\tilde p} \cr
 A \ \ \ \ar[r]^-{\tilde q} &  \ \ Spec(k)  \cr } \ ,$$
a comparison  gives 
$$ \bigoplus_{\nu+\mu \geq\dim(B)} H^{\nu}(A,{}^p R^\mu p_*(K_{\chi_0\chi}))\ \cong \ \bigoplus_{i\geq \dim(B)} H^{i}(A\times B,K_{\chi_0\chi}) \ \cong\  \bigoplus_{i\geq \dim(B)} H^{i}(B,M) \neq 0 \ .$$
Now, we can still vary the character $\chi$. By our constructions this conclusions 
remains valid for most characters $\chi$ of $\pi_1(A,0)$. 
On the other hand, by the vanishing theorem \ref{relVT} for all $\mu$ and all $\nu\neq 0$ the group $$H^\nu(A,{}^p R^\mu p_*(K_{\chi_0\chi})) =
H^\nu(A,{}^p R^\mu p_*(K_{\chi_0})_{\chi})$$ vanishes for most $\chi$. 
If we insert this information into the above nonvanishing result, for most $\chi$ we obtain 
$$\bigoplus_{i\geq \dim(B)} H^0(A,{}^p R^{i} p_*(K_{\chi_0\chi}))\neq 0\ $$ and hence ${}^p R^i p_*(K_{\chi_0\chi}) \neq 0$ for some $i\geq \dim(B)$. Since $p$ is smooth of relative dimension $\dim(B)$ and $K_{\chi\chi_0}$ is perverse, [BBD], 4.24  
therefore implies $i=\dim(B)$ so that  
$K_{\chi\chi_0}$ contains a nontrivial $B$-invariant constituent. Since $K$ is clean, by contradiction this  proves $depth(T)=\dim(X)$. Hence $Stab(T)=X$.
\qed

\bigskip\noindent

\goodbreak

\section{From $\C$ to finite fields}\label{Drinfeld}

\bigskip\noindent
In this section we reduce the proof of theorem \ref{abelian} to the case
where $k$ is the algebraic closure of a finite field $\kappa$. This amounts to
discuss strongly $H$-regular clean semisimple perverse sheaves
on abelian varieties by the induction procedure given on page \pageref{induction}.
By lemma \ref{newlemma} our task will be to  show the following assertion:  

\medskip
{\it There do not exist strongly $H$-regular clean semisimple perverse sheaves for (nontrivial) perfect groups $H$}.

\medskip
Since the case of abelian groups $H$ is completely discussed in [KrW], 
this implies the crucial induction step $$ Reg(n) \Longrightarrow Tor(n)$$ described
on page \pageref{induction}ff.
Combined with proposition 1 this implies that $Tor(m)$ for all $m< n$ gives $Tor(n)$, which
by induction proves theorem \ref{complexcase} resp. theorem \ref{abelian}.

\medskip
The aim in this section is to show that it suffices to prove the crucial assertion above
over finite fields.   
For this we use Drinfeld's approach [Dr] to produce from a  counterexample over $\C$ an analogous counterexample for an abelian variety over a finite field $\kappa$. The idea is the following:

\medskip
Suppose there exists a counterexample $K$ over $\C$, i.e. a strongly $H$-regular clean perverse sheaf $K\in Perv(X,\C)$  for a nontrivial perfect group $H$.
Then the irreducible perverse constituents $K_i$ become smooth
on some open dense smooth subsets $U_i$ of the support $Y_i$ of $K_i$.
There exists a scheme $I=\prod_{i=1}^h Irr_{d_i}^{U_i}$ of finite type over $Spec(\Z)$ 
representing the sheaf functor on the category of commutative rings $A$ associated to the presheaf functor $A \mapsto \underbar I(A) = \prod_{i=1}^h \underbar{Irr}_{d_i}^{U_i}(A)$, where $\underbar {Irr}_d^U(A)$ is the set of isomorphism classes of rank $d$ locally free sheaves of $A$-modules $N$ on $U(\C)$ such that
$N\otimes_A k$ is irreducible for every field $k$ equipped with a homomorphism $A\to k$.   
If $A$ is a local complete ring with finite residue field, for example the completion of a local ring of a closed point of $I$,  then $ \underbar{Irr}_{d_i}^{U_i}(A) \cong Irr_{d_i}^{U_i}(A)$.

\medskip 
Consider the \lq{bad}\rq\ subset $B^{\Q} \subset I\otimes \Q$ of $F$-points (for extension fields $F$ of $\Q$) of strongly $H$-regular clean perverse sheaves
 $K=\bigoplus_{i=1}^h d_i \cdot K_i$, where $K_i$ is the perverse intermediate extension
of the smooth perverse sheaf on $Y_i \subset X$ defined by the underlying etale sheaf of $F$-modules $M_i$, where $\prod_i M_i \in I(F)$. Using the results of [Dr, section 3], in a similar way as for the proof of lemma 2.5. in loc. cit. one shows (now using [Dr], lemma 3.10, lemma 3.11 and the proof of lemma 3.1) that $B^{\Q}$ is a constructible subset of $I\otimes \Q$. Then 
define $B$ to be the Zariski closure of $B^{\Q}$ in $I$. Since by assumption there exists a counterexample, we have $B^{\Q}\neq \emptyset$ and there exists
a Zariski open subset $V\subseteq I$ such that $B^{\Q}\cap (V\otimes\Q)$ is closed in $V\otimes\Q$ and $B\cap V$ is nonempty and smooth over $\Z$.  For any closed point $z\in B\cap V$ 
consider the completion $\widehat I_z$ of $I$ at $z$ and the locus $\widehat B_z$ defined by $B$ in $\widehat I_z$; the complete local ring $A_z$ of $\widehat B_z$ has a finite residue field $\kappa_l$. 
The closed point $z$ can be chosen so that $\kappa_l$ is a 
finite field of arbitrary characteristic $l \geq l_0$ for some $l_0$, so we may assume $l > d=\# H$.

\medskip
Choose a suitable finitely generated field $E$ with algebraic closure $\overline E \subset \C$
over which $U_i,Z_i,X,$ are defined so that $z$ is $Gal(\overline E/E)$-invariant. Then $Gal(\overline E/E)$ acts on $A_z$, i.e. $\widehat I_z$. If $E$ is chosen big enough, 
then $Gal(\overline E/E)$ acts on $\widehat B_z$ as in [Dr, lemma 2.7]. Now we recall the key
point of the argument in [Dr]: The fixed point locus of $F^k$ on $\widehat B_z$ of any Frobenius
substitution $F=F_{\kappa}$ for closed points v of a model of $X$ with residue field 
$\kappa$ 
defined over a finitely generated ring $R$ with quotient field $E$, 
is finite and flat over $\Z_l$ and  nonempty for $k$ large enough.
To show this, Drinfeld used the (by now proven [BK], [G]) de Jong's conjecture which implies that the reduction mod $l$ of the fixed point scheme is finite over $\kappa_l$ and from which then easily follows the assertion above (as in [Dr, lemma 2.8]).  If nonempty, this fixed point scheme defines an $l$-adic perverse sheaf $\tilde K$ on a model of $X$ defined over a suitable localizations of $R$ with closed point v, so that $\tilde K$ and 
its reduction $\tilde K_{0}$ defined over some finite extension of the residue field $\kappa=\kappa_q$ of v (of characteristic $p$) are \lq{bad}\ , i.e. satisfy $K*K\cong d\cdot K \cong d\cdot K^\vee$ without being skyscraper sheaves  
such that their associated Tannaka group is $H$. This part of the argument
follows as in [BBD, section 6] and [Dr, section 4 - 6]. The residue characteristic $p=char(\kappa)$ of the point v can be chosen arbitrarily large $p \geq p_0$, since v can be chosen to be an arbitrary closed point of the spectrum of the ring $R$, which is finitely generated over $\Z$; in fact we can choose  $p > d= \# H$ and $l\neq p$ suitably. 

\medskip
This construction of [Dr] and [BBD] reduces the proof of the characteristic zero assertion $Reg(n)$
for strongly $d$-regular clean perverse sheaves $K\in Perv(X,\C)$ to the proof of the corresponding
assertion for strongly $d$-regular clean semisimple perverse $l$-adic sheaves $K_0$ on abelian varieties over finite fields $\kappa$. In positive characteristic a semisimple perverse sheaf $K$ is called clean, if  no character twist $K_\chi$ contains acyclic irreducible constituents  (over $\C$ this is equivalent to the previous notion). 

\medskip
So to complete the proof
of theorem \ref{complexcase}, via the induction argument using proposition \ref{regular} for the relevant
conclusion $Reg(n) \Longrightarrow Tor(n)$ (see section \ref{S4}), in view of the next lemma \ref{descend} finally it will be enough to show 

\begin{Theorem} \label{4} Let $X_0$ be an abelian variety over a finite extension field
of a finite field $\kappa$ of characteristic $p$ and let $K_0$ be a perverse $\overline{\Q_l}$-adic
sheaf on $X_0$ for $l\neq p$. Let $X$ be the scalar extension of $X_0$ to $k$, and let $K$ be the extension of $K_0$ to a perverse sheaf on $X$. 
Then $K$ is not a strongly $H$-regular clean semisimple perverse sheaf on $X$ for some finite nontrivial perfect group $H$.
\end{Theorem}

\medskip
This theorem will follow from the subsequent lemma \ref{descend} and the arguments given in the sections \ref{Fourier} and
\ref{congruences}.

\bigskip\noindent

{\it $H$-regular descent to finite fields}. \label{arith}
For an abelian variety $X_0$ over a finite field $\kappa$ let
$X$ be a fixed base extension to the algebraic closure $k=\overline\kappa$.
In what follows, a perverse sheaf $K_0\in Perv(X_0,E_\lambda)$ with coefficients
in $E_\lambda$ will often be viewed as a perverse sheaf with coefficients in $\Lambda =
\overline{\Q}_l$. Suppose now, that the base field extension $K\in Perv(X,\Lambda)$ of $K_0$
is a strongly $H$-regular semisimple clean perverse sheaf for a nontrivial perfect group $H$. 
(This is what Drinfeld's construction gave us in the last section). We want to descend $K$ to a
strongly $H$-regular semisimple clean perverse sheaf over some finite field.

\medskip
{\it Conventions and notations}. 
Fix
$\Lambda= \overline{\Q_l}$ and some isomorphism $\tau:\Lambda \cong \C$, which allows  to define complex conjugation on $\Lambda$. Usually we suppress to write $\tau$ and write $\overline{\alpha}$ instead of $\tau^{-1}\overline{\tau(\alpha)}$.
Let $K$ be a perverse $\Lambda$-adic Weil sheaf on $X$, which is $\tau$-pure of integral weight $w$. This means that $K$ is equivariant with respect to the Frobenius 
$F_X$ so that exists an isomorphism $F^*: F_X^*(K) \cong K$. If $K$ is an irreducible perverse sheaf on $X$, then two such Weil sheaf structures $(K,F^*_1)$ and $(K,F_2^*)$
on $K$ yield an automorphism $F_2^*\circ (F_1^*)^{-1}: K\cong K$ which is
given by $\alpha \cdot id_K$ for some $\alpha\in \Lambda^*$ because $End_{Perv(X,\Lambda)}(K)=\Lambda \cdot id_K$.   

\medskip \label{twis}
We write $\alpha \in \Lambda^*_{mot}$
if  $\alpha$ is contained in a finite number field $E \subset \Lambda$ and $\alpha, \alpha^{-1}$
are integral over $\Z[p^{-1}]$, where $p$ is the characteristic of $\kappa$.
We say that $\alpha, \alpha'\in \Lambda^*$ are equivalent if $\alpha'/\alpha$ is a root of unity. Notice that for a simple perverse Weil sheaf $P$ on $X$ (with respect to $\kappa$)
the generalized Tate twists $P(\alpha')$ and $P(\alpha)$ become isomorphic as perverse Weil sheaves over some finite extension field of $\kappa$ iff $\alpha'$ and $\alpha$ are equivalent.
For the scalar extension $P$ of a perverse sheaf $P_0$ over $\kappa$ there exists $\alpha\in \Lambda^*$ such that for $P_0=S_0(\alpha)$ satisfies the \lq{determinant condition}\rq\, i.e.
the determinant of the underlying smooth
coefficient system of $S_{0}$ has finite order. This follows from the corresponding statement for the smooth coefficient system defining $P_0$ (defined on an normal open subset of the support of $P_0$) and follows from the known structure of abelianized Weil groups of normal schemes
[D3]. Up to equivalence, $\alpha$
is uniquely defined by $P_0$ and by [L,cor.VII.8]
the perverse sheaf $S_{0}$ is pure of weight zero.

\medskip
{\it Weakly motivic sheaves and complexes}. Following [Dr2, def.1.7 and 1.8] 
a $\Lambda$-adic sheaf on $X_0$, and similar a complex of $D_c^b(X_0,\Lambda)$
will be called weakly motivic if all eigenvalues 
of Frobenius are Weil numbers in $\Lambda^*_{mot}$ for all closed points
$x$ of $X_0$. As shown in [Dr2, 1.4.2 and app. B], the derived category 
of weakly motivic complexes in  $D_c^b(X_0,\Lambda)$ is stable under the \lq{six functors}\rq ,
under the perverse truncations functors ${}^p\tau_{\leq n}$ and the passage to perverse subquotients.
Furthermore, a smooth $\Lambda$-adic sheaf is weakly motivic if it satisfies the determinant condition.

\medskip
{\it Mixedness for Weil sheaves}. Although we are only interested in the case, where $K$ is an irreducible perverse sheaf obtained from a $\Lambda$-adic
perverse sheaf $K_0$ on $X_0$ by extension of scalars, it is convenient to view them as Weil sheaves over $X$. The scalar extensions $K$ of $\Lambda$-adic
perverse sheaf $K_0$ on $X_0$ are mixed Weil sheaf by a deep result of Lafforgue [L].
If the perverse sheaf $K$ is irreducible on $X$, then it is a pure Weil sheaf. By 
a formal half integral Tate twist 
we can and will always assume that the weight is $w=0$. 
Two Weil sheaf structures (both of weight $0$) for an irreducible perverse
sheaf $K$ on $X$ differ by a generalized Tate twist defined by some $\alpha\in \Lambda$ such that $\vert \tau(\alpha) \vert = 1$. For simplicity we also write $\vert \alpha\vert =1$,
if $\tau$ is fixed.

\begin{Lemma} \label{descend} Let $H$ be a nontrivial perfect finite group.
Then for a strongly $H$-regular clean semisimple perverse sheaf $K$ in $Perv(X,\overline\Q_l)$ on an abelian variety $X$ over $k$ that is defined over a finite subfield of $k$, there exists 
a strongly $H$-regular weakly motivic
pure perverse sheaf 
$K_0' \cong \bigoplus_{i=1}^{h'} \ d_i \cdot S_{0,i}(\alpha_i) $ 
of weight zero in $Perv(X_0,\overline\Q_l)$ defined over a finite subfield $\kappa$ whose base extension $K'$ to $k$ is a strongly $H'$-regular perverse sheaf in $Perv(X,\overline\Q_l)$, such that $H'$ is a nontrivial perfect quotient group of $H$. Furthermore  we can assume that $F_\kappa$ acts trivially on $End_{Perv(X,\Lambda)}(K')$.
\end{Lemma}   

{\it Proof}.
Then $K \cong \bigoplus_{i=1}^h d_i\cdot K_i$ for simple perverse sheaves $K_i$. Without restriction of generality we can
assume  that each $K_i$ is defined over $\kappa$, suitably
enlarging $\kappa$ if necessary. We choose absolutely irreducible
perverse sheaves $K_{0,i}$ on $X_0$ whose extension to $X$ is $K_i$ for $i=1,..,h$ and define 
$$K_0 = \bigoplus_i \ d_i \cdot K_{0,i} $$ so that $K_{0,1}=\delta_0$ and
$   K_{0,i} \cong K_{0,j}  \ \Longleftrightarrow \ K_i \cong K_j $.

\medskip
For the induced Weil sheaf structure on $K$, the Frobenius morphism $Fr=Fr_{\kappa}$ acts trivially on $End_{Perv(X,\Lambda)}(K)$.  Let ${\cal T}_0$ denote convolution tensor subcategory  of $D_c^b(X_0,\Lambda)$  generated by the perverse sheaf $K_0 \oplus K_0^\vee$. Since
$K$ is $H$-regular, we have ${\cal T}_0 \subseteq Perv(X_0,\Lambda)$. 
The categorial dimension $d_i$ of each irreducible object $K_i$ is nonnegative, 
so the same holds for all irreducible objects in ${\cal T}_0$.  
hence by a theorem of Deligne, ${\cal T}_0$ is a Tannakian category
over the algebraically closed field $\Lambda$. Let $G$ denote its Tannaka group.

\medskip
The indecomposable elements in $Perv(X_0,\Lambda)$ are of the form 
$S_0(\gamma) \otimes_\Lambda E_n$, where $E_n$ is a $\Lambda$-vectorspace of dimension $n$
on which $Fr_\kappa$ acts by a nilpotent matrix with one Jordan block ([BBD], p.139) and where $S_0$ is an absolute simple perverse sheaf with determinant of finite order and $\gamma\in \Lambda^*$ defines a generalized Tate twist. Since $K_{0,i}\cong S_{0,i}(\gamma_i)$ for some
$\gamma_i\in \Lambda^*$,  we can a priori assume $\gamma_i=1$.
Then $$  S_{0,i} * S_{0,j}  \cong \bigoplus_k\bigoplus_l  \ S_{0,k} (\beta_{ij}^{kl})\otimes_\Lambda E_{n(i,j,k,l)} $$ follows
from $K_{0,i}*K_{0,j} = \sum_{k} c_{ij}^k \cdot K_{0,k}$. Similarly $K_{i}^\vee \cong
K_{j}$ (for some $j$ depending on $i$) gives
$ S_{0,i}^\vee \cong  S_{0,j}(\beta_i)$
for certain $\beta_{ij}^{kl}\in \Lambda^*_{mot}$ and $\beta_i\in \Lambda^*_{mot}$. 
Replacing $\kappa$ by a suitable finite base field extension, we get rid of all Jordan blocks
and we can achieve $\beta_i = 1$, since $\beta_i$ must be a root of unity.
So let us assume this.

\medskip
The Tate twists $\delta_0(\beta_{ij}^{kl})$ of the unit perverse sheaf $\delta_0$ generate a tensor category ${\cal T}_{Weil}$ contained in ${\cal T}_0$. Since
$ L * \delta_0(\beta) = L(\beta) $ and $\delta_0(\beta)*\delta_0(\beta')=\delta_0(\beta\beta')$,
all the simple objects in the tensor category ${\cal T}_{Weil}$ are invertible. Hence
the Tannaka group $G_{Weil}$ of ${\cal T}_{Weil}$ is a diagonalizable commutative algebraic group over $\Lambda$. Its connected component $G_{Weil}^0$ therefore is a torus isomorphic to $(\Gm)^\nu$ for some $\nu$, so there exists an exact sequence
$$ 0 \to (\Gm)^\nu \to G_{Weil} \to \pi_0(G_{Weil}) \to 0 \ .$$
This sequence is a split exact sequence by the structure theory of diagonalizable groups so that
$G_{Weil} \cong (\Gm)^\nu \times \pi_0(G_{Weil})$. Passing to a suitable finite extension field 
of $\kappa$, we may then assume 
$$ G_{Weil} \cong (\Gm)^\nu \ .$$ 
Since  ${\cal T}_{Weil}$ is a tensor subcategory of ${\cal T}_0$ there exists a surjective group homomorphism $$ p : G \twoheadrightarrow G_{Weil} \ .$$

\medskip
Extension of scalars from $\kappa$ to its algebraic closure $k$ defines an exact faithful tensor functor $Res$ from ${\cal T}_0$ to $Rep_\Lambda(H)$. Since by construction up to isomorphism any representation of $H$ is in the image,
by Tannaka duality [DM, prop 2.21], [D2] this functor induces an injective group homomorphism
$$ r: H \hookrightarrow G \ $$
such that $Res$ is realized by the restriction functor with respect to $r$.

\medskip 
The category ${\cal T}_0$ is semisimple.
All irreducible objects in ${\cal T}_0$
are isomorphic to $S_{0,i}* \delta_0(\alpha)$ for some $i=1,..,h$ and some $\alpha\in \Lambda^*$ where $\delta_0(\alpha)$ represents an object in ${\cal T}_{Weil}$.
A representation of $G$ becomes trivial on $H$ if and only if it is contained in the subcategory ${\cal T}_{Weil}$.  

\medskip
Since $H$ is perfect, $H \subseteq F:=Kern(p)$.
By the structure theory of reductive groups there exists a torus or rank $\nu$ in $G^0$ that
surjects onto $(\mathbb G_m)^\nu$ under $p$.  Notice $F$ commutes with
$T$, since $F$ commutes with $p(T)$. Hence there exists 
a surjective central isogeny $\varphi: F \times T \to G$. Its kernel is $$Z:= \{ (x,x^{-1}) \in F\times T \ \vert \ x\in F\cap T\}\ .$$ An irreducible representations $\pi \boxtimes \chi$ of $F \times T$ 
comes from $G$ iff the central character $\omega_\pi$ of $\pi$ is trivial on $Z$. Since $G$ has only finitely  many irreducible representations up to character twist, this implies that $F$ is a finite group.
Since the subgroup $H$ of $F$ is perfect, its intersection with $Z$ is trivial and there
is a canonical injection $i: H \hookrightarrow F/Z$. By direct inspection of the irreducible objects, $H$-regularity implies that the induced restriction map
$Rep_\Lambda(F/Z) \to Rep_\Lambda(H)$ is a bijection. Hence $H$ is isomorphic to $F/Z$ and  $\varphi: H \times T \to G$ induces an isomorphism.
By the way we defined ${\cal T}_{Weil}$, the characters
$\chi \boxtimes 1$ of $H \times T \cong G$ are just the irreducible objects of ${\cal T}_{Weil}$. 
This immediately implies that $p: T \to (\mathbb G_m)^\nu$ is an isomorphism. Thus
$$   G =  H \times (\mathbb G_m)^\nu \ $$
so that $p$ is the projection onto the second factor and $r$ is the inclusion of the first factor.   
The projection onto the second factor realizes $Rep_\Lambda(H)$ as a tensor subcategory 
of ${\cal T}_0$.  Its irreducible objects are $S_{0,i}(\alpha_i)$ for $i=1,..,h$ and certain
$\alpha_i\in \Lambda^*_{mot}$ of weight zero. 
In particular, $K_0' $ is weakly motivic.
This proves our claim except for the last
assertion.

\medskip
For the latter observe
$End_{Perv(X,\Lambda)}(K') = \bigoplus_{i=1}^{h'} End_{\Lambda}(\Lambda^{d_i}) \otimes_{\Lambda} 
End_{\Lambda} (S_{0,i}(\alpha_i))$ and that Frobenius $F_\kappa$ acts trivially on $End_\Lambda(\Lambda^{d_i})$ and trivially on $End_\Lambda(S_{0,i}(\alpha_i)) = End_\Lambda(S_{0,i})$,
since $S_{0,i}$ were absolutely simple perverse sheaves.
\qed

\medskip
Lemma \ref{descend} contradicts proposition \ref{W} that is later proved in section \ref{congruences}. This will complete the argument for the induction step in the induction used for the proof of  theorem
\ref{complexcase}.

\bigskip\noindent

\goodbreak

\section{The Cebotarev density theorem} \label{CEB}

\bigskip\noindent
Suppose $X_0$ is a variety defined over  a finite field $\kappa$. Let $k$ denote the algebraic closure of $\kappa$ and let $X$ be the extension of $X_0$ to $k$. On $X_0$ the geometric Frobenius endomorphism  $Fr=Fr_\kappa$ acts by the $q$-power map on coordinates, for $q=\# \kappa$.
Let $\kappa_m$ be an extension field of $\kappa$
of degree $m$. Attached
 to a Weil complex $K$ there is the
function
$$   f_m^K: \ X(\kappa_m)  \longrightarrow   \Lambda  \ $$ 
defined by the supertraces $$f_m^K(x) = \tau Tr(Fr_x^{\frac{m}{d(x)}};K_{\overline{x}})
= \tau Tr(Fr_\kappa^m;K_{\overline{x}}) \ , $$ where $Fr_x\in Gal(\overline{\kappa}/\kappa_m(x))$ is the geometric
Frobenius at the closed point $x$ acting on the stalk $K_{\overline{x}}$ of $K$
of a geometric point ${\overline{x}}$ over $x$. 
If we replace $\kappa$ by $\kappa_m$, we may assume $m=1$ and then
simply write $f^K_1(x)=f^K(x)$.

\bigskip
Let $K,K'$ be pure perverse Weil sheaves of weight $w$ on a variety $X$ over $k$, where $k$ is the algebraic closure of a finite field $\kappa$ with $q$ elements.
The Frobenius $Fr =Fr_\kappa $ acts on the cohomology groups $H^\bullet(X,D(K) \otimes^L  K')$. 
There exist complexes $RHom(K',K)$ of abelian groups
and $R{\cal H}om(K',K)$ of sheaves such that 
$   R\Gamma(X, R{\cal H}om(K',K)) \ = \ RHom(K',K) $ and 
$ Hom_{D_c^b(X,\Lambda)}(K',K) \ = \ H^0(RHom(K',K)) $ holds. %(e.g. [FK, p.300ff]).
Since by definition $ Hom_{Perv(X)}(K',K) = Hom_{D_c^b(X,\Lambda)}(K',K)$, we get
$$ Hom_{Perv(X)}(K',K) = H^0(R\Gamma(X, R{\cal H}om(K',K))) = H^0(X,R{\cal H}om(K',K)) \ .$$
Recall $R{\cal H}om(K',K)) \in D^{\geq 0}(X,\Lambda)$ from [KW, lemma 4.3], since $K$ and $K'$ are perverse sheaves.
On the other hand $R{\cal H}om(K',K) \cong D(K' \otimes^L D(K))$, so by Poincare duality 
$$ Hom_{Perv(X)}(K',K) = H^0_c(X, K' \otimes^L D(K))^* \ .$$
The support conditions for the perverse sheaves $K'$ and $D(K)$,
to be in  ${}^p D^{\leq 0}(X)$, imply that the cohomology
groups $H^i_c(X, K' \otimes^L D(K))$ vanish for $i>0$. 
Since by assumption
$K$ and $K'$ are pure of weight $w$, we furthermore get $ w(K' \otimes^L D(K)) \leq 0$. Hence
$$    H^\bullet_c(X, K' \otimes^L D(K)) $$
has weights $\leq 0$, and the weight $0$ only occurs as eigenvalue
of the zero-th cohomology group 
$H^0_c(X, K' \otimes^L D(K))^*$.

\bigskip
By twisting the perverse sheaves $K$ and $K'$ we may now suppose $w=0$. By the dictionary [KW, theorem III.12.1 (6)] and
the Grothendieck-Lefschetz trace formula,  the characteristic function $\tau trace(Fr^m, H^\bullet_c(X, K' \otimes^L D(K)))$
is
$$  (f^{K'}_m,f^K_m) := \sum_{x\in X(\kappa_m)} \ f^{K'}_m(x) \overline{ f}^K_m(x) = 
\sum_{x\in X(\kappa_m)} \ f^{K' \otimes^L DK}_m(x) \ .$$
Hence by the vanishing result from above
$$ (f^{K'}_m,f^K_m) = \tau Tr(Fr^m; Hom_{Perv(X)}(K',K)) + \sum_{\nu < 0} (-1)^\nu \tau Tr(Fr^m; H^\nu_c(X, K' \otimes^L D(K)))\ .$$
The eigenvalues of $Fr^m$ on $H^\nu_c(X, K' \otimes^L D(K)))$ are
$ \leq q^{\nu m/2}$, hence the sum on the right side can be estimated
by $C \cdot q^{-m/2}$ for a constant $C = \sum_{\nu <0} dim_\Lambda
H^\nu_c(X,K' \otimes^L D(K))$
 depending only on $K$ and $K'$ but not on $m$. In other words

\begin{Corollary} \label{lastcor} If $K$ and $K'$ are pure perverse Weil sheaves of weight $w=0$, then
$$ \vert (f^{K'}_m,f^K_m) - \tau Tr(Fr^m; Hom_{Perv(X)}(K',K)) \vert  \ \leq\
 C \cdot q^{-m/2} \ .$$
 For any $\varepsilon >0$ there exists an integer $m_0=m(\varepsilon, K,K')$ such that  
$$ dim_\Lambda(Hom_{Perv(X)}(K,K')) \ > 0 $$
holds provided $(f^{K'}_m,f^K_m) > \varepsilon$ holds for some $m\geq m_0$.
\end{Corollary}

\bigskip
{\it Proof}. Indeed $\big\vert \tau Tr(Fr^m;Hom_{Perv(X)}(K,K')) - (f^{K'}_m,f^K_m) \big\vert < C q^{-m/2}$. \qed

\bigskip\noindent

\begin{Corollary}\label{limit}
For an irreducible perverse sheaf of weight zero defined over a finite field $\kappa$ 
for $\Vert K \Vert^2_m := (f_m^K,f_m^K)$ we get $$\lim_{m\to \infty} \Vert K\Vert^2_m = \lim_{m\to \infty} \sum_{x\in X(\kappa_m)}\vert f^K_m(x) \vert^2 = 1\ .$$
\end{Corollary}

\bigskip
{\bf Remark}. For a finite dimensional $\Lambda$-vectorspace $V$ with a continuous action of $Fr$  the traces $Tr(Fr^m,V) = \sum_{i=1}^{\dim(V)} \alpha_i^m$ are given by the eigenvalues
$\alpha_1,..,\alpha_{\dim(V)}$ of $Fr$ on $V$. If $w(V)=0$, then by the Kronecker-Weyl approximation theorem 
applied to $\log(\tau \alpha_i)/2\pi i\log(q) \in \R/\Z$ there exist infinitely many integers $m$ such that $\vert \tau Tr(Fr^m,V) - \dim_{\Lambda}(V) \vert < \varepsilon$ holds for fixed $\varepsilon >0$. 

\bigskip
{\it Two easy consequences}.
For a smooth $\Lambda$-adic sheaf $E_0\neq 0$ on $X_0$ and arbitrary fixed $x\in X(\kappa)$ the values $f_m^K(x)$ for $K_0=E_0[\dim(X)]$ are nonzero for infinitely many integers $m$. Secondly, suppose $K$ has weight zero and $E$ decomposes into inequivalent irreducible smooth $\Lambda$-adic Weil sheaves $E_i$ with multiplicity $m_i$. Then there exist infinitely many integers $m$ such that $\vert \tau Tr(Fr^m, End(K)) - \sum_i m_i^2 \vert < \varepsilon$ holds for some given $\varepsilon >0$. One immediately reduces to the case $E= V \otimes_\Lambda  E_1$ with $\dim_\Lambda(V)=m_1$.
Since $ Tr(Fr^m;Hom_{Perv(X)}(K_1,K_1))=1$, the claim now follows from $Tr(Fr^m, End(K)) = Tr(Fr^m,End_\Lambda(V)) \cdot Tr(Fr^m;Hom_{Perv(X)}(K_1,K_1))$. 

\goodbreak

\section{Fourier transform}\label{Fourier}

\bigskip\noindent
Suppose $X_0$ is an abelian variety defined over  a finite field $\kappa$. Let $a_0:X_0\times X_0 \to X_0$ be the group law. Let $k$ denote the algebraic closure of $\kappa$ and let $X$ be the extension of $X_0$ to $k$. On $X_0$ the geometric Frobenius {\it endomorphism}  $Fr=Fr_\kappa$ acts by the $q$-power map on coordinates for $q=\# \kappa$.

\bigskip
{\it Lang torsors}. The Lang torsor for $X_0$ is
defined by the etale  homomorphism ${\wp}_0(x)= Fr(x) - x$
$$  {\wp}_0: X_0 \to X_0 \ .$$
The morphism ${\wp}_0$ defines a finite etale geometrically irreducible Galois covering of $X_0$, whose
Galois group by [S2, p.116] is  
the abelian finite group $X_0(\kappa)$, i.e. the kernel of ${\wp}_0$.   
Considered over the algebraic closure $k$ of $\kappa$, this defines an etale
covering  with Galois group $\Delta$ denoted ${\wp}: X \to X$
$$ \xymatrix{ 0 \ar[r] & \Delta \ar[r] & X \ar[r]^{\wp} & X \ar[r] & 0 }    \ .$$
Since $\wp_0$ is geometrically irreducible, we get $\Delta = X_0(\kappa)$.
Let $\Delta^*$ denote the group of  characters of $\Delta$ with values in $\Lambda^*$ 
$$ \chi: \Delta \longrightarrow  \Lambda^* \ .$$
The direct image ${\wp}_{0,*}(\Lambda_{X_0}) $ of the constant sheaf $\Lambda_{X_0}$ on $X_0$ decomposes into a direct sum $\bigoplus_{\chi \in \Delta^*} L_{\chi,0}$  of smooth rank one  $\Lambda$-adic sheaves $L_{\chi,0}$ on $X_0$. Let denote $L_\chi$ the scalar extension of $L_{\chi,0}$ to $k$, defining the corresponding smooth etale Weil sheaf on $X$.  By class field theory, see [S2, p.142], 
$$ \fbox{$ Tr(Fr_x^m, L_\chi) = \chi(x)^{-m} $}  $$
holds for all points $x\in X(\kappa)$. Hence the functions $$f^{L_\chi}(x)=\chi(x)^{-1}$$
separate points in $X(\kappa)$. 

\bigskip
Let $m\geq 1$ be an integer. Again, $Fr^m (x) - x$  defines a geometrically irreducible etale morphism $X_0 \to X_0$, which however becomes a Galois covering only after a base field extension by passing to the finite extension field $\kappa_m$ of $\kappa$ of degree $m$, where
it induces the Lang torsor $$\wp^{(m)}_0: X_0\times_{Spec(\kappa)} Spec(\kappa_m) \to X_0\times_{Spec(\kappa)} Spec(\kappa_m)$$ of the abelian variety $X_0\times_{Spec(\kappa)} Spec(\kappa_m)$ over $Spec(\kappa_m)$.
The Frobenius {\it automorphism} $F=F_\kappa$ acts on $X_0\times_{Spec(\kappa)} Spec(\kappa_m)$ via its Galois action on $\kappa_m$, inducing an action of $F$ on $X_0(\kappa_m) = Hom_{Spec(\kappa_m)}(Spec(\kappa_m),
X_0\times_{Spec(\kappa)} Spec(\kappa_m))$. This action of $F$ coincides with the action of the Frobenius endomorphism $Fr$ on $X_0(\kappa_m)$. 
$$ \xymatrix@-0,5cm{   &    X_0\times_{Spec(\kappa)} Spec(\kappa_m) \ar[dd]^{\wp_0^{(m)}} \ar[dl]\cr
X_0 \ar[dd]_{Fr^m - id} &   \cr
&   X_0 \times_{Spec(\kappa)} Spec(\kappa_m) \ar[dd] \ar[dl]\cr
X_0 \ar[dd] &  \cr
&  Spec(\kappa_m) \ar[dl]\cr
Spec(\kappa) &    } $$

By $F^m-1 = (F-1) \sum_{i=0}^{m-1} F^i$  the trace $\sum_{i=0}^{m-1} F^i$ defines a homomorphism 
$$   S_m=\sum_{i=0}^{m-1} F^i:  \ X_0(\kappa_m) \longrightarrow X_0(\kappa) \ ,$$
which is surjective by [S2,VI,\S 1.6].    
Any character $\chi: X(\kappa) \to \Lambda^*$ can be 
extended to a character 
$$  \chi_m:  X(\kappa_m) \longrightarrow \Lambda^* \ ,$$
where $\chi_m = \chi\circ S_m$ is defined by the composite of the trace $S_m:  X(\kappa_m) \to X(\kappa)$ 
and the character $\chi$. Then, by definition, $\chi_m(x) = \chi(x^m) = \chi(x)^m$ holds for
$x\in X(\kappa)$. 

\medskip
More generally $S_{r,rm}: X(\kappa_r)^* \to X(\kappa_{rm})^*$,
defined by $$\chi_r \mapsto \chi_r \circ S_{r,rm}= \sum_{i=0}^{m-1} F^r\ ,$$ is injective.
Using these transition maps, any character $\chi_r \in X(\kappa_r)^*$ defines a collection
of characters $(\chi_{rm})_{m\geq 1}$ such that $\chi_{rm} \in X(\kappa_{rm})^*$. 
Any such $(\chi_{rm})_{m\geq 1}$ defines a translation invariant sheaf, say $L_\psi$, on $X$.
In this sense, we view $\psi$ as a character on all $X_0(\kappa_{rm})$ for $m=1,2,..$. For large enough $n$ all $L_\psi$, for torsion characters $\chi$ of $\pi_1(X,0) \to \Lambda^*$, arise in this way. Indeed, for $\chi^k=1$ chose $r$ large enough so that $X[k] \subset X(\kappa_r)$. Then there exist etale isogenies $f: X'\to X$ and $g:X\to X'$ such that $f_s \circ g = \wp^{(r)}$, where 
$f_s: X'\to X$ is the etale part of the isogeny $k \cdot id_X = f_s\circ f_{ins}$ and $f_{ins}$ is the 
inseparable part.  Since $\Lambda_{X'}$ is a direct summand
of $g_*(\Lambda_X)$, the sheaf $L_\psi$, as a constituent of $f_{s*}(\Lambda_{X'})$,  is also
a direct summand of $\wp^{(r)}_*(\Lambda_X)$.

\bigskip
The Frobenius automorphism $F$ acts on each $X_0(\kappa_n)$ and hence on the characters $\psi: X_0(\kappa_n) \to \Lambda^*$; and of course $F^n$ acts trivially. Conversely, suppose $\psi^{F^r}=\psi$ or $\psi((F^r - id)(x))=1$ holds for all $x\in X_0(\kappa_n)$ and all $n$. We may enlarge $n$ and hence assume that $r$ divides $n$, hence $n=rm$. We may replace $\kappa$ by $\kappa_r$
and $F^r$ by $F_{\kappa_r}$. Then $\psi$ factorizes over the quotient $X_0(\kappa_n)/(1-F_{\kappa_r})X_0(\kappa_n)$, which is isomorphic to $X_0(\kappa_r)$ via the trace homomorphism. Hence $\psi$ comes from a character $\psi': X_0(\kappa_r) \to \Lambda^*$ by the trace extension
$\psi = \psi'_{n/r}$, defined above. 

\bigskip {\it Induction}.
For a finite field extension $\kappa_r$ of $\kappa$ and a character $\psi: X_0(\kappa_r) \to
\Lambda^*$ we consider the Weil sheaf $L_\psi$ on $X$. It is $F_{\kappa_r}$-equivariant, and there exists an isomorphism $F_{\kappa_r}^*(L_\psi) \cong L_\psi$ for the Frobenius automorphism $F_{\kappa_r}$ of the field $\kappa_r$. Let $L$ be an $F_{\kappa_r}$-equivariant
Weil sheaf on $X$. Since $F_{\kappa_r}=F^r$ holds for the Frobenius automorphism $F$ for the field $\kappa$, 
$K=\bigoplus_{i=0}^{r-1} (F^i)^*(L_\psi)$ defines an $F$-equivariant Weil sheaf on $X$ over $\kappa$, i.e. one has an isomorphism $F^*(K) \cong K$ and we write $K=Ind_{\kappa_r}^\kappa(L)$. Any $F$-equivariant Weil sheaf $K$ on $X$, which as a perverse sheaf on $X$ is translation invariant under $X$ and multiplicity free, is of the form $K \cong \bigoplus_\psi \ L_\psi(\alpha_\psi)$ for certain characters $\psi : X_0(\kappa_{r(\psi)}) \to \Lambda^*$ and certain Tate twists by $\alpha_\psi\in \Lambda^*$; we can assume that the integers $r(\psi)$ are chosen to be minimal. Then $F^*(K)\cong K$ implies  $$K \cong \bigoplus_{\varphi} \ Ind_{\kappa_r}^\kappa(L_\varphi(\alpha_\varphi))\ $$  for certain
 $\varphi\in \{\psi \}$,  so that $F^i(L_\varphi)\cong L_\varphi$ holds if and only if $r(\varphi)$ divides $i$. 

\bigskip
{\it Character Twists}. For perverse sheaves $K_0$ on $X_0$ and $\chi_0 \in \Delta^*$ we can define the twisted complex $K_0\otimes L_{0,\chi}$. Again this is a perverse sheaf on $X_0$. 
Let $K\otimes \chi$ or $K_\chi$ denote the corresponding perverse Weil complex on $X$. Its 
associated function (over $\kappa_m$) is 
$$  f_m^{K\otimes \chi}(x) \ = \ f_m^K(x) \cdot \chi_m(x)^{-1} \ .$$

\bigskip
{\it Fourier transform}. By varying the characters $\chi: X(\kappa_m) \to \Lambda^*$, we define the Fourier transform $\widehat f^K_m$ of $f_m^K$ for each $m$ by the following summation
$$  \widehat f^K_m: X(\kappa)^* \to \Lambda \ .$$
$$ \fbox{$  \widehat f^K_m(\chi) = \sum_{x\in X(\kappa_m)}  f^K_m(x) \cdot \chi_m(-x) $} \ .$$
Fourier transform is additive with respect to the perverse sheaf $K$ in the sense that  $$ \widehat f^{K\oplus L}_m(\chi) \ = \ \widehat f^{K}_m(\chi)  \ + \
\widehat f^{L}_m(\chi) \ .$$

\smallskip
{\it Example 1}. For the skyscraper sheaf $K_0=\delta_0$ concentrated at the origin
%one has $f_m^K(x) =1$ for $x=0$, and $f_m^K(x) =0$ otherwise.
the Fourier transform is constant, i.e. $ \widehat f^K_m(\chi)= 1$ holds for all $m$ and all $\chi$. 

\bigskip
{\it Example 2}. The Fourier transform $ \widehat f^{K}_m(\chi) $ for the character sheaf $K_0=L_{\varphi}$  of a character $\varphi: X_0(\kappa_m)\to \Lambda$  is $\# X(\kappa_m)$ at $\chi=\varphi^{-1}$, and it is zero otherwise. The following perverse sheaf $K=\delta_X^\varphi $ is of weight $w=0$   $$\delta_X^\varphi :=  L_\varphi[\dim(X)](\alpha) \quad , \quad \alpha = q^{-\dim(X)/2} \ .$$ 
Its  Fourier transform $\widehat f^K(\chi)$ is $(-1)^{\dim(X)} q^{-m\dim(X)/2}\# X(\kappa_m)$ for $\chi=\varphi^{-1}$, and it is zero otherwise. 
 
\bigskip
{\it Example 3}. Given $\psi: X_0(\kappa_r)\to \Lambda^*$, for an extension $\kappa_r$ of $\kappa$ of degree $r$ suppose for the Frobenius automorphism $F$ of $\kappa$ that $F^i(L_\psi)\cong L_\psi$ holds iff $r$ divides $i$. Then $K= Ind_{\kappa_r}^\kappa(L_\psi)$ is defined over $\kappa$.  If $r$ divides $m$, then $f^K_m(x)= \sum_{i=0}^{r-1} f^{F^i(L_\psi)}_m(x)$  and hence $f^K_m$ can be computed by example 2); otherwise
$$  f_m^K(x) = 0 \quad \mbox{ for } \quad r \not\vert \ m \ .$$  
By base extension this last assertion can be easily reduced to the case $m=1$. For $m=1$ then $f^K(x)= Tr(Fr_x; K_{\overline x})$ holds for $x\in X_0(\kappa)$. Since the substitution $Fr_x$
permutes the summands of  $K_{\overline x} = \bigoplus_{i=0}^{r-1} F^i(L_\psi)_{\overline x}$
in the sense that $Fr_x: F^i(L_\psi)_{\overline x} \to F^{i+1}(L_\psi)_{\overline x}$ for $i<r-1$ and
$Fr_x: F^{r-1}(L_\psi)_{\overline x} \to F^{r}(L_\psi)_{\overline x} \cong (L_\psi)_{\overline x}$, the trace of $Fr_x$ is zero unless $r=1$.
 
\bigskip 
 By the Grothendieck-Lefschetz trace formula the Fourier transform 
has the following interpretation \label{GLF}
$$  \widehat f^{K}_m(\chi) \ =\ 
\sum_{x\in X(\kappa_m)}  f^{K\otimes \chi}_m(x) \ = \ 
f_m^{Rp_*(K\otimes \chi)}(*) \ $$
for the structure map $p_0: X_0 \to Spec(\kappa)=\{*\}$. In other words,
$f_m^{Rp_*(K\otimes \chi)}(*)$ is the trace of $Fr^m$ on the etale cohomology group
$H^\bullet(X,K_{\chi})$. 

\bigskip
{\it Convolution}. For $\Lambda$-adic Weil complexes $K$ and $L$ on $X$ the convolution $K*L$
is the direct image complex $K*L=Ra_*(K\boxtimes L)$ for the group law $a:X\times X\to X$. Hence by the Grothendieck-Lefschetz trace formula 
 $$f_m^{K*L}(x)= \sum_{y \in X(\kappa_m)} f_m^K(x-y)f_m^L(y)$$ which is the usual convolution of the functions $f_m^K(x)$ and $f_m^L(x)$ on the finite abelian group $X(\kappa_m)$.
By elementary Fourier theory therefore
$$    \widehat f_m^{K*L}(\chi) =  \widehat f_m^{K}(\chi) \cdot \widehat f_m^{L}(\chi) $$
For perverse sheaves $K$ that are pure of  weight $w=0$ we furthermore have 
$$     \widehat f_m^{K^\vee}(\chi) =   \overline{\widehat f_m^{K}(\chi)} \ $$
by [KW], theorem III.12.1(6). Indeed 
$$ \widehat f_m^{K^\vee}(\chi) = \sum_x f^{D(K)}(-x)\chi_m(-x)=
\sum_x \overline{f^{K}(-x)}\chi_m(-x)= \sum_x \overline{f^{K}(x)}\overline{\chi_m(-x)}=
\overline{\widehat f_m^{K}(\chi)}\ .$$

\bigskip\noindent
For $y\in X(\kappa_m)$ one has the elementary {\it Fourier inversion formula}
$$  f_m^K(y) = \frac{1}{\# X(\kappa_m)} \sum_{\chi \in X(\kappa_m)^*}
\widehat f_m^K(\chi) \chi_m(y) \ $$
with summation over all characters $\chi: X(\kappa_m) \to \Lambda^*$
(i.e. the characters obtained from the Lang torsor). 
For a pure perverse sheaf $K$ of weight zero, in the limit $m\to \infty$  the sum 
$$  \Vert K \Vert_{X_m}^2 = \sum_{x\in X(\kappa_m)} \vert f_m^K(x) \vert^2  = (f^K_m,f^K_m)$$
converges to $Tr(Fr^m;End_{Perv(X)}(K))$ as shown in section \ref{CEB}.

\bigskip\noindent

\goodbreak

\section{The Plancherel formula}\label{Plancherel}

\bigskip\noindent
Let $X$ be an abelian variety $X$ over the algebraic closure $k$
of a finite field $\kappa$. Let $K$ be a pure perverse Weil sheaf of weight 0 on $X$, so $K$ is equivariant $F^*(K)\cong K$ for the Frobenius automorphism
$F=F_{\kappa}$.   
The elementary Plancherel formula expresses $ \Vert K \Vert_{X_m}^2 = \Vert f_m^K\Vert^2$ in terms of the Fourier transform $\widehat f^K_m(\chi)$ of $f^K_m(x)$
$$  \Vert K \Vert_{X_m}^2 \ = \ \frac{1}{\# X(\kappa_m)} \cdot \Vert \widehat f_m^K \Vert^2 \ .$$
Here by definition the $L^2$-norms are $ \Vert f^K_m \Vert^2 = (f^K_m,f^K_m)$ for $(f,g)= \sum_{x\in X(\kappa_m)} f(x)\overline{g(x)}$ resp. $ \Vert \widehat f^K_m \Vert^2 = (\widehat f^K_m,\widehat f^K_m)$ for
$(\widehat f ,\widehat g )= \sum_{\chi\in X(\kappa_m)^*} \widehat f(\chi)
\overline{\widehat g}(\chi) $. More generally
$$  (f^K_m,f^L_m) \ = \ \frac{1}{\# X(\kappa_m)} \cdot  (\widehat f_m^K,\widehat f_m^L ) \ .$$

\bigskip\noindent
{\it Example}. $(f^K_m,f^L_m)$  for the Weil sheaves $K=L_\psi(\alpha)$ and $L=L_\varphi(\beta)$ vanishes unless $\psi=\varphi$, and in this case 
 $(f^K_m,f^L_m)= \alpha\cdot \overline{\beta} \cdot \# X(\kappa_m)$.  
For the following recall that the perverse Weil sheaf $\delta_X^\varphi = L_\varphi[\dim(X)](q^{-\dim(X)/2})$ is pure of weight zero.
 
\begin{Lemma} \label{14}
An irreducible perverse Weil sheaf $K$ on $X$ of weight zero, for which for almost all torsion characters $\chi$ the cohomology groups $H^\bullet(X,K_\chi)$ are zero, is isomorphic to $\delta_X^\varphi$ for some character $\varphi$. A pure perverse Weil sheaf $K$ (of weight zero) on $X$ is translation invariant on $X$, if the cohomology groups $H^\bullet(X,K_\chi)$ are zero for almost all torsion characters. 
\end{Lemma} 
 
For the proof notice that the characters $\chi: \pi_1(X,0) \to \Lambda^*$ 
that are defined over some finite extension of the finite field $\kappa$
are precisely the torsion characters of $\pi_1(X,0)$.
 
 \medskip
{\it Proof of lemma \ref{14}}. If $K$ is irreducible and $K\not\cong \delta_X^\varphi$ holds for all characters $\varphi$, then $H^i(X,K_\chi)=0$ for all $\chi$ and all $i$ with $\vert i\vert \geq \dim(X)$.  
Suppose there exist only finitely many characters $\chi_1,..,\chi_r$ for which $H^\bullet(X,K_\chi)$
does not vanish. Then $\widehat f^K_m(\chi)$ is zero  except for $\chi=\chi_i$,  $i=1,..,r$. 
Furthermore there exist constants $c_i$ independent from $m$ so that $\vert \widehat f^K_m(\chi_i) \vert \leq c_i \cdot q^{m(\dim(X)-1)/2}$ by the Weil conjectures, since these values are the traces of $Fr^m$ on $H^\bullet(X,K_\chi)$.
Hence $\Vert K\Vert_{X_m}^2 \leq \frac{const}{\# X(\kappa_m)} \cdot q^{m(\dim(X)-1)}$ by the Plancherel formula, and this implies $\Vert K\Vert_m \to 0$  in the limit $m\to \infty$ contradicting corollary \ref{limit}.
For pure $K$ we can apply this argument to the simple constituents $P$ of $K$,  to show that $K$ is translation invariant under $X$. Notice that $H^\bullet(X,K_\chi)=0$ implies $H^\bullet(X,P_\chi)=0$ by the decomposition theorem. \qed 

\begin{Corollary} \label{key}
Suppose $K\neq 0$ is an irreducible perverse Weil sheaf on $X$ and acyclic in the sense that
$H^\bullet(X,K)=0$. Then $H^\bullet(X,K_\chi)\neq 0$ holds for some torsion character $\chi: \pi_1(X,0)\to \Lambda$. 
\end{Corollary}

\medskip
The last corollary is an ingredient of the proof of the main theorem in [W2], and this main theorem
implies the generic vanishing theorem for abelian varieties over finite fields.
Using these consequences of corollary \ref{key}, one can further improve the statement
(this also holds over $\mathbb C$ if we omit the torsion restrictions).  

\medskip
\begin{Corollary} \label{strongerversion} Suppose $K$ is a complex in $D_c^b(X_0)$.
If all torsion character twists $K_\chi$ of $K$ are acyclic in the sense that
$H^\bullet(X,K_\chi)=0$, 
then $K=0$. If $X$ is simple and $K_\chi$ is acyclic for almost all torsion characters $\chi$, then all irreducible constituents of the perverse cohomology sheaves of $K$ are translation invariant. 
\end{Corollary}

\medskip
{\it Proof}. 
By the generic vanishing theorem [W2], $H^\nu(X,K_\chi) \cong
H^0(X,{}^pH^\nu(K_\chi))$ holds for generic torsion characters $\chi$. Since ${}^pH^\nu(K_\chi) = {}^pH^\nu(K)_\chi$,  vanishing
$H^\bullet(X,K_\chi)=0$ for generic torsion $\chi$ implies  
$$H^\bullet(X,{}^pH^\nu(K_\chi))= H^0(X,{}^pH^\nu(K_\chi))=0\ $$
so that for these $\chi$ the ${}^pH^\nu(K_\chi)\cong {}^pH^\nu(K)_\chi$  are acyclic perverse sheaves for all $\nu\in \mathbb Z$.
Hence by the main theorem of [W2]  all Jordan-H\"older
constituents of the perverse
sheaf $K$ are translation invariant under some nontrivial abelian subvarieties of $X$. If $X$
is simple, therefore all  ${}^pH^\nu(K)$ are extensions of translation invariant irreducible perverse sheaves.  This proves the second assertion.

\medskip
To show the first assertion, now assume that $X$ is simple.  
If $K\neq 0$, choose $m$  maximal such that
${}^pH^m(K)\neq 0$. By the right exactness of the functor $H^n(X,-)$ on perverse sheaves on $X$
for $n=\dim(X)$ there exists a surjection $H^{n+m}(X,K_\chi) \to H^n(X,{}^pH^m(K_\chi))$. 
So $H^{n+m}(X,K_\chi)=0$ implies $H^n(X,{}^pH^m(K)_\chi)=0$, and $K_\chi$ is acyclic for all torsion $\chi$. Since ${}^pH^m(K)\neq 0$, there exists a nontrivial quotient morphism
${}^pH^m(K) \to \delta_X^\psi$ by the second assertion ($X$ is simple) and $\psi$ is defined over some finite field and hence torsion. By the right exactness of $H^n(X,-)$ then
$H^n(X,\delta_X^{\psi\chi})=0$ holds for all torsion $\chi$.  A contradiction, that implies  $K=0$.

\medskip
To prove the first assertion for non-simple abelian varieties, one can assume that
$X$ is a product of simple abelian varieties. By considering direct images to the quotient
factor abelian varieties and proper basechange, our claim then can be easily reduced to the simple case by induction on the number
of simple factors.  \qed

\begin{Lemma}
Suppose $X_0=A_0 \times B_0$ is a product of two abelian varieties defined over a finite field $\kappa$. Let $K$ be a pure Weil sheaf of weight 0 on $X$ and let  $p_0:X_0\to B_0$ be the projection onto the second factor. Then  we have the following relative Plancherel formula
$$  \Vert K \Vert_{X_m}^2 \ =\ \frac{1}{\#A(\kappa_m)^*} \sum_{\chi\in A(\kappa_m)^*}
\Vert Rp_*(K_\chi) \Vert^2_{B_m} \      $$
where the summation is over all characters $\chi: A(\kappa_m) \to \Lambda^*$.
\end{Lemma}

\bigskip\noindent
{\it Proof}. $\Vert f^K_m \Vert^2_{X(\kappa_m)}  = \sum_{(x,y)\in X(\kappa_m)} 
\vert f^K_m(x,y) \vert^2 = \sum_{y\in B(\kappa_m)} \bigl(\sum_{x\in A(\kappa_m)}
\vert f^K_m(x,y)\vert^2\bigr) $ by definition. Using the Plancherel
formula for $A$ the inner sum can be rewritten  so that we get $\sum_{y\in B(\kappa_m)}
\#A(\kappa_m)^{-1} \sum_{\chi\in A(\kappa_m)^*}\vert \sum_{x\in A(\kappa_m)} f^K_m(x,y)\chi(-x)\vert^2 $. For $L=Rp_*(K_\chi)$ the Grothendieck-Lefschetz
trace formula gives $\sum_{x\in A(\kappa_m)} f^K_m(x,y)\chi(-x) = f^L_m(y)$. Hence
$$ \Vert f^K_m \Vert^2_{X(\kappa_m)}   =
\frac{1}{\#A(\kappa_m)} 
\sum_{\chi\in A(\kappa_m)^*}\sum_{y\in B(\kappa_m)} \vert f_m^{Rp_*(K_\chi)}(y) \vert^2 =  \frac{1}{\#A(\kappa_m)} 
\sum_{\chi_1\in A(\kappa_m)^*} \Vert f_m^{Rp_*(K_\chi)} \Vert^2_{B(\kappa_m)}
   \ .$$ \qed

\goodbreak

\section{\bf Congruences and primary decomposition}\label{congruences}

\bigskip\noindent

\bigskip\noindent
Let $k$ be the algebraic closure of a finite field $\kappa$ of characteristic $p$. Recall that a perverse semisimple sheaf $K$ on a variety $X$  over $k$
is a direct sum of simple perverse sheaves $L$. The support $Y=Y(L)$ of each 
constituent $L$ is an irreducible subvariety $Y$ of $X$. On a Zariski open dense  subvariety $U\subseteq Y$ the restriction $L\vert U$ of $L$ to $U$ is isomorphic to $E[\dim(U)]$ for some smooth
$\Lambda$-adic local system $E$ on $U$.

\begin{Proposition} \label{W} Let $X$ be an abelian variety defined over $\kappa$.
Let $K$ be a weakly motivic perverse semisimple Weil sheaf on $X$ 
over $\kappa$.  
If there exists an isomorphism $K*K \cong d\cdot K$ 
of Weil sheaves over $\kappa$, then $K$ is a skyscraper sheaf
on $X$.
\end{Proposition}

\bigskip\noindent
{\it Proof}. The condition $K*K \cong d\cdot K$ implies
$f^{K*K}_m = f^K_m * f_m^K = d \cdot f_m^K$, hence
$$    \widehat f_m(\chi)^2 \ =\   d \cdot \widehat f_m(\chi) $$
for all characters $\chi\in X_m^*$. Thus the Fourier transform
$\widehat f_m$ of $f_m$ is constant and equal to $d$ on its
support. 

\medskip
If ${\frak o}$ denotes the ring of integers ${\frak o}$ in $\Lambda  = \overline\Q_l$, for a weakly motivic complex
$K$ the functions $f_m^K$ have their values in ${\frak o}[\frac{1}{p}] \subset \Lambda$ for all $m$.
For $l\neq p$ hence $f_m^K(x)$ has $l$-adic integral
values and therefore the same  holds for the Fourier transform $\widehat f_m^K$.
For two characters $\chi', \chi\in X_m^*$, now assume
$$  \chi' = \chi \chi_l \quad , \quad (\chi_l)^{l^r} = 1 \ .$$
Since $\chi_l \in X_m^*$ is a character of $l$-power order, its
values are $l^r$-th roots of unity and hence congruent to 1
in the residue field $\mathbb F$ of ${\frak o}[\zeta_{l^r}]$ modulo each prime ideal of ${\frak o}[\zeta_{l^r}]$ over $l$. Hence
we get the congruence
$$   \widehat f_m^K(\chi') \ \equiv \ \widehat f^K_m(\chi) \quad ( \mbox { in } \mathbb F\ ) \ .$$
Indeed 
$  \widehat f_m^K(\chi') = \sum_{x\in X_m} f_m^K(x) \chi'(-x) \equiv 
\sum_{x\in X_m} f_m^K(x) \chi(-x) = \widehat f_m^K(\chi)$, since 
$\chi'(-x) = \chi(-x)\chi_l(-x) \equiv \chi(-x)$ holds in $\mathbb F$.
Now suppose that the prime $l$ is larger than $d$ and different from $p$.
Then the congruence $   \widehat f_m^K(\chi')  \equiv  \widehat f^K_m(\chi) $
implies
$$   \widehat f_m^K(\chi') \  = \ \widehat f^K_m(\chi) $$
since $\widehat f_m^K$ was equal to the integer $d$ on its support.

\medskip
By theorem \ref{Drinfeld} below, for every prime $\tilde l \neq p$ there always exists a $\overline\Q_{\tilde l}$-adic semisimple perverse sheaf $\tilde K$ so that $f^{\tilde K}_m = f^K_m$ holds for some underlying isomorphisms $\tilde\tau: \overline\Q_{\tilde l} \cong \C$ and $\tau: \overline\Q_l \cong \C$.  So, replacing $l$ by some $\tilde l$ larger than $d$ and $p$, we see that
$\widehat f^K_m$ (which only depends on $K$ and $m$ and not on $\tilde l$)
is translation invariant on the $\tilde l$-primary component $X_m^*[\tilde l^\infty]$ of 
the finite group $X_m^*$. Actually, since this holds for all primes $l$ not dividing
$N=p\cdot d$, the function $\widehat f_m^K$ is constant on the $l$-primary subgroups for
$l\!\not\vert \, \, N$ and for all $m$
$$ \xymatrix{  X(\kappa_m)^* \ \ar[rr]^-{\widehat f^K_m} \ar[dr]  & & \ E \cr
& \ \bigoplus_{l\vert p\cdot d} X(\kappa_m)^*[l^\infty]\ \ar@{.>}[ur]^\exists & \cr} $$
Therefore $f_m^K$ as the Fourier inverse of $\widehat f_m^K$ for all $m$  must be a function with support
contained in the $N$-primary subgroup of $X_m$ 
$$   supp(f_m^K) \ \subseteq \  X_m[N^\infty] \ .$$
By the next lemma \ref{Jacobi}
this implies that $K$ is a skyscraper sheaf. Indeed, let $Y$ be the support of 
an irreducible perverse constituent $L$ of $K$. Then we claim $\dim(Y)=0$. 
For this gather the constituents $L$ with some fixed support $Y$, chosen so that $\dim(Y)$ is maximal. Then a Zariski open dense subset $U\subset Y$ is disjoint to the support
of all irreducible components with support different from $Y$, and the restriction
of the direct sum of irreducible components $L$ of $K$ to $U$
becomes isomorphic to $E[\dim(U)]$ for some smooth etale sheaf $E$ on $U$ 
of rank say $r$.  Assuming $E\neq 0$, 
then for any point $x\in U(k)$ there exists an integer $m$ such that $x\in X_m$.
Replacing $m$ by $mn$ for a suitable integer $n$, we get $f_{mn}^K(x) = (-1)^{\dim(U)} f^E_{mn}(x)
\neq 0$ since otherwise the stalk $K_{\overline x}$ would be zero contradicting
the smoothness of $E\neq 0$ at $x$. This argument and the previous support property implies 
$$U(k) \ = \ \bigcup_{m=1}^\infty \ U_m \ \subseteq  \ \bigcup_{m=1}^\infty supp(f_m^K) \ \subseteq \ \bigcup_{m=1}^\infty \ X_m[N^\infty] \ .$$
Hence $\dim(U)=0$ by lemma \ref{Jacobi}. In other words, $K$ is a skyscraper
sheaf on $X$. \qed

\medskip

\begin{Lemma} \label{Jacobi}
Let $U_0 \subseteq A_0$ be an open dense subset of an absolutely irreducible closed subset $Y_0$ of an abelian variety
$A_0$ over a finite field $\kappa$. Suppose there exists an integer $N$ such that
$$   U(\kappa_m) \subseteq A(\kappa_m)[N^\infty] $$
holds for all $m$. Then $\dim(Y)=0$. 
\end{Lemma}

{\it Proof}. $Y$ is projective. For $\dim(Y)>0$ there exists an irreducible curve in $U$ defined over
a finite field. If the assertion of lemma \ref{Jacobi} were false, without restriction of generality we can then assume $\dim(Y)=1$. We may replace $A$ by the abelian subvariety generated by the curve $Y$.
Let $Y_n$ be the irreducible closed subvariety $Y_n\subseteq A$ defined by the image of $a_n: Y^n \to A$
under the addition map $(y_1,..,y_n)\mapsto \sum_{i=1}^n y_i$. Then $\dim(Y_n)\leq n$ and $\dim(Y_n)=n$ holds for $n \leq d=\dim(A)$, since otherwise $\dim(Y_n)<n\leq d$ for some $n$  so that $Y_{n-1}$ is stable under $Y$, and hence stable under the simple abelian variety $A$ generated by $Y$ although $\dim(Y_{n-1}) < d$. This shows $Y_{d} =A$ and $Y_{d-1} \neq A$ and $a_d(U^d) + \sum_{s\in Y\setminus U} (s+ Y_{d-1}) = Y_d = A$. Hence the complement of  $a_d(U^d)$ is contained in finitely many translates of $Y_{d-1}$. So
$a_d(U^d)$ contains a nonempty Zariski open subset $V$. After replacing $V$ by $ V\cap -V$ we can assume $V=-V$. Then as a set the union $V \cup (V+V)  \cup (V+V+V) \cup ...$  is $A$, since otherwise
there would exist a $V$-invariant proper closed nonempty subset of $A$, and that is impossible. 
As in lemma \ref{closed}, by 
the noetherian property of the Zariski topology there exists an integer $m$ such  that $A$ is contained in a finite union of at most $m$ sums of $V$. 
Since $V$ is contained in $a_d(U^d)$,  therefore 
for any $a\in A$ there exist $k\leq md$ and $y_i\in U$ for $1\leq i \leq k$ 
such that $a=\sum_{i=1}^k y_i$. By our assumption
that $y_i \in U \subset A(k)[N^\infty]$ holds, this implies $a\in A(k)[N^\infty]$ and therefore
$A(k)=A(k)[N^\infty]$. This is impossible, and therefore $\dim(Y)=0$. \qed  

\medskip
Obviously, in the last lemma the integer $N$ could have been replaced by any supernatural number, provided it is not divisible by all primes.

\medskip
In the proof of the last proposition we use the following variant of a theorem of Drinfeld [Dr2]

\begin{Theorem} \label{Drinfeld}
Let $X_0$ be a variety over $\kappa$ and $E$ be a finite extension of $\Q$. Let $\lambda,\lambda'$
be nonarchimedian places of $E$ prime to $p=char(\kappa)$ and $E_{\lambda}$ and $E_{\lambda'}$ be the corresponding completions. Let $K_0$ be an irreducible $\overline E_{\lambda}$-adic perverse sheaf on $X_0$ so that for every closed point $x$ 
in an open dense smooth subset $U_0$ of its support $det(1-Fr_x t,K_{\overline x})$ has coefficients in $E$ and its roots are $\lambda$-adic integers. Then there exists a $\overline E_{\lambda'}$-adic perverse sheaf $K'_0$ on $X$ compatible with $K_0$ (i.e. having the same characteristic polynomials of the Frobenius operators $Fr_{\overline x}$ for all closed points $x$ in $X$ in some common algebraic number field extension of $E$).
\end{Theorem}

\medskip
{\it Proof}. In the case where $X=U$ is smooth and $K_0=F_0[\dim(X)]$ is a smooth etale $\overline E_{\lambda}$-adic sheaf $F_0$ on $X_0$ the assertion of the theorem is proven in [Dr2, thm.1.1]. 
In general, to construct $K'_0$ we may assume that the support of $K_0$ is dense in $X_0$. Shrinking $U_0$ we can then assume that $K_0\vert_{U_0} = F_0[\dim(X)]$ holds for
a smooth etale $\overline E_\lambda$-sheaf $F_0$ on $U_0$. 
We apply Drinfeld's theorem theorem over $U_0$ to construct $F'_0$ over $U_0$
and define $K'_0$ to be the intermediate extension of $F'_0[\dim(X)]$ to $X_0$.
It remains to show the compatibility of $K_0$ and $K'_0$ for all closed points $x$ in the
complement $Y_0 = X_0 \setminus U_0$. For any such (fixed) $x$ this amounts to show the compatibility conditions $f_m^{K}(x) = f_m^{K'}(x)$ in $E$ for all $m=1,2,..$.  

\medskip
{\it Preparation of $X$}. Since the problem is local around $x$ in $X_0$ for the etale topology, we may assume that $X_0$ is affine.
We can replace then $X_0$ by $\mathbb A_0^n$ (affine space) by some closed embedding
$i_X: X_0 \to  \mathbb A_0^n$ so that the point $x$ maps to the point $0$. The perverse sheaves $K_0$ resp. $K'_0$ can be replaced by the perverse sheaves $i_{X,*}(K_0)$ resp. $i_{X,*}(K'_0)$, which by assumption are compatible on the complement of $i_X(Y_0)$. 
On can find an elliptic curve $\mathbb E_0$ over $\kappa$ with a double covering $\mathbb E_0 \to \mathbb P_0^1$
that is etale over the point $0$ in $0\in \mathbb A^1(k) \subset \mathbb P^1(k)$. 
By embedding $\mathbb A_0^n$ into $(\mathbb P^1_0)^n$ we can pullback $K_0$ and $K'_0$
to pure perverse sheaves on an open dense subset of $\mathbb E_0^n$.

\medskip
The pullbacks of $K_0,K'_0$ are pure (of the same weight) but in general need not be irreducible any more. However, without restriction of generality we may assume that
$K_0$ and $K'_0$ remain irreducible, since for the underlying pure etale sheaves (considered on a suitable regular dense subset of the supports) the compatibility conditions are inherited to their
irreducible direct summands. Indeed, this last remark follows from the dictionary [KW, thm. 12.1]
after applying Drinfeld's theorem [Dr2, thm.1.1] to each of the summands. So we will for simplicity assume that $K_0$ and $K'_0$ are irreducible, and extend $K_0$ and $K'_0$ to irreducible perverse sheaves on $\mathbb E_0^n$ by intermediate extension. 
So from now on we assume $X_0 = \mathbb E_0^n$ and proceed with the proof by using induction on $n$.

\medskip
{\it Inductionstep.} Consider any projection $f: X_0 = \mathbb E_0^n \to \mathbb E_0$. The fibers of $f$ are isomorphic to 
$\mathbb E_0^{n-1}$ and there exists a Zariki dense open subset $V_0\subset \mathbb E_0$ with the following property: Up to a complex shift by $\dim(\mathbb E)=1$, for any closed point $t\in V_0$ of $\mathbb E_0$ the restrictions of $K_0,K'_0$ to $f^{-1}(t)$ are perverse sheaves and $Y_0$ does not contain the whole fiber so that $K_0\vert_{f^{-1}(t)}$ is the intermediate 
extension of its restriction to the complement of $f^{-1}(t) \cap Y_0$ and similar for $K'_0$; for this see lemma \ref{gener}. Since $K_0,K'_0$ are pointwise pure and smooth on their support on the complement of $Y_0$, the same holds for their restrictions to  $f^{-1}(t) - (f^{-1}(t) \cap Y_0)$. As intermediate extensions
then $K_0\vert f^{-1}(t) \cap Y_0$ and  $K_0'\vert f^{-1}(t) \cap Y_0$ are pure on $f^{-1}(t)$. As explained above, we can therefore apply the induction assumption to them. 
So for closed points $t$ in $V_0$, $K_0$ and $K'_0$ are compatible
at all closed points of the fiber $f^{-1}(t)$. The complement $S_0$ of $V_0$ in $\mathbb E_0$ is finite:  $\dim(S_0)=0$.
We can use this argument for all the projections $f=pr_i$, $i=1,...,n$. This shows
that $K_0$ and $K'_0$ satisfy the local compatibility conditions at all closed points $x$ of $\mathbb E_0^n$
in the complement of $S=\prod_{i=1}^n S_0(i)$. Since $\dim(S_0(i))=0$, this complement $S$
is finite and by enlarging the base field $\kappa$, we may assume that all points in $S$ are $\kappa$-rational points. Taking quotients, the Grothendieck
Lefschetz fixed point formula for $K$
$$ \prod_{x\in \vert X_0 \vert} \det(1- t^{d(x)}F_x ; {\cal H}(K)_{\overline x})^{-1} \ = \
\prod_{i=-n}^n \det(1 - t F; H^i(X, K))^{(-1)^i}  \ $$
and similar for $K'$, together with the local compatibility conditions for $K_0$ and $K'_0$ outside $S$, 
gives
$$ \prod_{x\in S} \frac{\det(1- t^{d(x)}F_x ; {\cal H}^\bullet(K')_{\overline x})}{\det(1- t^{d(x)}F_x ; {\cal H}^\bullet(K)_{\overline x})}\ = \
\prod_{i=-n}^n (\frac{\det(1 - t F; H^i(X, K))}{\det(1 - t F; H^i(X, K'))})^{(-1)^i}  \ .$$
Since $K$ is irreducible perverse on $X$ (we can assume it to be pure of weight $w=0$)
and $X$ is proper over $k$, all eigenvalues of Frobenius $F$ on  $H^i(X, K)$ are of weight $i$.
Discarding the trivial case where the support of $K,K'$ is of dimension zero,
on the left side all contribution to poles and zeros are of weight $< 0$, since $K$ is perverse and pointwise of weight $\leq 0$.
Hence by weight reasons, $H^i(X,K)$ and $H^i(X,K')$ have the same $F$-eigenvalues for all $i=0,1,....,n$. Indeed,
no cancellation are possible on the right side of the formula by purity. The hard Lefschetz theorem for the pure 
perverse sheaves $K$ resp. $K'$ then also implies that $H^i(X,K)$ and $H^i(X,K')$ have the same $F$-eigenvalues for all $i=-1,....,-n$. Thus the right side of the last displayed formula turns out to be one, hence
$$ \prod_{x\in S} \det(1- t^{d(x)}F_x ; {\cal H}^\bullet(K)_{\overline x})^{(-1)^i}\  = \ \prod_{x\in S} \det(1- t^{d(x)}F_x ; {\cal H}^\bullet(K')_{\overline x})^{(-1)^i}  \ .$$
This implies
$$  \sum_{x\in S} f^K_m(x)  \ = \ \sum_{x\in S} f^{K'}_m(x)  $$
for all $m$. Now we are free to replace $K_0$, $K'_0$ by twists with characters
$\chi$ of $\pi_1(X)$. In particular we then obtain  
$$  \sum_{x\in S} f^K_m(x)\cdot \chi_m(x)  \ = \ \sum_{x\in S} f^{K'}_m(x) \cdot \chi_m(x)  $$
for  $\chi_m= \chi\circ S_m$ and all characters $\chi\in X(\kappa)^*$. So, by Fourier inversion and enlarging $\kappa$ we finally obtain
$    f^K_m(x)  \ = \  f^{K'}_m(x) $
for all $m$ and all $x\in S$ since we already know this for $x\notin S$. We therefore conclude
$$   \det(1- t^{d(x)}F_x ; {\cal H}(K)_{\overline x}) \ = \ \det(1- t^{d(x)}F_x ; {\cal H}(K' )_{\overline x}) $$
for all $x\in \vert X_0 \vert$. Now replace $t^{d(x)}$ by $t$. \qed

\goodbreak

\section{Invertible objects}\label{invertible}

For an abelian variety  $X_0$  and 
a perverse sheaf $P_0$ on $X_0$ defined over a finite field $\kappa$ let $X$ and $P$ be the extensions of $X_0$ and $P_0$ to 
an algebraic closure $k$ of $\kappa$. 

\begin{Proposition} \label{units} Suppose $P$ is a semisimple perverse sheaf on $X$
without negligible constituents. 
Then the following assertions are equivalent
\begin{enumerate}
\item $\chi(P)=1$ holds for the
Euler-Poincare characteristic of $P$ on $X$. 
\item $P*P^\vee \cong \delta_0 \oplus K$ holds for some negligible semisimple complex $K$ on $X$.
\item $P$ is an irreducible perverse sheaf and $P*P^\vee \cong \delta_0$ holds.
\item $P$ is an irreducible skyscraper sheaf on $X$.
\end{enumerate}
\end{Proposition}
 
\goodbreak 
 
\medskip 
We remark that over fields $k$ of characteristic zero this is shown in [KrW, prop. 10.1]
using analytic methods. Here we consider the case of finite fields. 

\medskip
Put $g=\dim(X)$. For the proof of the proposition  \ref{units}
we will first assume that $X$ is a {\it simple} abelian variety
and in that case we argue as follows:

\medskip
{\it Proof of implication 1) $\Longrightarrow$ 2)}.
Since $P$ is a perverse sheaf, being a multiplier $P$ corresponds to an irreducible representation of the quotient group
${\bf G}_{geom}(X)$ of ${\bf G}(X)$. Since ${\bf G}_{geom}(X)$ is a classical reductive group and the categorial dimension of $\overline P(X)$ is the Euler-Poincare characteristic, the corresponding representation has dimension $\chi(P)=1$ and hence defines a character $\mu$ 
of ${\bf G}(X)$. Then $P^\vee$ corresponds to the dual character $\mu^{-1}$. Therefore $P*P^\vee$ corresponds to $\mu\otimes \mu^{-1}=1$ and this implies that $P*P^\vee$ becomes isomorphic to
$\delta_0$ in $\overline{D}^{ss}(X)$. This and rigidity implies  $P*P^\vee \cong \delta_0 \oplus K$  
for some negligible semisimple complex $K$. This proves assertion  2) and notice, we did not use
that $X$ is simple for this part of the argument.

\medskip
{\it Proof of implication 2) $\Longrightarrow$ 3).} 
By the hard Lefschetz theorem for perverse sheaves [BBD], the Laurent polynomial 
$$  h_*(P) = \sum_{i}  \dim(H^i(X,P)) \cdot x^i  $$
is invariant under the substitution $x\mapsto x^{-1}$. Furthermore, as an easy consequence of Poincare duality, 
the Laurent polynomial $h_*(P^\vee)$ coincides with $h_*(P)$. Hence, for this Laurent polynomial
$h(x)$ the assumption $P*P^\vee \cong \delta_0 \oplus K$ implies 
$$  h(x)^2 = 1 + (x^{-1}+2+ x)^{g} \cdot f(x)  $$
for some other Laurent polynomial $f(x)$. 
That $h_*(K)$ has the form $(x^{-1}+2+ x)^{g} f(x)$ follows from our temporary assumption that
$X$ is a simple abelian variety. 
In the local ring 
$ R_{\Q} = \Q [x,x^{-1}] /(x+x^{-1}+2)^{g}  \cong  \Q[y]/y^{2g} $ 
this implies $h(x)^2 \equiv 1$. Hence $h(x) \equiv 1$
because $\chi(P) = h(-1) = 1$.
In other words, $h(x)= 1 + (x^{-1}+2+ x)^{g} g(x)$ for some
Laurent polynomial $g(x)$. Since by our assumptions $P$ is not translation 
invariant, the cohomology $H^i(X,P)$ vanishes for $\vert i\vert \geq g$. This implies 
$g(x)=0$ and hence $h_*(P)=1$, and therefore $h_*(K)=0$. Similarly $h_*(K_\chi)=0$ follows for all character twists
$K_\chi$ of $K$. This implies $K=0$, using corollary \ref{key} respectively  its amplification formulated in lemma \ref{strongerversion}.  Thus we have shown assertion 3). It is clear that $P*P^\vee \cong \delta_0$ implies that $P$ must be an irreducible perverse sheaf on $X$.

\medskip
{\it Proof of implication 3) $\Longrightarrow$ 4).}
Consider the addition law
$$a:X\times X\to X\ .$$ 
For $x\in X$ the fibers $a^{-1}(x)$ of $a$ can be identified with
the abelian variety $X$ by the map $a^{-1}(x)\ni (x-y,y)\mapsto y\in X$.  For fixed $x\in X$ we use the abbreviation 
$$  Q = (P\boxtimes P^\vee)\vert_{a^{-1}(x)}[-g] $$
for the restriction of the perverse sheaf $P\boxtimes P^\vee$
on $X\times X$ to the fiber $a^{-1}(x)$ in $X\times X$, up to a shift by $g= \dim(X)$.
By the proper base change theorem
the cohomology $H^\bullet(a^{-1}(x), Q[g])$ of $Q[g]$ can be identified with the stalk
$(P*P^\vee)_x$ of the convolution product $P*P^\vee$. By the assumption  $P*P^\vee\cong \delta_0$
 both $Q[g]$ and $Q$ are acyclic complexes for $x\neq 0$.

\medskip
By Lemma \ref{gener} below there exists a 
nonempty Zariski dense open subset $U$ of $X$ such 
that the complex $Q$
is in $Perv(a^{-1}(x))$ for $x\in U$. Since we may assume $0\notin U$, 
$Q$ is an acyclic perverse sheaf on $a^{-1}(x)$ for all $x\in U$.
Identifying $X$ and $a^{-1}(x)$, we can apply the main result of [W2]  
that all irreducible Jordan-H\"older
constituents of an acyclic perverse
sheaf $Q$ are translation invariant. 
Indeed, since $Q$ is acyclic and therefore $\chi(Q)=0$ holds,
all Jordan-H\"older constituents of $Q$ have vanishing Euler-Poincare characteristic
and therefore are translation invariant by the main theorem of [W2] ($X$ is simple). 
%Since $Ext^1_{Perv(X)}(\delta_\psi, \delta_{\psi'})=0$
%for $\psi'\neq \psi$, such perverse sheaves are a direct 
%sum of their different $\psi$-block components. Each $\psi$-block component is an iterated
%extension of sheaves $\delta_X^\psi$ 
In particular, for $Q\neq 0$ there
exists a nontrivial quotient morphism $Q \to \delta_X^\psi$ for some character $\psi$. 
This implies $H^g(a^{-1}(x), Q_{\psi^{-1}})\neq 0$. Notice $Q_{\psi^{-1}} \cong  (P_{\psi^{-1}}\boxtimes P^\vee)\vert_{a^{-1}(x)}[-g]$. Hence the 0-th stalk cohomology 
${\cal H}^0(Ra_*(P_{\psi^{-1}})\boxtimes P^\vee)_x)$ of  $ Ra_*(P_{\psi^{-1}})\boxtimes P^\vee)_x = (P_{\psi^{-1}}*P^\vee)_x$ does not vanish. Since $P$, and hence both its twist $P_{\psi^{-1}}$ and its Tannaka dual $P^\vee$, are invertible perverse sheaves, the semisimple complex $P_{\psi^{-1}}*P^\vee= Ra_*(P_{\psi^{-1}}\boxtimes P^\vee)$ is the direct sum of a simple invertible perverse sheaf $I$ and a negligible semisimple complex $T$ by the same argument that we used in the first step of the proof.

\medskip
Suppose $Q\neq 0$ holds for some fixed $x\in U$, then there are two possibilities: Either the semisimple sheaf complex $P_{\psi^{-1}}*P^\vee= Ra_*(P_{\psi^{-1}}\boxtimes P^\vee)$ has the skyscraper sheaf $\delta_x$ as one of its irreducible constituents (these are irreducible perverse sheaves up to complex shift). In this case $\delta_x$ must be the above mentioned unique invertible summand $I$.   
Or otherwise ${\cal H}^0((P_{\psi^{-1}}*P^\vee)_x)\neq 0$ implies ${\cal H}^0(T_x)\neq 0$, so by [KrW] for some character $\chi$ there exists a translation invariant irreducible direct summand $L\cong \delta_X^\chi[-g]$  of the semisimple complex $T\subseteq P_{\psi^{-1}}*P^\vee$ such that ${\cal H}^0(L)\neq 0$. But then ${}^p H^g Ra_*(P_{\psi^{-1}}\boxtimes P^\vee) \neq 0$ for the $g$-th  perverse cohomology. By the perverse cohomology bounds for smooth morphisms [BBD, 4.2.4] and [BBD, prop. 4.2.5] and the irreducibility of $P_{\psi^{-1}}\boxtimes P^\vee$ on $X\times X$, this  forces $P_{\psi^{-1}}\boxtimes P^\vee \cong a^*(M[g])$ for some irreducible perverse sheaf $M$ on $X$. Hence $P_{\psi^{-1}} * P^\vee \cong Ra_*(\Lambda_{X\times X})[g] \otimes^L M$, and this would yield a contradiction since the left side has Euler-Poincare characteristic $\chi(P_{\psi^{-1}})\chi(P^\vee)$ equal to  $\chi(P)^2 = 1$, where as the right hand side has  Euler-Poincare characteristic
equal to zero independent from $M$ since $Ra_* \Lambda_X[g] \otimes M = \bigoplus_{i=-g}^g {2g \choose g+i} \cdot M[i]$. So this excludes the second possibility. For $Q\neq 0$
it therefore  remains to discuss the case  $P_{\psi^{-1}} * P^\vee \cong \delta_x \oplus T$.
By irreducibility and rigidity, for fixed $x\in U$  this implies $P_{x} \cong P_{\psi^{-1}}$. However, such isomorphisms $P_{x} \cong P_{\psi^{-1}}$ can not exist for all $x$ in a Zariski dense open subset $U$ of $X$. Indeed since the support (resp. ramification locus of $P_{\psi^{-1}}$ within its support) are the same as the support $S$ (resp. the ramification locus of $P$), the support of $P$ and the ramification locus of $P$ would otherwise be invariant by translation under all $x\in U$. This would imply $S=X$ and $P$ would be unramified on all of $X$, and hence $P$ would be translation invariant so that $\chi(P)=0$ holds contradicting $\chi(P)=1$. This being said, it is clear that the discussion above implies that $Q$ vanishes for all $x$ in a dense Zariski open subset of $X$. Now replace $U$ by this Zariski dense open subset. For simplicity we again call this $U$. 
We then know that $P\boxtimes P^\vee$ vanishes on $a^{-1}(U)$. In other words,
all stalks $(P\boxtimes P^\vee)_{(x-y,y)} \cong P_{x-y} \boxtimes D(P)_{-y}$ vanish for $x\in U$ and $y \in X$.
Thus $-y \in S$ implies $x-y \notin S$; i.e. $x\notin S-S$ holds
for $x\in U$. As a consequence $\dim(S-S) <\dim(X)$ and
in particular $\dim(S) < \dim(X)$. 

\medskip
Recall, the points $(x-y,y)$ in the fiber $a^{-1}(x)$ are parametrized by the points $y\in X$.
The stalk of $Q$ at $(x-y,y)$ vanishes unless $x-y\in S$ and $-y \in S$.
In particular, the support of the complex $Q$ on $a^{-1}(x)$ has dimension not larger 
than $\dim(S)$. So the support of $Q$ is contained in a proper 
closed subset of the fiber $a^{-1}(x)$ for any $x\in X$ by dimension reasons.

\medskip
Let us analyze the last arguments more carefully. If $S\subset X$ is the irreducible support of $P$, then $Y=S\times -S$ is the support of the irreducible perverse sheaf $P\boxtimes P^\vee$ on $X\times X$. Of course, we can view $P\boxtimes P^\vee$ as a perverse sheaf on $Y$.
The image $a(Y)$ under $a: X\times X \to X$ is the irreducible set $Z:= S-S$. Consider the restriction $f: Y\to Z$
of the morphism $a$.
By lemma \ref{gener} there exists
a Zariski open dense subset $V$ of $Z$ such that the restriction of $P\boxtimes P^\vee$
to $f^{-1}(x)$ is perverse up to a shift by $\dim(Z)$. For $x\neq 0$ in $V$ it is in addition
an acyclic complex on $f^{-1}(x) \subset a^{-1}(x) \cong X$. So, by the previously used argument,
all irreducible Jordan-H\"older constituents of this perverse sheaf are translation invariant by $X$ on 
$a^{-1}(x)$. So the support includes the whole fiber. On the other hand $(S\times -S)\cap f^{-1}(x) =\{y\in X\ \vert\ x-y\in S\ , \ -y\in S\}$, which  is a subset of $a^{-1}(x)$ of dimension $\leq \dim(S) < g$. A contradiction which shows that no point $x\neq 0$ can be contained in $V$.
Hence $dim(Z)=0$ and then $Z=S-S=0$. In particular $\dim(S)=0$ and $P$ is a skyscraper
sheaf.

\medskip
{\it The non-simple case}. Our previous arguments prove proposition \ref{units} for simple abelian varieties.
To discuss general abelian varieties we use induction on $\dim(X)$. So, let us assume that proposition \ref{units} holds for all abelian varieties of dimension less than
$\dim(X)$. Since our previous proof of the implication 1) $\Longrightarrow $ 2) did not use
that $X$ is simple, we can assume  $P*P^\vee \cong \delta_0 \oplus K$ for some negligible semisimple complex $K$. Since $P$ is clean, 
in particular $P$ must be simple.
If $X$ is not simple, replacing $\kappa$ by some finite field extension we may assume that there exist a surjective homomorphism $f_0: X_0 \to Y_0$ to some abelian
variety $Y_0$ over $\kappa$. We can assume that the kernel of $f$ is a simple 
abelian subvariety $A$ of $X$ and using an isogeny, we may furthermore assume 
$X_0 \cong A_0 \times Y_0$.  This allows us to view the characters of $\pi_1(A)$
as characters of $\pi_1(X)$.
%To show that $P$ is a skyscraper sheaf, we are free to replace $P$ by some
%character twist $P_\chi$. 

\medskip
Notice, $K = \bigoplus_\nu K_\nu$ with $K_\nu = P_\nu[-m_\nu]$ for irreducible perverse sheaves $P_\nu$ and $m_\nu < depth(P_\nu)$, where $P_\nu$ is translation invariant with respect to an abelian subvariety $B_\nu$ of dimension $\dim(B_\nu)= depth(P_\nu)>0$; since $P$ is clean, this follows
from lemma \ref{elong} in the appendix. For generic characters $\chi$ of $\pi_1(A)$ the direct image 
$Rf_*(P_\chi)$ is perverse by the generic vanishing theorem [W], and the semisimple complex  $Rf_*(K_\chi)$ is  negligible by lemma \ref{nice}
(since $A$ is simple, \lq{generic}\rq\ in both instances can be replaced by the condition almost all). More precisely we can achieve $Rf_*(K_{\nu,\chi})=0$ whenever $A$ is contained in $B_\nu$, and otherwise $Rf_*(K_{\nu,\chi})= Rf_*(P_{\nu,\chi})[-m_\nu]$ for perverse negligible
sheaves $\tilde P_\nu := Rf_*(P_{\nu,\chi})$ of $depth(\tilde P_{\nu}))= depth(P_\nu)$. Hence
$Rf_*(K_\chi) = \bigoplus_\nu \tilde P_\nu[-m_\nu]$ and $m_\nu < depth(\tilde P_\nu)$.
Since $P_\chi*P_\chi^\vee \cong \delta_0 \oplus K_\chi$ implies
$$Rf_*(P_\chi)*Rf_*(P_\chi)^\vee \cong \delta_0 \oplus Rf_*(K_\chi)\ $$ 
and since $End_{D(X)}(Rf_*(P_\chi)) = Hom_{D(X)}(Rf_*(P_\chi)*Rf_*(P_\chi)^\vee,\delta_0)$  by  rigidity, we get $End_{D(X)}(Rf_*(P_\chi)) \cong \Lambda$. Indeed  
$Hom_{D(X)}(\tilde P_\nu[-m_\nu],\delta_0) \cong {\cal H}^0(\tilde P_\nu[-m_\nu])_0^\vee  $ and this stalk cohomology vanishes because $m_\nu < depth(\tilde P_\nu)$.
Hence the perverse sheaf $Rf_*(P_\chi)$ is not only semisimple, it is simple. Then clearly  
it can not be 
negligible. So we can apply the induction assumption for the perverse sheaf $Rf_*(P_\chi)$ on $Y$ to  show $Rf_*(P_\chi) \cong \delta_y$ for some 
closed point $y$ in $Y$. Without restriction of generality we may assume $y=0$ and $\chi=1$. 
As a consequence the 
sheaf complex $P\vert_{f^{-1}(y)}$ is an acyclic sheaf complex on
$f^{-1}(y) \cong A$  for all $y\in W=Y - \{ 0\}$ (proper base change theorem). Hence for almost all characters $\chi$ of $\pi_1(A)$
the direct image complex $Rf_*(P_\chi)$ vanishes on the Zariski open dense subset
$W$ of $Y$.  For (the closure of) the support $S$ of $P$ in $X$ we claim $f(S) = \{ 0\}$.
For this it suffices to show $S+A=S$ if $f(S) \neq \{ 0\}$, since
then [W, prop. 2] would imply that $P$ is translation invariant by $A$ and hence negligible, contradicting our assumptions on $P$. Suppose $f(S) \neq \{ 0\}$ and choose some point $y\in f(S)$ in general position. Then $y\neq 0$ and the restriction of $P$ to $f^{-1}(y)$ does not vanish. However, since the twisted complex
$(P\vert_{f^{-1}})_\chi$ is acyclic for almost all characters $\chi$ of $\pi_1(A)$, lemma \ref{strongerversion} implies that all irreducible perverse constituents of the perverse cohomology sheaves of this complex are translation invariant by $A$. Hence for all
$y\neq 0$ in $f(S)$, $f^{-1}(y)$ is contained in $S$. This immediately implies $S+A=S$ if $f(S) \neq \{ 0\}$. Again, [W, prop.2] 
gives a contradiction that implies
$f(S)=\{ 0\}$. So $P$ can be viewed as a perverse sheaf on the fiber $f^{-1}(0)\cong A$. Since $A$ is simple by construction, we can apply the first part of our proof to show that $P$ is a skyscraper sheaf. This completes the proof
of proposition \ref{units}. \qed
 
 \bigskip\noindent
 
\goodbreak

\begin{Lemma}\label{gener}
Let  $Z$ be a variety over $k$ of dimension $n$ and let $f: Y \to Z$ be a morphism between algebraic varieties over $k$.  Then
\begin{enumerate}
\item
For a perverse sheaf $K$ on $Y$ there exists a Zariski dense open subset $U$ of $Z$ such that 
for closed points $x\in U$ the shifted restrictions $Q=(K[-n])\vert_{f^{-1}(x)}$ of $K$ to the fibers $f^{-1}(x)$   are perverse sheaves. 
\item
If $K$ is irreducible and defined as the intermediate extension of a local system on some dense Zariski open subset $Y\setminus S$ of $Y$, then we can choose $U\subset Z$ as before so that in addition for all $x\in U$ the perverse sheaf $Q$ is the intermediate extension of its restriction from $f^{-1}(x)$
to the complement of $S_x=S \cap f^{-1}(x)$.
\end{enumerate}
\end{Lemma}

\medskip
{\it Proof}. 
We can shrink $Z$ to some smooth open dense subset $U$ and replace $Y$ by $f^{-1}(U)$ so that the supports $S_\nu$ of the cohomology sheaves
${\cal H}^\nu(K)$ satisfy that $S_\nu \to U$
is surjective unless $S_\nu$ is empty. In addition we can achieve that
the dimensions of $S_\nu \cap f^{-1}(x)$ do not depend on $x\in U$.
Fix some $x\in U$.
Notice that $K\vert_{f^{-1}(x)}$ is in ${}^p D^{[-n,0]}(f^{-1}(x))$ by the Artin-Grothendieck affine vanishing theorem and the existence of a local regular sequence $x_1,..,x_n$ at $x\in U$ 
for $\dim(U)=n$. Indeed, the (locally defined) inclusion $i:Y_1 \to Y$ of the zero locus $x_1=0$
has the property $i^*({}^p D^{[a,b]}(Y))\subseteq {}^p D^{[a-1,b]}(Y_1)$; see [KW,p.154].
Inductively this implies ${}^p D^{[a,b]}(Y)\vert_{f^{-1}(x)} \subseteq {}^p D^{[a-n,b]}(f^{-1}(x))$.
So it suffices to show $K\vert_{f^{-1}(x)} \in {}^p D^{\leq -n}(f^{-1}(x))$. If this were not true,
then there exists $\nu$ with $\dim supp({\cal H}^\nu(Q)) + \nu >0$. Since
$\dim supp({\cal H}^\nu(Q)) = \dim supp({\cal H}^\nu(K[-n])) - \dim(U) = 
\dim supp({\cal H}^{\nu-n}(K)) - n$ by the equidimensionality of the supports,
then $\dim supp({\cal H}^{\nu-n}(K)) + (\nu-n) >0$.
Hence $K \notin {}^p D^{\leq 0}(Y)$. Since $K$ is a perverse sheaf on $Y$, this
gives a contradiction and proves our first claim. 

\medskip
From this we obtain for fixed $M \in {}^p D^{<m}(Y)$ that $Q_M = (M[-n])\vert_{f^{-1}(x)}$ satisfies $Q_M \in {}^p D^{< m}(f^{-1}(x))$ for all $x\in U$, after suitably shrinking the base $U$.

\medskip 
For $f:Y\to Z$ let $D_{Y/Z}$ denote the relative duality functor [KL, 1.1].
Then, in addition to the previous choices, we can also assume that $K$ and 
$D_{Y/Z}(K)$ are reflexive relative over $U$ in the sense of [KL, prop.1.1.7] and [KL, def.1.1.8].
Then by the first part of lemma \ref{gener} the complex  $K[-n]$ is relatively perverse over $U$
[KL, 1.2.2(ii)]. For simplicity of notation, from now on suppose $U=Z$. 
We claim $L[-n]:=D_{Y/Z}(K[-n]) \in Perv(Y)[-n]$ for an irreducible perverse sheaf $L$ on $Y$.
Since $K[-n]$ is relatively perverse over $Z$, so is $D_{Y/Z}(K[-n])$. This follows from the definition of relative perversity and the fact that $D_{Y/Z}(D_{Y/Z}(K[-n])) = K[-n]$ holds by  [KL, prop. 1.1.6] since $K$ (now) is reflexive over $Z$. Because in our situation relative duality  commutes with base change over $Z$, the supports
of $K[-n]$ and $D_{Y/Z}(K[-n])$ in the fibers are the same. By the equidimensionality of the fibre supports over $U$, we therefore obtain $\dim(supp({\cal H}^i(D_{Y/Z}(K[-n]))) = \dim(Z) + \dim_Z(supp({\cal H}^i(D_{Y/Z}(K[-n]))) \leq n - i$ and hence $L[-n]=D_{Y/Z}(K[-n]) \in {}^p D^{\leq n}(Y)$. For $L_0= {}^p H^0(L)$ we get a distiguished triangle $(M,L,L_0)$ with $M\in {}^p D^{< 0}(Y)$.
By the remark preceding this paragraph and the first assertion of lemma \ref{gener} therefore
$$    (L_0[-n])\vert_{f^{-1}(x)} \in Perv(f^{-1}(x)) \quad , \quad     (M[-n])\vert_{f^{-1}(x)} \in {}^p D^{< 0}(f^{-1}(x))$$
for all $x$ in  the base $Z$, after suitably shrinking $Z$. 
On the other hand $    (L[-n])\vert_{f^{-1}(x)} = D((K[-n])\vert_{f^{-1}(x)}) = DQ$
is the Verdier dual of the perverse sheaf $Q$ since relative duality commutes with base change in our situation.
Hence $    (L[-n])\vert_{f^{-1}(x)} = DQ$  is in $Perv(f^{-1}(x))$.
The distinguished triangle $(M[-n], L[-n], L_0[-n])\vert_{f^{-1}(x)}$
then implies $ (M[-n])\vert_{f^{-1}(x)} = 0$. This holds for all $x$ and therefore the restriction of $M$ to all fibers vanishes. Therefore $M=0$, and
this implies that $L$ is a perverse sheaf on $Y$.

\medskip
Since $D_{Y/Z}(L[-n])=D_{Y/Z}^2(K[-n]) = K[-n]$, the number of perverse Jordan-H\"older constituents of $K$ and $L$ must
coincide, after suitably shrinking the base. Hence if $K$ is an irreducible perverse sheaf, also $L$ is an irreducible perverse sheaf. The support and smooth locus of $K$ and $L$
coincide (for this it is enough to see that this is true in the fibers, where this follows from the definition of $L$).  

\medskip
By assumption, $K$ is irreducible with support not contained in $S$.
For the inclusion $i_S: S \to Y$ therefore ${}^p H^0(i_S^*(K))=0$ holds,
or equivalently $\dim(S_\nu\cap S)= \dim( supp ({\cal H}^{-\nu}(i_S^*(K))) \leq \nu -1 $ for all $\nu$ [KW, prop III.9.3].
We now shrink $Z$ so that both $S$ and $S\cap S_\nu$ become equidimensional over $Z$.
The equidimensionality implies  $ \dim(f^{-1}(y)\cap  supp\ ({\cal H}^{-\nu}(i_S^*(K))) \leq \nu -1 - n $. So for all $\nu$
$$\dim(supp\ {\cal H}^{-\nu}(i_S^*(K[-n]\vert_{f^{-1}(y)}))) \leq \nu -1 \ .$$ 
$(i_S^*(K[-n])\vert_{f^{-1}(x)})$ is the pullback of $Q=K[-n]\vert_{f^{-1}(x)}$ to
$S_y:=S \cap f^{-1}(x)$. 
Again by [KW, prop.III.9.3] for $i_{S_x}: S_x \hookrightarrow f^{-1}(x)$ therefore $${}^p H^0(i_{S_x,*}(Q))=0 \ .$$
In other words: $Q$ does not have nontrivial perverse quotients with support
in $S_x \subset f^{-1}(x)$.
The Verdier dual $DQ$ of $Q$ is obtained by the restriction of $L[-n]=D_{Y/Z}(K[-n])$ to the fibre $f^{-1}(x)$ by [KL, prop. 1.1.7 and def. 1.1.8]. We have seen that the supports of $K$ and $L$ coincide. Hence $L$ is an irreducible  perverse whose support is not contained in $S$. Therefore the same argument applied to $L[-n]$ instead of $K[-n]$ now
also shows that $L[-n]\vert_{f^{-1}(x)}$  does not have nontrivial perverse quotients with support
in $S_x \subset f^{-1}(x)$. In other words, neither $Q$ nor $DQ$ has perverse quotients with support
in $S_x$. Hence by [KW, lemma 5.1] the perverse sheaf $Q$ is the intermediate extension
of its restriction to $f^{-1}(x) \setminus S_x$. This proves the second assertion.
\qed

\bigskip\noindent

\goodbreak

\section{Appendix}%: Translation invariant perverse sheaves}
\label{appendix}

Let $X$ be an abelian variety over $k$.
The simple constituents of a semisimple complex
$K =\oplus_\nu {}^p H^\nu(K)[-\nu]$ with perverse cohomology sheaves ${}^pH^\nu(K)$  are the simple perverse constituents of the perverse sheaves ${}^p H^\nu(K)$. Let $N_{Euler} \subset D^{ss}(X)$ denote the additive subcategory generated by negligible objects where
for an abelian variety over $k$
{\it negligible} objects are semisimple complexes $K$ on $X$ whose simple perverse constituents $K_i$  have Euler-Poincare characteristic $$\chi(K_i)=\sum_\nu (-1)^\nu \dim(H^\nu(X,K_i)))=0\ .$$ 
%For characters $\chi$  of the fundamental group $\pi_1(X,0)\to \Lambda^*$ of an abelian variety $X$, the character twist $K_\chi$ of a complex  $K\in N_{Euler}$ is in $N_{Euler}$. 
A complex in $D^{ss}(X)$ will be called {\it clean} if it does not contain simple constituents from $N_{Euler}$.
A complex $K$ in $D^{ss}(X)$ decomposes into a direct sum $K=M\oplus R$ of a clean complex $M$ and a negligible complex $R$, and $M$ and $R$ are unique up to isomorphism. 
 
 \medskip
A simple perverse sheaf $K$ in $N_{Euler}$ is of the form
$$      K_\chi = q^*(L)[dim(A)]   $$
for some character $\chi$ of $\pi_1(X,0)$, some quotient $q: X\to X/A$ with respect to an abelian
subvariety $A\subseteq X$ of dimension $>0$, and some simple perverse sheaf $L$ on $X/A$ (see [KrW], [W2]).  We can assume that the perverse sheaf $L$ on $X/A$ is clean. If $A$ is chosen in this way, the dimension $l=\dim(A)$ will be called $depth(K)$. In general, for a complex $K\in N_{Euler}$, we define $depth(K) = \min_i(depth(K_i))$ for $K_i$ running over the simple perverse constituents of $K$. 

\medskip
For a semisimple complex $K$ the stabilizer $Stab(K)$, considered as an abstract group, is the subgroup of $X(k)$ defined by all $x\in X(k)$ for which
$T_x^*(K)\cong K$ holds. 

\begin{Lemma} \label{closed} For a semisimple perverse sheaf $K$ the stabilizer 
$Stab(K)$ is Zariski closed in $X$.
\end{Lemma}

{\it Proof}. Assume $K$ is simple. The  Zariski closure $Z$ of $Stab(K)$ is a closed subgroup of $X$. So, $A=Z^0$ is an abelian subvariety of $X$. Since $K$ is simple, $x \in Stab(X) \Leftrightarrow T_x^*(K) \cong K \Leftrightarrow {\cal H}^0(K*K^\vee)_x \neq 0$ holds by [BN]), hence $Stab(X)$ is a constructible subgroup of $X$ and  $Stab(K)\cap A$ is a constructible dense subgroup of $A$.  Therefore it contains a Zariski open subset $U$ of $A$.  If we replace $U$ by $V=U\cap -U$, then $V_1=V$, $V_2=V+V$, $V_3=V+V+V,..$ etc. are dense open subsets of $X$ contained in $Stab(K)$. Therefore $Stab(K)$ contains $W=\bigcup_{i=1}^\infty V_i$, a  Zariski open dense subgroup of $A$. Therefore finitely many translates of $W$ cover $A$, and
$A$ is the union of finitely many cosets of $W$. This implies that $W$ is Zariski closed, hence $W=A$ and $Stab(K)=Z$.  By the next lemma \ref{dep} the proof in the general case is easily reduced to the case where $K$ is simple. \qed

\begin{Lemma} \label{dep} If K is a semisimple complex
and $A \subseteq Stab(K)$ is an abelian variety, then $A$ is contained in $Stab(K_i)$ for all simple perverse constituents $K_i$ of $K$. 
\end{Lemma}

{\it Proof}.  For representatives $K_1,..,K_n$ of the isomorphism classes of the simple constituents of ${}^p H^\nu(K)$ consider their $n$ isotopic blocks in ${}^p H^\nu(K)$. Looking at Jordan-H\"older series, for $x\in A$ we get $T_x^*(K_i) \cong K_{i(x)}$ so that $i\mapsto i(x)$ defines a permutation in $S_n$. Varying $x$, we get a homomorphism
of $A$ to the permutation group $S_n$,  hence a subgroup $A'$ of $A$ of finite index stabilizes each $K_i$ up to isomorphism: $A' \subseteq A\cap Stab(K_i) \subseteq A$ for $i=1,..,n$. 
Then $A\cap Stab(K_i)$ has finite index in $A$. But $A\cap Stab(K_i)$ is Zariski closed (see lemma \ref{closed}), this implies $A\cap Stab(K_i) = A$ for all $i$. Hence $A'=A$ and $A\subseteq Stab(K_i)$.  \qed

\medskip
The proof of lemma \ref{dep} also shows that $Stab(K)$ contains $\bigcap_{i=1}^n Stab(K_i)$
as a normal subgroup of finite index, with quotient group contained in $S_n$. Hence $Stab(K)$
is a finite union of translates of  $\bigcap_{i=1}^n Stab(K_i)$. In particular, $Stab(K)$ is Zariski closed.
This shows the claim at the end of lemma \ref{closed}.
 
By lemma \ref{closed}, the connected component $A=Stab(K)^0$ of $Stab(K)$ is an abelian subvariety of $X$, that in the following is called the connected stabilizer of $K$.
For any simple perverse sheaf $K$ the dimension of its connected stabilizer is $depth(K)$.

\begin{Lemma} \label{7} For an abelian variety $X$ over $k$ let $f:Y\to X$ and 
$g:X\to Z$ denote isogenies between abelian varieties.
For negligible semisimple objects $K\neq 0$ in $D^{ss}(X)$ the following holds: 
\begin{enumerate}
\item Any character twist $K_\chi$ of $K$ is negligible. 
\item $f^*(K)$ is a negligible semisimple complex with $depth(f^*(K)) = depth(K)$. 
\item $g_*(K)$ is a negligible semisimple complex with $depth(g_*(K)) = depth(K)$. 
\item $g: Stab(K)^0 \to Stab(g_*(K))^0 $ and $f:Stab(f^*(K))^0 \to Stab(K)^0$ are isogenies for
simple perverse sheaves $K$.
\item For most characters $\chi$ of $\pi_1(X,0)$ the twist $K_\chi$ is acyclic, i.e.
$H^\bullet(X,K_\chi)=0$.
\item If $K$ is perverse, there exists $\chi$ so that  $H^\nu(X,K_\chi)\neq 0$ holds for some $\nu\neq 0$.
\end{enumerate}
\end{Lemma}

\medskip{\it Proof}. For the proof we may assume that $K$ is perverse. For the first assertion see [KrW, cor.10]. For the second and third assertion use 
that $\chi(K)=0$ for perverse sheaves $K$ implies $\chi(K_i)=0$ for all its simple perverse constituents $K_i$ (for references see [KrW] and [W2]): For an isogeny 
$f:Y\to X$ with kernel $F$ the complex $f^*(K)$ is perverse and $Rf_*(f^*K) = \bigoplus_{\chi \in F^*} K_\chi$.  Hence $\chi(Rf_*(f^*(K)))=\chi(f^*(K))=0$ by assertion 1, so all perverse constituents of $f^*(K)$ are negligible by the above remark. 
Similarly $g_*(K)$ is perverse and $\chi(g_*(K))=\chi(K)=0$. This proves the first part of the
assertions 2 and 3.

\medskip
For the remaining part of the assertions 2 and 3 on the depth we may assume that $K$ is a simple perverse sheaf.  
For the abelian variety $A=Stab(K)^0$ and $y\in f^{-1}(A)$, then
$\varphi: K \cong T^*_{f(y)}(K)$ implies $f^*(\varphi): f^*(K) \cong f^*(T^*_{f(y)}(K)) =
T_y^*(f^*(K))$. Thus $y\in Stab(f^*(K))$,
and the abelian variety $f^{-1}(A)^0$ stabilizes $f^*(K)$ and then also each simple component of $f^*(K)$ (lemma \ref{dep}). This implies $depth(f^*(K)) \geq depth(K)$. 
Conversely, if an abelian variety $B\subseteq Y$ stabilizes a constituent of $f^*(K)$, then for each $y\in B$ there exists a nontrivial map in $$ Hom(f^*(K),T^*_y(f^*(K)))
= Hom(f^*(K),f^*(T^*_{f(y)}(K)))= Hom(K,\oplus_{\chi \in F^*} T^*_{f(y)}(K_\chi))\ .$$ Hence
$T^*_{-f(y)}(K) \cong K_\chi$ for some $\chi\in F^*$, since $K$ is simple. Thus, for a subgroup $U\subseteq B$ of finite index and $y\in U$ we get $f(y)\in Stab(K)$. But $f(U) \subseteq Stab(K)\cap f(B)\subseteq f(B)$ implies that $Stab(K)\cap f(B)$ is of finite index in $f(B)$. Hence $Stab(K)$ contains $f(B)$ and $depth(f^*(K)) \leq depth(K)$. This shows  $depth(f^*(K))=depth(K)$ and $f: Stab(f^*(K))^0 \to Stab(K)^0$ is an isogeny.

\medskip
For the assertion on $depth(g_*(K))$, again let $K$ be a simple perverse sheaf. 
Put $A=Stab(K)^0 \subset X$. Then
$g(A)$ stabilizes $g_*(K)$ and therefore all its perverse constituents, hence 
$depth(K) \leq depth(g_*(K))$. Conversely, if $L$ is a simple constituent of
$g_*(K)$ and $\varphi: T_z^*(L) \cong L$ is an isomorphism for $z\in Stab(L)$, 
choose some $x\in X$ with $g(x)=z$. By the adjunction formula 
 $Hom(L,g_*(K)) = Hom(g^*(L),K)$  for perverse sheaves there exists a nontrivial morphism $g^*(L) \to K$. Hence
$K$ is a simple constituent of $g^*(L)$. Since $g^*(\varphi): g^*(T_z^*(L)) \cong g^*(L)$ and $g^*(T_z^*(L))= T_x^*(g^*(L))$ imply $g^*(L) \cong T_x^*(g^*(L))$, now 
$T_x^*(K) \cong K$ follows by lemma \ref{dep} for all $x\in g^{-1}(Stab(L))^0$. Hence $g^{-1}(Stab(L))^0 \subset Stab(K)^0$ and $depth(g_*(K)) \leq depth(K)$.

\medskip
For assertion 5 again we  may assume that $K$ is perverse and simple.
Then for the connected stabilizer $A\subseteq X$ of $K$ there exists an isogeny $f:A\times B\to X$
such that $H^\bullet(X,K) \subseteq H^\bullet(A\times B, L_\psi \boxtimes K_B) = H^\bullet(A,L_\psi) \otimes_\Lambda H^\bullet(B,K_B)$. For a twist of  $K$ with a character $\chi$, whose pullback
$(\chi_1,\chi_2)$ under $\pi_1(f): \pi_1(A,0) \times \pi_1(B,0) \to \pi_1(X,0)$ satisfies $\chi_1 \neq \psi^{-1}$,
the cohomology of $K_\chi$ vanishes. This implies assertion 5. For $\chi_1 =\psi^{-1}$ the same argument shows 
that the cohomology $ H^\bullet(A,\Lambda[\dim(A)]) \otimes_\Lambda H^\bullet(B,K_{B,\chi_2})$ does not vanish in degree $\nu=\dim(A)$ for most $\chi_2$ by the main result of [KrW].
Since any pair of characters $(\chi_1,\chi_2)$ can be written as the restriction of some character $\chi$ of $\pi_1(X,0)$, this implies
$H^{\dim(A)}(X,K_{\chi\chi_0})\neq 0$ for some of the finitely many characters $\chi_0$ of $\pi_1(X,0)$ which are trivial on the image of $\pi_1(f):  \pi_1(A,0) \times \pi_1(B,0) \to \pi_1(X,0)$. 
This proves assertion 6. \qed
  
\begin{Lemma} \label{nice}
Let $T$ be a semisimple complex on $X$ and let $A$ be an abelian subvariety of $X$.
For the quotient map  $f:X \to B=X/A$ suppose $depth(T)=\dim(B)>0$. Then for most $\chi$ the  complex $Rf_*(T_\chi)$ is either zero or negligible of depth $\dim(B)$.  If in addition $Stab(T)^0$ maps surjectively onto $B$, there exists a character $\chi$ so that  $Rf_*(T_\chi)$  is nonzero. If furthermore $X=A\times B$, then there exist a character $\chi_0$ so that for most
characters $\chi$ of $\pi_1(A,0)$, extended to characters of $\pi_1(X,0)$ that are trivial on $\pi_1(B,0)$, the complex $Rf_*(T_{\chi_0\chi})$ is negligible of depth $\dim(B)$, but not acyclic. 
\end{Lemma}

\medskip

{\it Proof}. For the proof we may replace $T$ by one of its semisimple perverse sheaves $K={}^p H^\nu(T)$ on $X$, without restriction of generality.  Put $D:=Stab(K)^0 \cap A$. By our assumptions then  $dim(D)=0$ holds iff $Stab(K)^0$ surjects onto $B$. 

\medskip
First suppose $dim(D)>0$. Then there exists an isogeny $g: D\times Y \to X$
so that $g^*(K)= L_\psi \boxtimes K_Y$ holds for some perverse sheaf $K_Y$ on $Y$ (see [W]). 
Since $K$ is a summand of $g_*g^*(K)$, as in the proof of lemma \ref{7}.5 we can replace $X$ by $D\times Y$ and $f$ by $f\circ g$ to show $Rf_*(K_\chi)=0$ for most $\chi$.
Indeed, then  $Rf_*(K_\chi)=0$ holds whenever the restriction of $\chi$ to $\pi_1(D,0)$ is different from $\psi^{-1}$, since $f\circ g$ factorizes over the projection $D\times Y \to Y$. This proves our first claim.

\medskip
If $dim(D)=0$, $B$ stabilizes $Rf_*(K_\chi)$. Therefore $Rf_*(K_\chi)$ is zero or negligible for all
$\chi: \pi_1(X,0)\to \Lambda^*$. 
By our assumptions $g(x,a)=a+x$  defines an isogeny $g: A \times Stab(K)^0 \to X$ with $f\circ g(x,a)= f(x)$. Choose an isogeny
$\tilde g:X \to A\times Stab(K)^0 $ so that $\tilde g \circ g$ is multiplication by an integer $n$
$$  \xymatrix{ A \times Stab(K)^0  \ar[d]^{pr_2} \ar[r]^-g  &  X \ar[d]^-f \ar[r]^-{\tilde g}  &  A \times Stab(K)^0  \ar[d]^{pr_2}\cr
Stab(K)^0  \ar[r]^-g  &  X/A  \ar[r]^-{\tilde g}  & Stab(K)^0  \cr} $$
By the arguments of [W] used for the proof of theorem 2 of loc. cit. it follows that $g^*(K) \cong K_A\boxtimes L_\psi $ holds for some invariant simple perverse sheaf $L_\psi$ on the factor $Stab(K)^0$ and some perverse sheaf $K_A$ on the factor $A$. 
By our assumptions on the depth, $K_A$ is not negligible for at least for one perverse cohomology degree $\nu$ (fixed above). 
Notice $g_*g^*(K) = g_*(  K_A\boxtimes  L_\psi )$ 
is a direct sum of twists $K_\chi$ of $K$ including $K$ itself.
To show $Rf_*(K_\chi)\neq 0$  for one $\chi$,  
we may therefore replace $f$ by $pr_2$, $X$ by $A \times Stab(K)^0 $ and $K$ by any simple clean constituent $L$ of $\tilde g_*(K)$. This immediately follows by the commutative diagram above. Now $\tilde g_*(K)$
is a retract of $\tilde g_* g_* g^*(K) = n_* K_A   \boxtimes n_* L_\psi  $, hence
isomorphic to a direct sum of sheaves $n_*K_A \boxtimes L_{\psi'} $. Notice $n_*K_A$ is not negligible and therefore can be replaced by its nonvanishing clean part.
Thus, without restriction of generality we are in the situation where $X=A \times B$ and $K= K_A\boxtimes L_{\psi'}  $
holds for a clean perverse sheaf $K_A\neq 0$. But then $ Rpr_{2,*}(K_{\chi_A\boxtimes \chi_B})= H^\bullet(A,K_{A,\chi_A}) \otimes L_{\psi'\chi_B}$. So, for $\chi_A\chi_0=(\chi_A,1)$ 
our claim follows since $\chi(H^\bullet(A,K_{A,\chi_A})) = \chi(K_{A,\chi_A}) = \chi(K_A) \neq 0$
holds for all characters $\chi_A: \pi_1(A,0)\to \Lambda^*$ by [KrW, cor.10].  \qed

\begin{Lemma} \label{elong}
For a clean semisimple perverse sheaf  $P$ on the abelian variety $X$
the complex $T=P*P^\vee$ is semisimple and the simple constituents 
$T_\nu$ of $T$ are of the form $T_\nu\cong P_\nu[-m_\nu]$ for simple perverse sheaves
$P_\nu$ on $X$ and integers $m_\nu$. For $m_\nu\neq 0$ we have $-depth(P_\nu) < m_\nu < depth(P_\nu)$. 
\end{Lemma}

\medskip
{\it Proof}. $T=\bigoplus_\nu P_\nu[-m_\nu]$ is semisimple (decomposition theorem). 
To show $m_\nu \leq depth(P_\nu)$ for fixed $\nu$, consider the abelian subvariety
$A = Stab(P_\nu)^0$ of $X$. By an isogeny $f: X\to A \times B$ we may replace $P$
by $f_*(P)$, $T$ by $f_*(T)$ and $P_\nu$ by $f_*(P_\nu)$ so that by lemma \ref{7}, without restriction of generality, 
we may assume $X = A\times B$ and 
$P_\nu = \delta_A^\psi \boxtimes Q$ for an irreducible clean perverse sheaf $Q$ on $B$. Twisting $P$ allows to assume $\psi=1$. Let $p: X \to A$ and $q: X\to B$ be the projections. By the generic vanishing theorems [KrW], for a generic character $\chi$ of $\pi_1(B)$ 
the direct image complexes $K=Rp_*(P_\chi)$ and $Rp_*(P_{\nu,\chi})= \delta_A \otimes_\Lambda H^0(B,Q_\chi) \neq 0$ are perverse sheaves on $A$. Hence  $K*K^\vee = \bigoplus_\nu  Rp_*(P_{\nu,\chi})[-m_\nu]$ implies $m_\nu \leq \dim(A)$ by the cohomological bounds [BBD, 4.2.4]. %Therefore
%${\cal H}^i(\delta_A[-m_\nu]) \otimes_\Lambda H^0(B,Q_\chi) =0$ for all $i >0$, and hence $m_\nu \leq \dim(A)= depth(P_\nu)$. 
Since $H^{2\dim(A)}(A,\delta_A[-\dim(A)])$ is nonzero,  
$H^{2dim(A)}(X,P_\chi*P_\chi^\vee)\neq 0$ if $m_\nu = \dim(A)$.
Then duality, the hard Lefschetz theorem and  $H^\bullet(X,P_\chi*P_\chi^\vee) = H^\bullet(X,P_\chi) \otimes_\Lambda 
H^\bullet(X,P_\chi^\vee)$ imply  $H^m(X,P_\chi)\neq 0$ for some $m \geq \dim(A)$. But $H^\bullet(X,P_\chi) = \bigoplus_{i,j\in \Z} H^j(B, {}^pR^iq_*(P)_\chi)$, using the decomposition theorem. Since $H^j(B, {}^pR^iq_*(P)_\chi)=0$ holds for $j\neq 0$  and  generic $\chi$ (generic vanishing theorem), ${}^pR^mq_*(P) \neq 0$ follows for some 
$m \geq \dim(A)$ and hence  $P = q^*(\tilde P)$ for some complex $\tilde P$ on $B$ by [BBD, prop.4.2.5]. Because $P$ is clean, this gives a contradiction unless $A=0$. So, for $A\neq 0$ this shows  $m_\nu < \dim(A)$. The lower inequalities now follow from the hard Lefschetz theorem.
\qed

\bigskip\noindent

\bigskip\noindent

\goodbreak
\bigskip\noindent
{\bf References}

%\bigskip\noindent
%[Atlas]  Conway J., Curtis R., Norton S., Parker R., Wilson R., {\it ATLAS of Finite Group Representations}, Oxford University Press (2003)

\medskip\noindent
[BBD] Belinson A., Bernstein J., Deligne P., {\it Faiscaux pervers}, Asterisque 100 (1982) 

\medskip\medskip\noindent
[BK] B\"ockle G., Khare C., {\it Mod $l$ representations of arithmetic fundamental groups. II.
A conjecture of A.J. de Jong}, Compos. Math. 142 (2006)

%\bigskip\noindent
%[C] Carter R.W., {\it Finite groups of Lie type}, John Wiley (1993)

\medskip\noindent
[D] Deligne P., {\it Finitude de l'extension de $\Q$ engendree par les traces de Frobenius, en characterisque finie}, Moscow Mathematical Journal 12 (2012)

\medskip\noindent
[D2] Deligne P., {\it Categories tannakiennes}, in: The Grothendieck Festschrift, Volume II, Birkh\"auser (1990) 

\medskip\noindent
[D3] Deligne P., {\it La Conjecture de Weil II}, Inst. Hautes Etudes Sci., Publ. Math. 52 (1980)

\medskip\noindent
[DM] Deligne P., Milne J., {\it Tannakian Categories}, in: Lecture Notes in Mathematics 900,
Springer Verlag (1982)

\medskip\noindent
[Dr] Drinfeld V., {\it On a conjecture of Kashiwara}, Math. Res. Lett. 8 (2001), 713 - 728

\medskip\noindent
[Dr2] Drinfeld V., {\it On a conjecture of Deligne}, Moscow Mathematical Journal 12 (2012)
%arXiv:1007.4004v5 (2011)

\medskip\noindent
[G]  Gaitsgory D., {\it On de Jong's conjecture}, Israel J. Math. 157, no.1 (2007), 155 - 191

\medskip\noindent
[KL] Katz N.M., Laumon G., {\it Transformation de Fourier et majoration de sommes
exponentielles}, Pub. math. de l'I.H.E.S., tome 62 (1985), 145 - 202 

\medskip\noindent
[KrW] Kr\"amer T., Weissauer R., {\it Vanishing theorems for constructible sheaves on abelian varieties},
J. Algebraic Geom. 24 (2015), 531-568 

\medskip\noindent
[KrW2] Kr\"amer T., Weissauer R., 
{\it On the Tannaka group attached to the Theta divisor of a generic principally 
polarized abelian variety}, arXiv:1309.3754 [math.AG] (to appear in 
Mathematische Zeitschrift)
 	
\medskip\noindent
[KrW3] Kr\"amer T., Weissauer R., {\it  Semisimple Super Tannakian Categories with a Small Tensor Generator}, Pacific Journal of Math. vol. 276 (2015), No. 1, 229 - 248

\medskip\noindent
[M] Matsushima Y., Espaces homogenes de Stein des groupes de Lie  complexes,
Nagoya Math. 18, 153 - 164 (1961)

\medskip\noindent
[GIT] Mumford D., Fogarty J., Kirwan F., {\it Geometric Invariant theory}, Third enlarged Edition, 
Erg. Math. Grenzg. 34, Springer (1991)

\medskip\noindent
[KW] Kiehl R., Weissauer R., {\it Weil conjectures, perverse sheaves and $l$-adic Fourier transform}, Erg. Math. Grenzg. (3. Folge) 42, Springer (2000)

\medskip\noindent
[L] Lafforgue L., Chtoucas de Drinfeld et correspondance de Langlands, Invent. Math. 147 (2002), no.1, 1 - 241  

\medskip\noindent
[S] Serre J.P., Linear Representations of Finite Groups, Springer (1977)

\medskip\noindent
[S2] Serre J.P., {\it Algebraic Groups and Class Fields}, Graduate Texts in Mathematics 117, Springer (1988)

\medskip\noindent
[VS] Viale L.B., Srinivas V., {\it Albanese and Picard 1-motives},  Memoire SMF 87, Paris (2001)

\medskip\noindent
[BN] Weissauer R., {\it Brill-Noether sheaves}, arXiv:math/0610923v4 (2006)

\medskip\noindent
[BN2] Weissauer R., {\it On the rigidity of BN-sheaves}, arXiv:1204.1929 (2012)  

\medskip\noindent
[W] Weissauer R., {\it  Degenerate Perverse Sheaves on Abelian Varieties}, 
arXiv:1204.2247 (2012)

\medskip\noindent
[W2] Weissauer R., {\it Vanishing theorems for abelian varieties over finite fields}, 
arXiv:1407.0844 (2014)

\end{document}